\documentclass[11pt]{amsart}
\usepackage{amsmath,amscd}
\usepackage{amssymb,xypic,graphicx}
\usepackage{amsthm,verbatim}

\usepackage[hidelinks]{hyperref}
\usepackage{tikz}
\usepackage{tikz-cd}
\usetikzlibrary{matrix,arrows,decorations.pathmorphing}

\newcommand*{\rom}[1]{\expandafter\@slowromancap\romannumeral #1@}
\newcommand{\tr}{{\sf tr}}

\newcommand{\fkg}{{\mathfrak{{g}}}}
\newcommand{\fks}{{\mathfrak{{s}}}}

\newcommand{\dec}{{\sf dec}}

\newcommand{\fkl}{{\mathfrak{{l}}}}
\newcommand{\db}{\overline{\partial}}

\newcommand{\LS}{{{\mathcal{L}_\infty{\mathcal S}}}}

\newcommand{\gr}{\operatorname{gr}}

\newcommand{\OO}{\mathcal{O}}
\newcommand{\Sym}{\operatorname{Sym}}

\newcommand{\II}{{\mathcal I}}

\newcommand{\Sh}{\operatorname{Sh}}

\newcommand{\CC}{{\mathcal{C}}}
\newcommand{\WW}{{\mathcal{W}}}

\newcommand{\supp}{\operatorname{supp}}

\renewcommand{\ker}{\operatorname{ker}}

\numberwithin{equation}{subsection}

\newtheorem{thm}{Theorem}[subsection]
\newtheorem{prop}[thm]{Proposition}
\newtheorem{lem}[thm]{Lemma}
\newtheorem{cor}[thm]{Corollary}
{  \theoremstyle{definition}
\newtheorem{defi}[thm]{Definition}
\newtheorem{ex}[thm]{Example}

\newtheorem{rem}[thm]{Remark}

\newtheorem{cons}[thm]{Construction}
\newtheorem{cond}[thm]{Condition}
}

\newcommand{\Pf}{\noindent {\it Proof}}
\newcommand{\id}{\operatorname{id}}
\newcommand{\Lie}{\operatorname{Lie}}

\newcommand{\ra}{\rightarrow}
\renewcommand{\AA}{{\mathcal A}}
\newcommand{\fM}{{\mathfrak M}}

\newcommand{\FF}{\mathcal{F}}
\newcommand{\EE}{{\mathcal{E}}}

\newcommand{\UU}{{\mathcal U}}

\newcommand{\HH}{{\mathcal H}}
\newcommand{\PP}{{\mathcal P}}
\newcommand{\QQ}{{\mathcal Q}}
\newcommand{\VV}{{\mathcal V}}
\newcommand{\cS}{{\mathcal S}}

\newcommand{\NL}{{\mathcal{NL}_\infty}}

\newcommand{\Om}{\Omega}

\newcommand{\Hom}{\operatorname{Hom}}

\newcommand{\Ext}{\operatorname{Ext}}
\newcommand{\End}{\operatorname{End}}
\newcommand{\EEnd}{{\mathcal E}{\sf nd}}

\newcommand{\Aut}{\operatorname{Aut}}

\newcommand{\C}{{\Bbb C}}

\newcommand{\R}{{\Bbb R}}

\newcommand{\ot}{\otimes}

\newcommand{\ad}{\operatorname{ad}}

\newcommand{\ed}{\qed\vspace{3mm}}

\newcommand{\Ad}{\operatorname{Ad}}

\newcommand{\MC}{\sf MC}
\newcommand{\mc}{\sf mc}

\newcommand{\at}{\makeatletter @\makeatother}

\title{Homotopy L-infinity spaces}
\author{Junwu Tu}
\thanks{Mathematics Department, University of Missouri, Columbia, MO, 65211, USA, {\em e-mail:}{\texttt{ tuju\at missouri.edu}}}

\begin{document}
\begin{abstract}
In this paper, we introduce a new class of structured spaces which is locally modeled by Costello's L-infinity spaces. This provides an alternative approach to study the derived geometric structures in the algebraic, analytic, or smooth category. Our main result asserts that on a compact complex manifold, the moduli space of simple holomorphic structures on a complex vector bundle with a fixed determinant bundle, admits a natural analytic homotopy L-infinity enhancement.
\end{abstract}

\maketitle

\section{Introduction}\label{intro-sec}

\subsection{Backgrounds} In their seminal papers~\cite{DT} and~\cite{Tho}, Donaldson and Thomas introduced the now named Donaldson-Thomas invariants by {\sl virtual} counting of compact moduli spaces of stable coherent sheaves on Calabi-Yau 3-folds, using the construction of virtual fundamental cycles (see~\cite{BF} and~\cite{LT}). In the past few years, there have been some spectacular work towards deformation quantization or categorification of the DT theory. This line of research begins with the paper~\cite{PTVV} of Pantev-To\" en-Vaqui\' e-Vezzosi where the important notion of shifted symplectic structures was introduced, which is further based on the foundational work of To\" en-Vezzosi on derived algebraic geometry~\cite{TV},~\cite{TV2}. Later on, Joyce's school produced several papers~\cite{BBBJ},~\cite{BBJ},~\cite{BBDJS},\cite{BJM},~\cite{Joyce}, which reveal extremely rich geometric structures on the moduli spaces. In particular, it follows from their work that it is possible to construct a perverse sheaf on the moduli space (under suitable orientation data) whose Euler characteristic gives the DT invariants, which categorifies the DT theory on the level of chain complexes. The construction of such a perverse sheaf was also independently studied by Kiem-Li~\cite{KL}. Recently in his ICM paper~\cite{Toen}, To\" en proposed, in collaboration with Calaque, Pantev,Vaqui\' e, and Vezzosi,  a general approach to study the deformation quantization problem for shifted symplectic spaces (which include the DT moduli space as the $-1$-shifted special case).

\subsection{Structured spaces.} A guiding principle for all these developments is that while the moduli spaces are often rather singular, there are natural homological enhancements over the underlying geometric spaces so that these structured spaces behave like smooth varieties. This allows us to study the geometry of moduli spaces. For example, we may obtain virtual fundamental classes, or consider vector fields/differential forms, and even study symplectic geometry/deformation quantization, as Pantev-To\" en-Vaqui\' e-Vezzosi suggested. This is usually referred to as the ``hidden smoothness" principle, first noted by Kontsevich in~\cite{Kont}. The main idea, in the case of moduli space of coherent sheaves, is that the tangent space at a point $\EE$ in the moduli should be a complex with cohomology groups equal to
\begin{equation}~\label{hs-eq}
\Ext^*(\EE,\EE)[1]=T_\EE,\end{equation}
whose ``dimension" equals to $\sum_i (-1)^{i+1} \dim \Ext^i(\EE,\EE)=-\chi(\EE^\vee\ot \EE)$ which is a topological invariant of $\EE$, hence a locally constant function on the moduli. This gives a hint on the existence of certain smoothness in the sense of homological algebra. There are, at least, four different approaches to set up the foundations for studying this hidden homological structure:
\begin{itemize}
\item[--] (Ciocan-Fontanine, Kapranov): Differential graded schemes~\cite{CK};

\item[--] (Joyce, Spivak): Derived differential geometry~\cite{Joyce2},~\cite{Spivak};

\item[--] (Lurie): Derived algebraic geometry~\cite{Lurie};

\item[--] (To\" en, Vezzosi): Homotopical algebraic geometry~\cite{TV},~\cite{TV2}.
\end{itemize}
Among them, Ciocan-Fontanine and Kapranov's theory of dg schemes is certainly the most accessible one. In order to avoid using an ambient smooth space, it is necessary to develop a more flexible theory as done by Lurie and To\" en-Vezzosi where homotopy sheaf theory is used to define structured spaces that are locally modeled by dg schemes. Finally, Spivak and Joyce's theory deals with smooth manifolds (possibly with boundaries/corners), hence in principle can be applied to deal with more general moduli spaces, such as various types of moduli spaces of J-holomorphic curves in symplectic geometry. A common feature of all four approaches is that the local models used are always commutative algebras.

\subsection{Costello-Kapranov's resolution of the structure sheaf} In~\cite{Costello}, Costello proposed to use Lie structures to study derived geometry. In a recent work~\cite{GG}, Grady and Gwilliam further developed Costello's ideas. The local models in this approach are called ``$L_\infty$ spaces". This $L_\infty$ approach may be considered as a Koszul dual description of the commutative approaches mentioned above. 

Let us explain the main ideas more precisely. By the hidden smoothness principle described above (see~\ref{hs-eq}), we already know what should the derived tangent bundle be. Indeed, in the case of moduli spaces of sheaves, the ``tangent bundle" ought to be a complex whose cohomology sheaves, after restriction on a point $\EE$ in the moduli, computes the shifted self $\Ext$-groups of $\EE$. So we ask the following

\medskip
{\bf Question.} {\sl How could we recover the structure sheaf of a space from its tangent bundle?}

\medskip
Our goal is to answer this question on the moduli space of coherent sheaves. But, in order to get a hint of how to proceed, it is helpful to first consider the case when the underlying space is a smooth space $M$ with the usual tangent bundle $T_M$. It turns out that even in this simplest form, the question above is highly non-trivial. An elegant answer is given in a beautiful paper of Kapranov~\cite{Kap}.

Let us explain Kapranov's solution. First, observe that the infinite jet bundle $J_M$ has a natural flat connection $D$ called the Grothendieck connection from which we may recover the structure sheaf by taking the kernel of the connection map
\[ J_M \stackrel{D}{\ra} \Om_M\ot J_M.\]
To answer the question, it remains to find a relationship between $J_M$ and $T_M$. For this, note that the bundle $J_M$ has a 
natural decreasing filtration (by jets vanishing in a finite order) which is complete and with associated graded bundles the symmetric powers $S^k\Om_M$. Thus we conclude that
\[ \prod_{k=0}^\infty \gr^k(J_M)\cong \hat{S} \Om_M.\]
This filtration, most of the time, does not split, which makes it difficult to recover the algebra bundle $J_M$ from $T_M$. The main observation of Kapranov is that after taking appropriate resolution of both sides, the induced filtration always splits. Via the above isomorphism, the differential on the resolution of $J_M$ induces a ``deformation" of the differential on the resolution of $\hat{S}\Om_M$. Such a differential is by definition a {\sl $L_\infty$ structure} on $T_M[-1]$. The characteristic property of this $L_\infty$ structure is that its Chevalley-Eilenberg algebra recovers the algebra bundle $J_M$.

Formalizing the structures appeared in the discussion above, Costello introduced the notion of $L_\infty$ spaces which plays an important role in his geometric study of the Witten genus.

\subsection{Homotopy L-infinity enhancements} In the case of moduli space of sheaves, we need to give a modified version (see Definition~\ref{diff-lie-defi}) of Costello's definition of $L_\infty$ spaces. On the mathematical level, the main difference is that we allow the $L_\infty$ structure to have a curvature term, even after modulo  differential forms of positive degrees. This modification is crucial in the setting of moduli spaces of sheaves since the dimension of the $\Ext$-groups may jump when considered as a function on the moduli spaces, known as the semicontinuity theorem.
On the expository level, the definition we gave emphasized the formal geometry point of view. We refer to Remark~\ref{diff-rem} for a detailed discussion of the differences. 

Roughly speaking, a $L_\infty$ space is a triple $(V,\fkg,D)$ where $V$ is a smooth space, $\fkg$ is a (possibly curved) $L_\infty$ algebra bundle over $V$ whose Chevalley-Eilenberg algebra plays the role of the ``jet bundle", and $D$ is a flat connection 
\[D: C^*\fkg \ra \Om_V \ot C^*\fkg\]
on the Chevalley-Eilenberg algebra bundle satisfying certain compatibility conditions. We refer to Definition~\ref{diff-lie-defi} for a precise version. Furthermore, we say that a given $L_\infty$ space $(V,\fkg,D)$ is an enhancement of a possibly singular space $U$, if there exists a closed embedding $\fks: U\ra V$ such that the natural map $C^*(V,\fkg,D) \ra \fks_*\OO_U$ is a chain map, and that it induces an isomorphism 
\[ \underline{H}^0\big(C^*(V,\fkg,D)\big) \cong \fks_*\OO_U.\]
Here the notation $C^*(V,\fkg,D)$ denotes the de Rham complex of the Chevalley-Eilenberg algebra $C^*\fkg$ endowed with the flat connection $D$. See Definition~\ref{lie-enh-defi} for a precise definition of $L_\infty$ enhancements. 

It turns out that locally (in the analytic topology) on the moduli space of simple holomorphic bundles, we can construct $L_\infty$ enhancements. This essentially follows from Kuranishi theory, interpreted in the language of $L_\infty$ algebras. However, these local pieces do not paste in the strict sense, i.e. with isomorphisms on double intersections that satisfy the cocycle condition on triple intersections. Since $L_\infty$ structure admits homotopy theory, it is only natural to seek for a homotopy version of the cocycle condition. Namely, there are transition morphisms on double intersections, which satisfy the cocycle condition up to chosen homotopies on triple intersections, which should still be bounded by chosen $2$-homotopies on quadruple intersections, and so forth on all multiple intersections. This leads us to the notion of a {\sl homotopy $L_\infty$ enhancement}, formulated precisely in Definition~\ref{h-lie-enh-defi}. An important technical point in this definition is that we use {\sl hypercoverings} instead of C\v ech coverings as relevant convergence properties can only be proved in a neighborhood of every point in a multiple intersection. The main result of the paper is the following
\medskip
\begin{thm}~\label{main-thm}
Let $X$ be a compact complex manifold, and $E$ be a smooth complex vector bundle on $X$. Let $M$ be a (locally) universal moduli space of simple holomorphic structures on $E$ with fixed determinant bundle. Then $M$ admits a homotopy $L_\infty$ enhancement.
\end{thm}

The existence of homotopy gluing data immediately implies that the Chevalley-Eilenberg cohomology of the locally defined $L_\infty$ enhancements glue into a global object on $M$. Thus we obtain 

\medskip
\begin{cor}~\label{main-cor}
Let $X$ and $M$ be as in the previous theorem. Then there exists a canonical non-positively graded sheaf of $\OO_M$-algebras $\OO_M^\bullet=\oplus_{k=0}^{\infty} \OO_M^{-k}$ such that $\OO_M^0=\OO_M$, and with each graded component $\OO_M^{-k}$ coherent over $\OO_M$. Furthermore, if $X$ is a Calabi-Yau 3-fold, then the whole algebra $\OO_M^\bullet$ is coherent as a $\OO_M$-module.
\end{cor}

\begin{rem}
Assuming that $X$ is an algebraic Calabi-Yau variety,  Joyce and Song~\cite{Joy-Song} introduced a trick to use Seidel-Thomas twist functor to show that a bounded family of coherent sheaves is isomorphic as an analytic space to a moduli space of vector bundles. Because of this, Theorem~\ref{main-thm} and Corollary~\ref{main-cor} also apply for any bounded moduli space of coherent sheaves.

Assume that $X$ is a Calabi-Yau $3$-fold, and furthermore that the moduli space $M$ is compact. By the previous paragraph, this compactness assumption is fine since we may allow coherent sheaves in $M$. Following Ciocan-Fontanine, Kapranov~\cite{CK2}, and Kontsevich~\cite{Kont}, one may define the (coherent sheaf) K-theoretic virtual fundamental class to be
\[ [\OO_M^\bullet] \in K(M).\]
One can show, by a local computation (using Proposition~\ref{D-cohomology-prop}), that the class $[\OO_M^\bullet]$ has zero dimensional support, and hence the support map yields a class
\[ \supp [\OO_M^\bullet] \in A_0(M).\]
This gives an alternative construction of the virtual fundamental class of $M$. 

The deformation invariance of this virtual cycle construction would follow if one could carry out the construction of a homotopy $L_\infty$ enhancement on the relative moduli space over a one dimensional domain. This is left for future investigation.
\end{rem}

\subsection{Results of the paper} Finally we give a detailed section-wise summary of the paper.

In Section~\ref{loc-enh-sec}, we derive the Kuranishi theory from the viewpoint of deformation theory of normed $L_\infty$ algebras. The main scheme of the section is standard, and indeed we follow Fukaya's paper~\cite{Fukaya}. The section is based on the convergence properties proved in the Appendix~\ref{ht-sec}. Local existence of $L_\infty$ enhancements on moduli spaces follows essentially from this interpretation of Kuranishi theory. We also give some applications of this more algebraic point of view to give local descriptions of various moduli spaces. For example, we prove that
\begin{itemize}
\item[--] on $K3$ surfaces, the moduli spaces are smooth, recovering Inaba's result~\cite{In};
\item[--] on Calabi-Yau $3$-folds, the moduli spaces are locally critical loci, a result originally due to Joyce-Song~\cite{Joy-Song}.
\item[--] on CY $4$-folds, the local models are given by zero locus of a holomorphic map
\[ \kappa : V \ra E,\]
where $V$ is an open subset of a complex vector space, and $E$ is another complex vector space endowed with a symmetric nondegenerate bilinear form $\langle-,-\rangle$. The map $\kappa$ satisfies the condition $\langle \kappa, \kappa \rangle=0$. This result is due to Borisov-Joyce~\cite{BJ}. See also the recent work of Cao-Leung~\cite{CL} in a different setting.
\end{itemize}

In Section~\ref{lie-space-sec}, we introduce the notion of $L_\infty$ spaces by putting a flat connection (analogous to the Grothendieck connection) on the Chevalley-Eilenberg algebra of a $L_\infty$ algebra bundle. Then we study the Chevalley-Eilenberg cohomology of $L_\infty$ spaces, which in particular clarifies the relationship between $L_\infty$ spaces and dg-schemes/derived schemes.

In Section~\ref{functor-sec}, we construct a ``functor"
\[ \cS: \NL \ra \LS\]
from the simplicial category of normed $L_\infty$ algebras to the simplicial category of $L_\infty$ spaces. It sort of {\sl propagates} the given $L_\infty$ algebra $g$ to nearby fibers in a neighborhood of the origin in the degree one part of $g$. Hence we refer to it as the propagation ``functor".  A difficulty of the simplicial category $\LS$ is that it is not clear if its Hom simplicial sets are Kan. But in order to construct higher homotopies in Theorem~\ref{main-thm}, it is crucial to have certain homotopy lifting property. The way out of this is to prove that the Hom simplicial sets in $\NL$ are Kan. This is done in the Appendix~\ref{kan-sec}. In the end of the section, we formulate a precise definition of homotopy $L_\infty$ spaces, using Lurie's language of $\infty$-categories (see~\cite{Lurie}). 

In Section~\ref{moduli-sec}, we prove the main Theorem~\ref{main-thm} and Corollary~\ref{main-cor}. All the higher homotopies required in the proof are in the image of $\cS$. This enables us to use the Kan property of $\NL$ to deduce inductively the existence of higher homotopies on the next level of intersections. The construction of this section works in the exact same way for elliptic complexes appearing in other gauge theories.

Appendix~\ref{ht-sec} deals with homological transferring of $A_\infty$ or $L_\infty$ algebras in the normed setting. We give a detailed proof of relevant convergences for infinity algebras, homomorphisms, and homotopies, which are essential for applications in this paper. The explicit transferring formulas written down by Markl~\cite{Markl} play an essential role.

Appendix~\ref{kan-sec} is devoted to the proof of Proposition~\ref{kan-prop} which asserts the Kan property of the Hom sets in $\NL$. In the case of $L_\infty$ algebras without norms, this is a consequence of the Kan property of nerves of pronilpotent $L_\infty$ algebras, proved by Hinich~\cite{Hinich} in the differential graded setting, and by Getzler~\cite{Getzler} in general. The proof in the normed setup follows Getzler's paper.

\medskip
{\bf Acknowledgments.} I am grateful to Alexander Polishchuk for numerous helpful discussions and his guidance during my postdoc years at University of Oregon. It is only after our joint work~\cite{PT} in the case of Jacobians did I realize that the algebraic structures on Ext-groups can give an alternative description of the derived structures on moduli spaces. I also thank Weiyong He for his constant encouragement to learn some analysis, and for tolerating me to ask quite a few naive questions about heat kernels.

\section{Deformation theory via infinity algebras}\label{loc-enh-sec}

It is a well established result that the \emph{formal} deformation theory of a holomorphic vector bundle is governed by a differential graded (Lie) algebra, or its homotopy version, $A_\infty$ ($L_\infty$) algebras. On the other hand, Kuranishi theory~\cite{Miy} gives the actual deformation theory, as opposed to the formal one. In this section, we use certain convergence properties on $A_\infty$/$L_\infty$ algebras to recover the Kuranishi theory from the algebraic deformation theory. The section is based on Appendix~\ref{ht-sec} where detailed proofs of the analysis involved are presented. 

\subsection{Terminologies} Let $X$ be a smooth projective variety over $\C$. We consider $X$ as a complex analytic manifold. Let $E$ be a fixed smooth vector bundle over $X$. For a vector bundle $V$ on $X$, we shall denote by $C^\infty(V)$ the space of its smooth sections.

Recall a semiconnection on $E$ is a first order differential operator 
\[\db: C^\infty(E) \ra C^\infty(\Omega_X^{0,1}\ot E)\]
such that $\db(fs)=\db(f) s +f\db(s)$ for all $f\in C^\infty_X$ and $s\in C^\infty(E)$. A semiconnection $\db$ extends to a unique operator (which we still denote by $\db$) acting on $\Omega_X^{0,*}\ot E$ by requiring the following equation
\[ \db(\alpha s)=\db(\alpha) s+(-1)^k\alpha \db(s)\]
holds for any $k$, and any homogeneous Dolbeault form $\alpha\in \Omega_X^{0,k}$. A semiconnection $\db$ on $E$ is called integrable if $\db\circ\db=0$. A holomorphic structure on $E$ is an integrable semiconnection.

Let $\EE$ be a holomorphic structure on $E$. By definition, $\EE$ corresponds to an integrable semiconnection $\db_\EE$. The deformation theory of $\EE$ is governed by the differential graded (Lie) algebra $\Omega_X^{0,*}(\EEnd(\EE))$ whose underlying vector space is $\Omega_X^{0,*}(\EEnd(E))$, and is endowed with the differential $\ad(\db_\EE)$. Indeed, given a \emph{Maurer-Cartan} element $B\in \Omega_X^{0,1}(\EEnd(\EE))$, i.e. $B$ satisfies the equation
\[ [\db_\EE, B] + B\circ B = 0,\]
we get another complex structure on $E$ defined by the operator $\db_\EE+B$. The integrability of $\db_\EE+B$ follows from the Maurer-Cartan equation above.

We shall consider families of holomorphic structures on $E$, deforming $\EE$, parametrized by local analytic spaces. Recall that a local analytic space is a locally ringed space $(U,\OO_U)$ such that 
\begin{itemize}
\item[(1.)] $U=V\cap Z(f_1,\cdots, f_r)$ for some open subset $V \subset \C^N$, and $Z(f_1,\cdots, f_r)$ is the zero locus of $f_1,\cdots, f_r\in \OO_V$;
\item[(2.)] $\OO_U=\OO_V/\langle f_1,\cdots, f_r \rangle$.
\end{itemize}
A pointed local analytic space is the pair $(U,\OO_U)$ together with a point $o\in U$. 

We also need a few terminologies in~\cite{Miy}. For a fixed integer $l>\dim_\R X+1$, we denote the Sobolev completion of order $l$ by a subscript $l$. We set $$\FF:\Omega_X^{0,1}(\EEnd(\EE))_l \ra \Omega_X^{0,2}(\EEnd(\EE))_{l-1}$$ to be the analytic map between Banach spaces defined by 
\[\FF(\alpha)=\ad(\db_\EE)(\alpha)+\alpha\circ\alpha~\footnote{We required that $l>\dim_\R X +1$ to ensure the product morphism is bounded.}.\]
A deformation of $\EE$ parametrized by a pointed local analytic space $(U,\OO_U,o)$ is a morphism of pointed ringed spaces from $(U,\OO_U,o)$ to the germ of Banach analytic space $(\FF^{-1}(0),0)$, which is of class $C^\infty$ in the $X$-direction.

\subsection{Kuranishi family via A-infinity algebras} Homological transferring applied to the differential graded algebra $\Omega_X^{0,*}(\EEnd(\EE))$ yields an $A_\infty$ structure on its cohomology $\Ext^*(\EE,\EE)$, together with an $A_\infty$ homomorphism
\[ I: \Ext^*(\EE,\EE)\ra \Omega_X^{0,*}(\EEnd(\EE)).\]
For an element $b\in \Ext^1(\EE,\EE)$, we define two infinite series
\begin{align*}
 \kappa(b)&:= \sum_j m_j(b^j),\\
 I_*(b)&:=\sum_j I_j(b^j).
\end{align*}
Lemma~\ref{norm-m-lem} and Lemma~\ref{norm-i-lem} imply that the radius of convergence of $\kappa$ and $I_*$ are both positive. By degree reasons, provided convergence, $\kappa(b)$ lies in $\Ext^2(\EE,\EE)$, and $I_*(b)$ lies in $\Omega_{X}^{0,1}(\EEnd(\EE))$~\footnote{Smoothness of $I_*(b)$ is proved in Lemma~\ref{smooth-lem}.}.

The zeros of $\kappa(b)$, provided convergence, are called Maurer-Cartan elements of the $A_\infty$ algebra $\Ext^*(\EE,\EE)$. Let $V\subset \Ext^1(\EE,\EE)$ be a small neighborhood of the origin (in the sense that $V$ is contained in a disc of radius $\frac{1}{C}$ where $C$ is as in Lemma~\ref{norm-m-lem}) in the vector space $\Ext^1(\EE,\EE)$. The analytic map
\[ \kappa: V \to \Ext^2(\EE,\EE)\]
cuts out a pointed local analytic space $(U,\OO_U,o):=(\kappa^{-1}(0), 0)$. 

\medskip

\begin{cons}~\label{a-kura-cons}
We next construct a deformation of $\EE$ parametrized by $(\kappa^{-1}(0),0)$. Consider the following diagram of maps
\[\begin{CD} 
V @> I_* >>  \Om_X^{0,1}(\EEnd(\EE))_l \\
@VV\kappa V                 @VV\FF V \\
\Ext^2(\EE,\EE)       @.        \Om_X^{0,2}(\EEnd(\EE))_{l-1}.
\end{CD}\]
To show that the above diagram induces a morphism $\kappa^{-1}(0)$ to $\FF^{-1}(0)$, it suffices to show that for every bounded linear functional $\lambda$ on $\Om_X^{0,2}(\EEnd(\EE))$, the composition
\[ \lambda\circ \FF \circ I_*\]
is analytic, and furthermore lies inside the ideal defining $\kappa^{-1}(0)$. This follows from the following calculation:
\begin{align}\label{regular-eq}
(\FF \circ I_*)(b)= \ad(\db_\EE)(I_*(b))+I_*(b)\circ I_*(b) \nonumber\\
 = \sum_{j_1,j_2\geq 0} I_{j_1+j_2+1}(b^{j_1}\ot \kappa(b) \ot b^{j_2}).
\end{align}
Here the second equality uses the fact that $I$ is an $A_\infty$ homomorphism. Applying a bounded linear functional $\lambda$ to $\sum_{j_1,j_2\geq 0} I_{j_1+j_2+1}(b^{j_1}\ot \kappa(b) \ot b^{j_2})$ produces an analytic function in the defining ideal of $\kappa^{-1}(0)$. Note that here the convergence follows from the fact that $\lambda$ is bounded and estimates in Lemma~\ref{norm-i-lem}.

It remains to show that this analytic family is smooth in the $X$-direction, i.e. $I_*(b)$ is smooth for all $b\in V$.  This is proved in Lemma~\ref{smooth-lem}.\ed
\end{cons}

The following Theorem was first proved by Miyajima~\cite{Miy} in a different setup, and was rediscovered by Fukaya~\cite[Sections 1,2]{Fukaya} in the current setting. We include a proof here to add into Fukaya's proof a more detailed account of relevant convergences, using Lemmas~\ref{norm-m-lem},~\ref{norm-i-lem},~\ref{norm-p-lem}, and Corollary~\ref{bound-h-cor}.
\medskip

\begin{thm}~\label{a-kura-thm}
After appropriately shrinking the open subset $V$, the local deformation constructed in Construction~\ref{a-kura-cons} is a versal family. If $\EE$ is furthermore simple, then it is a universal family.
\end{thm}

\Pf. Let $(U',o)$ be a pointed local analytic space. Assume that
\[ B: U' \ra \Omega^{0,1}_X(\End(\EE))_l\]
is an analytic family of deformation of $\EE$. Thus the point $o$ gets mapped to $0$. As shown in the Appendix~\ref{a-pert-subsec}, there is an $A_\infty$ homomorphism
\[ P: \Omega_X^{0,*}(\EEnd(\EE))_l \ra \Ext^*(\EE,\EE)\]
which also satisfies certain boundedness condition discussed in Lemma~\ref{norm-p-lem}. The boundedness condition implies the series
\[ P_* (\alpha):= \sum_{k\geq 1} P_k(\alpha^k)\]
is absolutely convergent if $||\alpha||$ small enough. Since $B(o)=0$ which has norm zero, the composition
\[ U'\stackrel{ B}{\longrightarrow} \Omega^{0,1}_X(\End(\EE))_l \stackrel{P_*}{\longrightarrow} \Ext^1(\EE,\EE)\]
is well-defined by appropriately shrinking $U'$. Since the push-forward of a Maurer-Cartan element is again Maurer-Cartan element, the composition
\[ U'\stackrel{ B}{\longrightarrow} \Omega^{0,1}_X(\End(\EE))_l \stackrel{P_*}{\longrightarrow} \Ext^1(\EE,\EE) \stackrel{\kappa}{\longrightarrow}\Ext^2(\EE,\EE)\]
is equal to zero, inducing a morphism
\[ \phi: U'\ra \kappa^{-1}(0).\]

Next we prove that the pull-back family via $\phi$ is isomorphic to the family $B$. For this, we need to show that for any $x\in U'$ sufficiently close to $o$, the two Maurer-Cartan elements $B(x)$ and $I_*P_*B(x)$ are gauge-equivalent. We use the homotopy $H$ between $\id$ and $I\circ P$, defined in the Appendix~\ref{a-pert-subsec}. Indeed, by Corollary~\ref{bound-h-cor}, the infinite series
\[ H_*( B(x) ) = \sum_{k\geq 1} H_k ( B(x)^k )\]
is absolutely convergent in the $C^\infty$ topology of $\Om_{[0,1]}^*\big(\Omega^{0,*}_X(\End(\EE))_{l-1}\big)$. Being push-forward of a Maurer-Cartan element, $H_*( B(x) )$ is also Maurer-Cartan. Writing in components we have
\[ H_*( B(x) ) = B(t,x) +C(t,x) dt\]
for some
\begin{align*}
B(t,x)&\in C^\infty\big([0,1], \Omega^{0,1}_X(\End(\EE))_{l-1}\big),\\
C(t,x)&\in C^\infty\big([0,1], \Omega^{0,0}_X(\End(\EE))_{l-1}\big).
\end{align*}
To find the gauge equivalence between $B(0,x)=B(x)$ and $B(t,x)$, one solves the following ordinary differential equation
\[ \frac{dg(t,x)}{dt}= g(t,x) \cdot C(t,x), \mbox{\;\; with \;\;} g(0,x)=\id.\]
For $x$ sufficiently close to $o$, due to the boundedness of $H$, the norm of $C(t,x)$ can be made small enough so that the above ODE is solvable for all $t\in [0,1]$. The endomorphism $g(1,x)$ gives us the desired gauge equivalence between $B(x)$ and $I_*P_* B(x)$.

Assume that the holomorphic bundle $\EE$ is furthermore simple, i.e. $\End(\EE)=\C$. We need to show that if $\phi': U'\ra \kappa^{-1}(0)$ is another morphism such that for any $x\in U'$ the Maurer-Cartan elements $I_*\phi'(x)$ and $I_*\phi(x)$ are gauge-equivalent, then $\phi'=\phi$. We may shrink the open subset $V\subset \Ext^1(\EE,\EE)$ small enough so that there exists an element $\eta\in \Om_X^{0,0}(\EEnd(\EE))$ such that
\[ I_*\phi'(x)= e^{\eta}* I_*\phi(x)\]
where $*$ on the right hand side denotes the gauge-group action. This implies that we have a Maurer-Cartan element
\[ \theta:=e^{t\eta}* I_*\phi(x) + \eta dt \in \Om_X^{0,*}(\EEnd(\EE))\ot \Om_{[0,1]}^*\]
such that $\theta(0)=I_*\phi(x)$ and $\theta(1)=I_*\phi'(x)$. Pushing forward $\theta$ via $P\ot \id$ (well-defined since $P$ is bounded, see Lemma~\ref{norm-p-lem}), we obtain a Maurer-Cartan element $(P\ot\id)_*\theta$ of the $A_\infty$ algebra $\Ext^*(\EE,\EE)\ot \Om_{[0,1]}^*$, such that
\[ (P\ot\id)_*\theta (0)= P_*I_*\phi(x)=\phi(x), \mbox{\;\; while \;\;} (P\ot\id)_*\theta (1)= P_*I_*\phi'(x)=\phi'(x).\]
But by simpleness, the degree zero part of the $A_\infty$ algebra is $\Ext^0(\EE,\EE)=\End(\EE)=\C$ consisting of multiples of the strict unit $\id_\EE$ which acts trivially on the set of Maurer-Cartan elements of $\Ext^*(\EE,\EE)$. Hence we conclude that $\phi=\phi'$.\qed

\subsection{Local descriptions of moduli spaces} The previous subsection illustrates that locally (in the analytic topology) deformations of $\EE$ is completely determined by the $A_\infty$ structure on $\Ext^*(\EE,\EE)$. In this subsection, we give some applications of this observation to give a local description of the moduli spaces when $X$ is compact and Calabi-Yau. Under these assumptions, the following algebraic properties are proved in~\cite{Poli}:
\begin{itemize}
\item The $A_\infty$ algebra $\Ext^*(\EE,\EE)$ is always strict unital. Namely, denote by $1\in \Ext^0(\EE,\EE)$ the identity morphism of $\EE$. Then we have
\begin{align*}
 m_2(1,a)&=(-1)^{|a|}m_2(a,1)=a, \mbox{\;\; and}\\ 
 m_k & (a_1,\cdots,1,\cdots,a_{k-1})=0, \;\;\forall k\geq 3.
\end{align*}
\item The Calabi-Yau assumption implies the $A_\infty$ algebra $\Ext^*(\EE,\EE)$ is also {\sl cyclic}. That is, it admits a non-degenerate graded symmetric pairing $\langle-,-\rangle$ such that
\[ \langle m_k(a_0,\ldots,a_{k-1}),a_k\rangle =(-1)^{k+\mid a_0\mid\sum_{l=1}^k\mid a_l\mid}\langle m_n(a_1,\ldots,a_k),a_0\rangle .\]
\end{itemize}

\medskip
These algebraic structures can be used to describe local models of the moduli spaces, which we now illustrate in low dimensions.

\medskip
(1.) {\bf $\dim X=2$.} We assume furthermore that the bundle $\EE$ is simple. Under this assumption, we claim that the moduli space is smooth near $\EE$. To prove this, it suffices to show that the Kuranishi map $\kappa\equiv 0$. By the simpleness assumption, we have
$\Ext^0(\EE,\EE)=1\cdot \C$, which is dual to $\Ext^2(\EE,\EE)$. Thus we calculate the pairing
\[\langle\kappa(b),1\rangle = \sum_{k\geq 2} \langle m_k(b^k), 1\rangle.\]
By cyclic symmetry and strict unital properties, the pairing $\langle m_k(b^k), 1\rangle$ vanishes for $k\geq 3$. The remaining term 
\[ \langle m_2(b,b),1\rangle =\pm \langle b, b\rangle=0,\]
since the pairing is graded-symmetric, and hence anti-symmetric on degree one elements. By non-degeneracy we conclude that $\kappa\equiv 0$.

\medskip
(2.) {\bf $\dim X=3$.} This is arguably the most interesting case. In this dimension, locally around $\EE$, the moduli space can be written as the critical locus of a holomorphic function 
\[ \Psi: V \ra \C,\]
where $V\subset \Ext^1(\EE,\EE)$ is a small enough open subset containing the origin. This result is proved in~\cite[Theorem 5.4 and 5.5]{Joy-Song}. In the following we give a simple proof of this result. Using the pairing, we have a natural isomorphism 
\[ \iota: \Ext^2(\EE,\EE) \cong \Ext^1(\EE,\EE)^\vee.\]
Pre-composing with the Kuranishi map gives an analytic morphism
\[ \iota\circ \kappa: V \ra \Ext^1(\EE,\EE)^\vee.\]
The point is that this map can be viewed as a holomorphic section of $\Om_V^1$ since the space $V$, being an open subset of the vector space $\Ext^1(\EE,\EE)$, its tangent bundle is $\OO_V\ot_\C \Ext^1(\EE,\EE)$. Furthermore, if we define a holomorphic  function $\Psi:V\ra \C$ by
\[ \Psi(b):=\sum_{k\geq 2} \frac{1}{(k+1)} \langle m_k (b^k),b \rangle,\]
it is straightforward to verify that $d\Psi= \iota\circ \kappa$. As one might expect, the potential function $\Psi$ defined above is simply the pull-back of the holomorphic Chern-Simons functional on $\Om_X^{0,1}(\EEnd(\EE))$ via the analytic map
\[ I_*: V \ra \Om_X^{0,1}(\EEnd(\EE)).\]
Finally, we note that in ~\cite{Joy-Song}, Joyce and Song also proved that both $V$ and $\Psi$ can be chosen so that it is invariant under the gauge action of
$$G^c\subset \Ext^0(\EE,\EE),$$ 
the complexification of a maximal compact subgroup $G\subset \Ext^0(\EE,\EE)$. In our setting, one can deduce the same result using the fact that the $A_\infty$ Chern-Simons functionals pull-back via cyclic quasi-isomorphisms. 

\medskip
(3.) {\bf $\dim X=4$.} Again locally around a bundle $\EE$, the moduli space is the zero locus of the Kuranishi map $\kappa: \Ext^1(\EE,\EE)\supset V \ra \Ext^2(\EE,\EE)$. In this case, the pairing on $\Ext^2(\EE,\EE)$ is symmetric and we have $\langle \kappa,\kappa\rangle\equiv 0$. To see this, we use the construction of Subsection~\ref{loc-fam-subsec} to get a family of $A_\infty$ algebra over $V$. Explicitly this curved differential graded bundle is given by
\begin{equation}~\label{dim4}
 T_V \stackrel{d\kappa}{\longrightarrow} \OO_V\ot_\C \Ext^2(\EE,\EE) \stackrel{\big(d\kappa\big)^\vee}{\longrightarrow} \Om_V,
\end{equation}
where we have used the non-degenerate pairing to get identifications
\[ \OO_V\ot_\C \Ext^2(\EE,\EE) \cong \big( \OO_V\ot_\C \Ext^2(\EE,\EE)\big)^\vee,\;\;\mbox{ and }\;\; \OO_V\ot_\C\Ext^3(\EE,\EE) \cong \Om_V.\]
Note that the composition in Equation~\ref{dim4} does not equal to zero since there is the curvature term $\kappa$. On the other hand, the $A_\infty$ identity $m_1m_0=0$ implies that
\[ (d\kappa)^\vee (\kappa)=0, \;\; \mbox{ or equivalently} \;\; \langle d\kappa, \kappa\rangle=0.\]
Since $\kappa(0)=0$, we have $\langle \kappa,\kappa \rangle\equiv 0$. From the identity $\langle d\kappa, \kappa\rangle=0$, we can also conclude that the moduli space is never transversal unless $\kappa=0$ in which case it is smooth. This is because if $\kappa\neq 0$, then there exists a tangent vector $v$ such that $d_v \kappa\neq 0$. The identity
\[ \langle d_v\kappa, \kappa\rangle=0\]
then implies that components of $\kappa$ do not form a regular sequence (as $d_v\kappa$ lies its kernel). 

Another simple observation is that if $\dim \Ext^2(\EE,\EE)=1$, then $\kappa=0$ because the pairing is symmetric on degree $2$ elements. The above local description was originally due to Borisov-Joyce~\cite{BJ}, which was also discussed by  Cao-Leung~\cite{CL} in the differential geometric setting.

\section{L-infinity spaces}~\label{lie-space-sec}

In this subsection, we introduce the notion of $L_\infty$ spaces which serve as local models in this Lie approach to derived geometry. We also study the Chevalley-Eilenberg cohomology of $L_\infty$ spaces through which the relationship between the Lie approach and the classical commutative approach is clarified.

\medskip
\subsection{Local family of A-infinity algebras}~\label{loc-fam-subsec} Again we let $V\subset \Ext^1(\EE,\EE)$ be a small neighborhood of zero. For each point $b\in V$ we define
\[ m^b_k(\alpha_1,\cdots,\alpha_k):=\sum_{j_0\geq 0,\cdots,j_k\geq 0} m_{k+j_0+\cdots+j_k}(b^{j_0},\alpha_1,\cdots,\alpha_k,b^{j_k}).\]
Lemma~\ref{norm-m-lem} implies the convergence of this series uniformly on $k$ once we choose $V$ small enough.

\begin{lem}
The structure maps $m_k^b$ defines an $A_\infty$ structure on the vector space $\Ext^*(\EE,\EE)$, for all $b\in V\subset \Ext^1(\EE,\EE)$.
\end{lem}

\Pf. Straightforward verification using the $A_\infty$ relation of $m_k$'s. \ed

Denote by $\Ext^*(\EE,\EE)^b$ this deformed $A_\infty$ structure. Note that the $A_\infty$ structure defined by $m_k^b$'s is in general curved, i.e. $m_0^b\neq 0$. In fact, by definition, it is precisely the Kuranishi map: $m_0^b=\kappa(b)$. 

By the Lemma above, we obtain a family of $A_\infty$ algebras on the bundle
\[ \AA_{\EE,V}:=\Ext^*(\EE,\EE)\otimes_\C \OO_V.\]
Its fiber over a point $b\in V$ is the $A_\infty$ algebra $\Ext^*(\EE,\EE)^b$.

\medskip
\subsection{Local family of L-infinity algebras}~\label{loc-lie-subsec} We shall mainly be concerned with $L_\infty$ algebras rather than $A_\infty$ algebras. 
Back to the context of $\Ext$-algebras, applying the functor $\Lie$ to the $A_\infty$ algebra bundle $\AA_{\EE,V}$ yields a local family of $L_\infty$ algebras $(\AA_{\EE,V})^{\Lie}$ over a small neighborhood $V\subset \Ext^1(\EE,\EE)$ of the origin.

We shall ultimately be interested in an $L_\infty$ sub-algebra of $(\AA_{\EE,V})^{\Lie}$ consisting of the {\sl traceless} part of it. More precisely, the natural morphism $$\OO_X\ra \EEnd(\EE)$$ admits a canonical splitting by the trace morphism $$\tr: \EEnd(\EE) \ra \OO_X,$$ which is in fact a splitting as Lie algebras $$\EEnd(\EE)\cong \EEnd_0(\EE)\oplus \OO_X$$ with $\EEnd_0(\EE)$ traceless endomorphisms, and $\OO_X$ endowed with the trivial Lie structure. Similarly, on cohomology, we have a direct sum decomposition $$\Ext^*(\EE,\EE)\cong \Ext^*_0(\EE,\EE)\oplus H^*(\OO_X).$$

In Corollary~\ref{lie-cor}, we proved that the subspace $\Ext^*_0(\EE,\EE)$ is a $L_\infty$ subalgebra of $\Ext^*(\EE,\EE)^{\Lie}$. This enables us to construct a family of curved $L_\infty$ algebras over a small open neighborhood $V\subset \Ext^1_0(\EE,\EE)$ of the origin as follows.

\begin{cons}\label{local-lie-cons}
Denote by $l_k$ the structure morphisms of the $L_\infty$ algebra $\Ext^*_0(\EE,\EE)$. Lemma~\ref{norm-m-lem} implies that there exists a constant $C>0$ independent of $k$, such that
\[ || l_k ||\leq k!\cdot C^k.\]
Let $V\subset \Ext^1_0(\EE,\EE)$ be a small open subset of origin, we can define a family of $L_\infty$ structure on the bundle $$\fkg_{\EE,V}:=\Ext^*_0(\EE,\EE)\ot_\C \OO_V$$ whose fiber over a point $b\in V$ is given by
\[ l_k^b(a_1,\cdots,a_k):=\sum_{j\geq 0} \frac{1}{j!} l_{k+j}(b^j,a_1,\cdots,a_k).\]
As in the $A_\infty$ case, we denote by $\Ext^*_0(\EE,\EE)^b$ this $L_\infty$ algebra. Note that $\Ext^*_0(\EE,\EE)^b$ is simply the $L_\infty$ subalgebra of the symmetrization of the $A_\infty$ algebra $\Ext^*(\EE,\EE)^b$. We shall use the notation 
$$\fkl_k: S^k ( \fkg_{\EE,V}[1] ) \ra \fkg_{\EE,V}[1] $$
for the family version of the higher brackets.
\end{cons}


\medskip
\subsection{Grothendieck connections}~\label{diff-lie-subsec} On the $L_\infty$ algebra bundle $\fkg_{\EE,V}=\Ext^*_0(\EE,\EE)\ot_\C\OO_V$ described in Construction~\ref{local-lie-cons}, there is a natural connection $\nabla$ whose flat sections are, by definition, given by $\Ext^*_0(\EE,\EE)\ot 1$. Furthermore, the $L_\infty$ structure on the bundle $\fkg_{\EE,V}$ satisfies certain compatibility condition with respect to $\nabla$. In this subsection, we explore this observation, which leads to a class of structured spaces: {\sl $L_\infty$ spaces}. This is a concept introduced by Costello in his geometric study of Witten genus~\cite{Costello}. The definition given below is more restrictive than that of Costello's~\cite[Definitions 2.1.1 and 2.1.2]{Costello}. See Remark~\ref{diff-rem} for a detailed discussion of the differences. 

Let us choose a basis $e_1,\cdots,e_m\in \Ext^1_0(\EE,\EE)$, and denote by $t_1,\cdots,t_m$ the dual linear coordinates on $V$, then differentiation shows that the $L_\infty$ structure and the connection $\nabla$ on $\fkg_{\EE,V}$ are compatible in the sense that the equation
\begin{equation}~\label{diff-eq} 
\nabla_{\partial/\partial t_i} \fkl_k (\alpha_1,\cdots,\alpha_k)= \fkl_{k+1} (e_i,\alpha_1,\cdots,\alpha_k)
\end{equation}
holds for flat sections $\alpha_1,\cdots,\alpha_k$ of $\fkg_{\EE,V}$.

We next give a more intrinsic way of describing this compatibility. 

\medskip
\begin{cons}~\label{main-local-cons}
Observe that the $L_\infty$ algebra bundle $\fkg_{\EE,V}$ is of the form
\[ \fkg_{\EE,V}=T_V[-1]\oplus (\fkg_{\EE,V})_2[-2] \oplus(\fkg_{\EE,V})_3[-3]\oplus\cdots\]
where $(\fkg_{\EE,V})_i$ denotes the degree $i$ component of $\fkg_{\EE,V}$. Thus its Chevalley-Eilenberg algebra is of the form
\[ C^*\fkg_{\EE,V}= \hat{S}(\fkg^\vee[-1])= \hat{S}\big( \Om_V\oplus (\fkg_{\EE,V})^\vee_2[1]\oplus\cdots\big).\]
The $L_\infty$ morphisms $\fkl_k$'s give rise to a square zero, degree one derivation $Q$ on $C^*\fkg_{\EE,V}$. The dual connection of $\nabla$ naturally extends, by Leibniz rule, to the Chevalley-Eilenberg algebra. We denote this extension by $\widetilde{\nabla}$. Furthermore, there is another natural connection $$\tau:C^*\fkg_{\EE,V} \ra \Om_V \ot C^*\fkg_{\EE,V}$$ defined on generators $\Om_V\oplus (\fkg_{\EE,V})^\vee_2[1]\oplus\cdots$ by
\[ \tau(\alpha)=\begin{cases}
\alpha\ot 1 &\mbox{if\;\; } \alpha \in \Om_V\\
0 &\mbox{if\;\; } \alpha \in (\fkg_{\EE,V})^\vee_i \mbox{ for $i\geq 2$.}
\end{cases}\]
\end{cons}
An unwinding of definitions proves the following
\medskip
\begin{lem}~\label{g-conn-lem}
Let us set $D:=\tau+\widetilde{\nabla}$. Then $D^2=0$, and Equation~\ref{diff-eq} is equivalent to the equation $[D,Q]=0$.
\end{lem}

This leads to the following definition of a class of structured spaces, which will serve as the local model of enhancements of moduli spaces of bundles.

\medskip
\begin{defi}~\label{diff-lie-defi}
Let $M$ be a complex manifold. Let $\fkg$ be a possibly curved $L_\infty$ algebra bundle, locally free of finite rank over $\OO_M$ which is of the form
\[ \fkg=T_M[-1]\oplus \fkg_2[-2] \oplus \cdots.\]
Denote by $(C^*\fkg,Q)$ its Chevalley-Eilenberg algebra. We refer to the grading of $C^*\fkg$ as the cohomological degree. A {\sl flat} connection $$D: C^*\fkg \ra \Om_M\ot C^*\fkg $$ on the bundle $C^*\fkg=\hat{S}(\Om_M\oplus\cdots)$, of {\sl cohomological degree zero}, is called a {\sl Grothendieck connection} if
\begin{itemize}
\item[(a.)] it is a derivation, i.e. $D(ab)=Da\cdot b + (-1)^{|a|} a\cdot Db$.
\item[(b.)] when restricted to the component $\Om_M \ra \Om_M\ot 1$, it is given by
\[ D(\alpha)=\tau(\alpha)=\alpha\ot 1.\] 
\item[(c.)] the connection $D$ is compatible with the $L_\infty$ structure in the sense that
\[ [D,Q]=0.\]
\end{itemize}
A $L_\infty$ space is a triple $(M,\fkg,D)$ with $D$ a Grothendieck connection on $C^*\fkg$.
\end{defi}

In this language, Lemma~\ref{g-conn-lem} implies that $D=\tau+\widetilde{\nabla}$ is a Grothendieck connection on $C^*\fkg_{\EE,V}$. Thus the triple $(V,\fkg_{\EE,V},D)$ forms a $L_\infty$ space.

\medskip
\begin{rem}~\label{diff-rem}
In~\cite[Definitions 2.1.1, 2.1.2]{Costello}, Costello introduced a class of structured spaces which he also called ``$L_\infty$ spaces". Roughly speaking, Costello's $L_\infty$ space (in the complex analytic setting) consists of a complex manifold $M$, together with a sheaf of $L_\infty$ algebras $\widetilde{\fkg}$ over the de Rham complex $\Om_M^*$. Via the bar construction (over the differential graded commutative ring $\Om_M^*$), one can show this structure is equivalent to a degree one, square zero derivation of $\hat{S}_{\Om_M^*} \big( \widetilde{\fkg}^\vee[-1]\big)$. Thus if $(M,\fkg,D)$ is a $L_\infty$ space in the sense of Definition~\ref{diff-lie-defi}, the sheaf $\widetilde{\fkg}:= \Om_M^*\ot_{\OO_M} \fkg$ forms a $L_\infty$ space in the sense of Costello, corresponding to the degree one, square zero derivation $Q+D$. Conversely, a $L_\infty$ space $(M,\widetilde{\fkg})$ in the sense of Costello gives rise to the data of a triple $(M,\fkg,D)$ provided that the structure morphisms of $\widetilde{\fkg}$ vanish in certain components. For this reason, the two definitions of a ``$L_\infty$ space" should cause no confusion.

On the mathematical level, the differences between Definition~\ref{diff-lie-defi} and Costello's are
\begin{itemize}
\item[(1.)] In Costello's version, he requires that the $L_\infty$ algebra $\fkg$, the reduction of $\widetilde{\fkg}$ by modulo the differential ideal $\langle\Om_M^{\geq 1}\rangle$, has no curvature term, while we do not require this condition. In the case of moduli spaces of bundles, this curvature term is precisely the Kuranishi maps. In a way, adding a curvature term to $\fkg$ allows us to use $L_\infty$ spaces as local resolutions for singular spaces.
\item[(2.)] We restrict ourselves to the case when $\fkg$ is strictly positively graded, with degree one component always the tangent bundle. This is similar to the condition of being negatively graded in Ciocan-Fontanine and Kapranov's theory of dg schemes~\cite{CK}. A precise relationship between $L_\infty$ spaces and dg schemes can be found in Subsection~\ref{ce-subsec} below.
\end{itemize}
\end{rem}

\medskip
\subsection{Chevalley-Eilenberg cohomology}~\label{ce-subsec} In this subsection, we explore the Chevalley-Eilenberg cohomology of a $L_\infty$ space $(M,\fkg,D)$. 

\begin{defi}\label{ce-lie-defi}
Let $(M,\fkg,D)$ be a $L_\infty$ space. Define its Chevalley-Eilenberg algebra
\begin{equation}~\label{ce-eq}
 C^*(\fkg,D):= \big( \Om_M^*\ot_{\OO_M} C^*\fkg, Q+D\big).
\end{equation}
 to be the de Rham complex with coefficients in the D-module $C^*\fkg$, endowed with the total differential.
\end{defi}

In order to calculate the Chevalley-Eilenberg cohomology $C^*(\fkg,D)$, we need some basic properties about the cohomology of the operator $D$. The proofs will be omitted since they are almost identical to results of~\cite[Section 2,3]{PT}. The moral of these properties is that the connection operator $D$ is simply the natural flat connection on an appropriate jet bundle, written in a trivialization.

\medskip
\begin{lem}~\label{koszul-lem}
Let $(\Lambda^* \Om_M \ot \hat{S} \Om_M,\tau)$ be the complex whose differential $\tau$ is the unique $\Lambda^*\Om_M$-linear derivation determined by
\[ \tau(1\ot \alpha)= \alpha\ot 1 \]
for $\alpha\in \Om_M$. Then the canonical maps
\begin{align*}
i: & \OO_M \rightarrow (\Lambda^* \Om_M \ot \hat{S} \Om_M,\tau)\\
\pi: & (\Lambda^* \Om_M \ot \hat{S} \Om_M,\tau) \rightarrow \OO_M
\end{align*}
are quasi-isomorphisms. Moreover, there exists a homotopy $h: \Lambda^* \Om_M \ot \hat{S} \ra \Lambda^{*-1} \Om_M \ot \hat{S}$ such that $\pi i=\id$ and $i\pi=\id+\tau h+h\tau$.
\end{lem}

Let $D: \hat{S}\Om_M \ra \Om_M \ot \hat{S}\Om_M$ be a connection on the bundle $\hat{S}\Om_M$ that is also a derivation. The derivation property implies that $D$ is uniquely determined by
\begin{align*}
&D_{-1}: \;\Om_M \ra \Om_M\ot 1\\
&D_{0}:  \; \Om_M \ra \Om_M \ot \Om_M\\
&D_{j}: \;\Om_M \ra \Om_M \ot S^{j+1} \Om_M, \;\;\forall j\geq 1.
\end{align*}

\medskip
\begin{lem}~\label{str-jet-lem}
Let $D$ be a connection on $\hat{S}\Om_M$ which is also a derivation. Assume that $\nabla$ is flat, and that $D_{-1}=\tau$. Then we have 
\[ \big( \Lambda^* \Om_M \ot \hat{S} \Om_M, D \big) \cong \OO_M\]
where $\cong$ denotes homotopy equivalence.
\end{lem}

More generally, we have the following

\medskip
\begin{lem}~\label{jet-lem}
Let $D$ be a connection on $\hat{S}\Om_M$ as in Lemma~\ref{str-jet-lem}. Let $(\FF,\nabla)$ be a bundle with a connection on $M$, not necessarily flat. Let $D_\FF$ be a degree one differential on $\Lambda^* \Om_M \ot \hat{S} \Om_M \ot \FF$, extending the connection map $\nabla$ on $\FF$, and making it into a left differential graded module over $\big(\Lambda^* \Om_M \ot \hat{S} \Om_M, D\big)$. Then we have
\[ \big( \Lambda^* \Om_M \ot \hat{S} \Om_M\ot \FF, D_\FF\big) \cong \FF.\]
\end{lem}

Let $(M,\fkg,D)$ be a $L_\infty$ space. Recall that the $L_\infty$ algebra bundle $\fkg$ is of the form
\[ \fkg=T_M[-1]\oplus \fkg^{\geq 2}\]
where $\fkg^{\geq 2}$ denotes the part of $\fkg$ with degrees bigger or equal to $2$. Its shifted dual decomposes as
\[ \fkg^\vee[-1]= \Omega_M^1 \oplus \WW \]
where $\WW=(\fkg^{\geq 2})^\vee[-1]$ is a strictly negatively graded bundle. This decomposition induces an isomorphism
\[ C^*\fkg=\hat{S}(\fkg^\vee[-1])\cong \hat{S}(\Omega_M^1)\ot \hat{S}(\WW).\]
Note that by degree reasons, the subbundle $\fkg^{\geq 2}$ is a sub-$L_\infty$ algebra, which implies that the differential $Q$ restricts to the subspace $$\hat{S}(\WW)\cong 1\ot \hat{S}(\WW)\subset C^*\fkg.$$

\medskip
\begin{thm}~\label{D-cohomology-prop}
Let $(M,\fkg,D)$ be a $L_\infty$ space. Then there exists a homotopy equivalence
\[ C^*(\fkg,D) \cong \big(\hat{S}(\WW),Q\big)\]
where the bundle $\WW$ is as described in the previous paragraph.
\end{thm}

\Pf. The complex $C^*(\fkg,\nabla)$ is the total complex of the bicomplex
\[\big( C^*(\fkg,\nabla), Q, D\big).\] 
Let us write down this bicomplex explicitly. The bidegree $(i,j)$ component is simply
\[ C^{i,j}:= \Om_M^j \ot \hat{S}(\Omega_M^1) \ot [\hat{S}(\WW)]_i, \]
The bicomplex then looks like

\medskip
\begin{equation}\label{bicomplex}
\begin{CD}
@AAA                    @A D AA @.   @A D AA  @A D AA  @.\\
\cdots  @>Q>> C^{i,j}  @>>>\cdots @>>> C^{-1,j} @>Q>>  C^{0,j} @>>> 0\\
@AAA                    @A D AA @.   @A D AA  @A D AA  @.\\
\vdots   @>Q>> \vdots  @>>> \cdots @>>> \vdots @>Q>> \vdots @>>> 0\\
@AAA                    @A D AA @.   @A D AA  @A D AA  @.\\
\cdots   @>Q>> C^{i,0}  @>>> \cdots @>>> C^{-1,0} @>Q>> C^{0,0} @>>> 0\\
\end{CD}
\end{equation}

\medskip
By Lemma~\ref{jet-lem} above, after taking vertical cohomology, we get a complex of the form
\begin{equation}\label{dg-eq}
\cdots [\hat{S}(\WW)]_i \ra [\hat{S}(\WW)]_{i+1}\ra \cdots \ra [\hat{S}(\WW)]_0=\OO_M\ra 0\ra \cdots.
\end{equation}
A diagram chasing shows that the differential on it matches with the restriction of $Q$.\ed

\begin{rem}
The proof of this theorem clarifies the relationship between our approach of derived structure versus Ciocan-Fontanine and Kapranov's theory of differential graded schemes: simply take the cohomology of $D$, i.e. flat sections with repsect to the Grothendieck connection. 
\end{rem}

Before we provide examples, we formalize the notion of ``$L_\infty$ enhancements".
\medskip

\begin{defi}~\label{lie-enh-defi}
Let $(U,\OO_U)$ be an analytic space. An $L_\infty$ enhancement of $(U,\OO_U)$ is a quadruple $(V,\fkg,D,\mathfrak{s})$ such that $(V,\fkg,D)$ is a $L_\infty$ space, and
\[ \fks: U \ra V\]
is an analytic morphism between analytic spaces. We require that the composition
\[ C^*(\fkg,D) \ra \OO_V \ra \fks_*\OO_U\]
is a morphism of complexes of sheaves, and induces an isomorphism 
\[ \underline{H}^0(C^*(\fkg,D)) \ra \fks_*\OO_U.\]
\end{defi}

\begin{ex}
Let us consider the case of moduli spaces of simple bundles  with fixed determinant on a $3$-fold, i.e. the $L_\infty$ spaces $(V,\fkg_{\EE,V},D)$ in Construction~\ref{main-local-cons}. The most interesting case is when the $L_\infty$ algebra
\[ \Ext_0^*(\EE,\EE)\]
consists of only two graded components in degree one and two. This is, for example, realized in the case of moduli spaces of bundles on projective Fano, or Calabi-Yau $3$-folds. Thus the bundle $\fkg=\fkg_{\EE,V}$ decomposes as
$$\fkg=T_V[-1] \oplus \fkg_2[-2].$$
Theorem~\ref{D-cohomology-prop} implies that the Chevalley-Eilenberg algebra $(C^*(\fkg,\nabla),Q)$ is quasi-isomorphic to
\[ \big(S(\fkg_2^\vee[1]), \lrcorner \kappa\big)\]
where $\kappa: \OO_V \ra \fkg_2$ is the curvature of $\fkg$ which is nothing but the Kuranishi map. In derived algebraic geometry, this latter algebra is often referred to as the derived zero locus of $\kappa$. By Kuranishi theory, this gives an enhancement of a neighborhood $U$ of $\EE$ in the moduli space.

If $X$ is a Calabi-Yau $3$-fold, then by Serre duality there is a non-degenerate pairing
\[ \Ext^1_0(\EE,\EE)\ot \Ext^2_0(\EE,\EE)\ra \C.\]
This implies the identification $\fkg_2\cong\Om_V$, and hence $\fkg_2^\vee\cong T_V$. Furthermore, the $L_\infty$ algebra structure is cyclic with respect to this pairing. Define the $L_\infty$ Chern-Simons functional by
\[ \Psi(b):=\sum_{j\geq 2} \frac{1}{(j+1)!} \langle l_j(b^j), b\rangle\]
for $b\in V$. Note that Lemma~\ref{lie-cor} implies the convergence of the series. By definition we have
\[ d\Psi=\kappa,\]
under the identification $\fkg_2\cong\Om_V$. Thus we obtain the commutative differential graded algebra 
\[ \big(S(T_V[1]), \lrcorner d\Psi\big),\]
known as the derived critical locus of $\Psi$.
\end{ex}

\section{The propagation functor}
\label{functor-sec}

Motivated by constructions appeared in the last section, we study in the current section some formal algebraic properties of a certain ``Functor" from the category of normed $L_\infty$ algebras to the category of $L_\infty$ spaces. We shall construct such a functor on the level of $\infty$-categories. The $\infty$-category language enables us to make the notion of {\em homotopy coherent gluing} of $L_\infty$ spaces precise.

\subsection{Normed L-infinity algebras}

We first formalize the convergence property appeared in the last section.

\begin{defi}\label{normed-alg-defi}
A finite dimensional $L_\infty$ algebra $g$ endowed with a norm $||\bullet||$ is called {\sl normed} if there exists a constant $C>0$, independent of $k$, such that
\[ ||l_k||\leq k!\cdot C^k.\]
\end{defi}

The following lemma shows that the normed condition is stable under perturbation.

\begin{lem}\label{propagation-lem}
There exists a neighborhood of zero $V\subset g_1$ such that, for each $b\in V$, the multi-linear maps defined by
\[ l_k^b(a_1,\cdots,a_k):=\sum_{j\geq 0} \frac{1}{j!} l_{k+j}(b^j,a_1,\cdots,a_k).\]
form another normed $L_\infty$ algebra structure on the graded vector space $g$ with the same norm.
\end{lem}

\Pf. The fact that $\{l_k^b\}_{k=0}^\infty$ form a $L_\infty$ algebra is a direct computation. It remains to prove that it is also bounded. For this we have
\begin{align*}
||l^b_k(a_1,\cdots,a_k)||&= || \sum_{j\geq 0} \frac{1}{j!} l_{k+j}(b^j,a_1,\cdots,a_k)||\\
& \leq \sum_{j\geq 0} \frac{1}{j!} (k+j)!\cdot C^{k+j} ||b||^j ||a_1||\cdots||a_k||\\
&= k!\big(\sum_{j\geq 0} {k+j \choose j} C^{k+j} ||b||^j\big) ||a_1||\cdots||a_k||
\end{align*}
Thus it suffices to show that
\[ \sum_{j\geq 0} {j+k \choose j} C^{k+j} ||b||^j \leq E^k\]
for some positive number $E$, independent of $k$. For this we use the inequality ${j+k\choose j} \leq 2^{j+k}$ to get
\[ \sum_{j\geq 0} {j+k \choose j} C^{k+j} ||b||^j\leq 2^{k}C^k\sum_{j\geq 0} (2C||b||)^j= \frac{2^{k}C^k}{1-2C||b||}.\]
If we choose $V$ small enough so that $\frac{1}{1-2C||b||}\leq 2$, then $E=4C$ would work.\ed

Assume that $g$ is strictly positively graded.  And let $V\subset g_1$ be a small enough neighborhood of zero. We consider the Kuranishi map $\kappa: V \ra g_2$ defined by
\[ \kappa(b):=\sum_k \frac{1}{k!} l_k(b^k).\]
The analytic space $\kappa^{-1}(0)$ admits a holomorphic $L_\infty$ enhancement by the exact same construction as in Construction~\ref{main-local-cons}.  Let us denote the $L_\infty$ space underlying this $L_\infty$ enhancement  by $(V,\fkg,D)$ where $\fkg:=g\otimes_\C \OO_V$.

In the following subsections, we would like to promote the association 
\[ g \mapsto \cS(g):=(V,\fkg, D)\]
to a ``functor" from normed $L_\infty$ algebras to $L_\infty$ spaces. We put functor in quote since the choice of $V$ is not canonical. One can introduce the notion of a {\sl germ} of a pointed $L_\infty$ spaces to solve this issue. However, we choose to not do so in order to avoid possible confusions. Instead we choose to describe in more concrete terms what we mean by being functorial.

\subsection{Bounded homomorphisms} In this subsection, we shall extend the construction of $\cS$ on the level of homomorphisms. First we define the notion of a bounded $L_\infty$ homomorphism. 

\medskip
\begin{defi}~\label{normed-hom-defi}
Let $g_\alpha$ and $g_\beta$ be two normed $L_\infty$ algebras.  An $L_\infty$ homomorphism
\[ f: g_\alpha \ra g_\beta\]
is called bounded if there exists a positive number $C>0$, independent of $k$, such that
\[ ||f_k||\leq k!\cdot C^k.\]
\end{defi}

\medskip
\begin{lem}~\label{bounded-compo-lem}
Let $f: g_\alpha \ra g_\beta$, and $h: g_\beta\ra g_\gamma$ be two bounded $L_\infty$ homomorphisms. Then the composition $h\circ f: g_\alpha\ra g_\gamma$ is also bounded.
\end{lem}

\Pf. Indeed, the degree $n$-component of $h\circ f$ satisfies
\begin{align*}
|| [h\circ f]_n || &= || \sum_{\substack{i_1,\cdots,i_k\geq 1\\ i_1+\cdots+i_k=n}} \sum_{\sigma\in \Sh(i_1,\cdots,i_k)} \frac{1}{k!} h_k ( f_{i_1}\ot\cdots\ot f_{i_k}) ||\\
&\leq \sum_{\substack{i_1,\cdots,i_k\geq 1\\ i_1+\cdots+i_k=n}} \frac{n!}{i_1!\cdots i_k!} \frac{1}{k!} k! C^k i_1! C^{i_1}\cdots i_k! C^{i_k}\\
&= (n!C^n) \sum_{\substack{i_1,\cdots,i_k\geq 1\\ i_1+\cdots+i_k=n}} C^k\\
&\leq (n!C^n) \sum_k {n-1\choose k-1} C^k\\
&\leq (n!C^n)(1+C)^n\\
&\leq n!\cdot (C(1+C))^n.
\end{align*}
Here $C$ is a positive constant such that $||f_k||\leq k!C^k$ and $||h_k||\leq k!C^k$.\ed

\medskip
\begin{defi}~\label{hom-defi}
Let $(V_\alpha,\fkg_\alpha, D^\alpha)$ and $(V_\beta,\fkg_\beta, D^\beta)$ be two $L_\infty$ spaces. A homomorphism between them consists of a pair
\[ (F,F^\sharp) : (V_\alpha,\fkg_\alpha,D^\alpha) \ra (V_\beta,\fkg_\beta,D^\beta),\]
where $F: V_\alpha \ra V_\beta$ is a morphism of analytic spaces, and
\[ F^\sharp : \fkg_\alpha\ra f^* \fkg_\beta\]
is a morphism of $L_\infty$ algebras, linear over $\OO_{V_\alpha}$. Furthermore, we require that $F^\sharp$ intertwines with $D^\alpha$ and $f^*D^\beta$.
\end{defi}

\medskip
\begin{cons}~\label{1-hom-cons}
Let $f: g_\alpha\ra g_\beta$ be a bounded $L_\infty$ homomorphism between two normed $L_\infty$ algebras. Associated to $g_\alpha$ and $g_\beta$ are the corresponding $L_\infty$ spaces $(V_\alpha,\fkg_\alpha, D^\alpha)$ and $(V_\beta,\fkg_\beta, D^\beta)$. We choose $V_\alpha$ small enough so that the convergent series
\[ F(b):= \sum_{k=1}^\infty f_k(b^k)\]
defines an analytic morphism $F: V_\alpha\ra V_\beta$. Finally we define the $L_\infty$ morphism $F^\sharp$ by
\[ F^\sharp_k(\alpha_1,\cdots,\alpha_k):= \sum_{j=0}^\infty \frac{1}{j!} f_{j+k}(b^j,\alpha_1,\cdots,\alpha_k).\]
It is straightforward to check that $(F,F^\sharp): (V_\alpha,\fkg_\alpha, D^\alpha) \ra (V_\beta,\fkg_\beta, D^\beta)$ is a morphism of $L_\infty$ spaces. Denote by $\cS(f)$ the morphism $(F,F^\sharp)$.
\end{cons}

The above construction is functorial in the sense that if we have two bounded $L_\infty$ homomorphisms $f: g_\alpha \ra g_\beta$, and $h: g_\beta\ra g_\gamma$. Then by shrinking $V_\alpha$ and $V_\beta$ if necessary, we can arrange so that
\[ F: V_\alpha \ra V_\beta, \mbox{\;\; and \;\;} H: V_\beta \ra V_\gamma.\]
Under this choice of the open sets $V$'s, it follows from the definitions that $\cS(h\circ f)=\cS(h)\circ\cS(f)$.

Since $L_\infty$ structures admit homotopy theory, it is natural to extend the construction of $\cS$ to the level of higher morphisms. To describe this extension precisely, we first describe the higher morphisms in both categories.

\medskip

\subsection{The simplicial category of normed L-infinity algebras}~\label{simp-lie-subsec} Let $A$, $B$ be two $L_\infty$ algebras (not necessarily normed for now). We define
\[ C(A,B):=\Hom\big( \overline{S}^c(A[1]), B\big) \]
where $\overline{S}^c(A[1])$ is the reduced bar construction of $A$. Observe that
\[ C(A,B)=\big(\overline{S}^c(A[1])\big)^\vee\ot B,\]
being the tensor product of commutative differential graded algebra and a $L_\infty$ algebra, is again a $L_\infty$ algebra. We denote by $L_k$ its structure morphisms. Since the reduced coproduct on $\overline{S}^c(A[1])$ is conilpotent, $C(A,B)$ is a pronilpotent $L_\infty$ algebra. Unwinding the definitions of $L_k$'s, we see that for a degree one element $\psi\in C(A,B)_1$, the $L_\infty$ Maurer-Cartan equation
\[ \sum_k \frac{1}{k!} L_k(\psi,\cdots,\psi)=0\]
is equivalent to the condition that $\psi$ is a $L_\infty$ homomorphism from $A$ to $B$.

We first describe the simplicial enrichment of the category of $L_\infty$ algebras, when the normed condition is not presented. Let $\Delta^n$ be the standard $n$-th simplex. We denote by $\Omega_{\Delta^n}^*$ the space of smooth differential forms on $\Delta^n$. For a simplicial space $Y$, we let $\Omega_Y^*$ to be the space of smooth simplicial differential forms on $Y$. Associated to a pronilpotent $L_\infty$ algebra $L$ is a simplicial set $$\MC_\bullet(L):=\MC( L\ot \Omega_{\Delta^\bullet}^*),$$
called the {\sl nerve of $L$}. In the case when $L=C(A,B)$, the simplicial set $$\MC_\bullet(C(A,B))$$ is usually considered as the mapping space between the $L_\infty$ algebras $A$ and $B$. Indeed, we have seen that its zero simplices are precisely $L_\infty$ homomorphisms.

Let us furthermore assume that both $A$ and $B$ are normed $L_\infty$ algebras. In order to accommodate the normed condition, we define a subspace of this mapping space $\MC_\bullet(C(A,B))$. Intuitively speaking, these are {\sl bounded} higher morphisms.

For a normed finite dimensional vector space $H$, denote by $||\bullet||_{m}$ the $C^m$-norm on the space $H\otimes \Om_{Y}^*$ of simplicial forms on $Y$ with values in $H$. 

\begin{defi}~\label{bounded-cochain-defi}
Let $Y$ be a simplicial space. For two normed finite dimensional vector spaces $A$ and $B$, we set the space of {\sl bounded cochains } from $A$ to $B\ot \Om_Y^*$ by 
\[
 C^b(A,B\ot \Om_{Y}^*):=\Big\{ (\psi_k)_{k=1}^\infty \in \Hom\big(\overline{S}^c(A[1]),B\otimes \Om_{Y}^*\big) \left.\middle|\right.  ||\psi_k||_{m}\leq D_m \cdot k!\cdot C^k\Big\}
\]
where $D_m>0$ only depends on $m$, and $C>0$ is independent of both $k$ and $m$. In other words, bounded cochains are absolutely convergent in the $C^\infty$ topology, under small perturbations.
\end{defi}

\medskip
 
\begin{lem}~\label{higher-comp-lem}
Let $\psi\in C^b(A,B\ot \Om_Y^*)$ and $\phi\in C^b(B,C\ot\Om_Z^*)$ be two bounded cochains. Then the composition $\phi\circ \psi$ defined by
\[ A\stackrel{\psi}{\ra} B\ot \Om_{Y}^* \stackrel{\phi\ot \id}{\longrightarrow} C\ot \Om^*_{Z}\ot\Om^*_{Y} \ra C\ot \Om^*_{Z\times Y}\]
is also a bounded cochain. Here the last arrow is the K\" unneth map.
\end{lem}

\Pf.  Analogous to Lemma~\ref{bounded-compo-lem}, plus the fact that pull-backs of differential forms are bounded in $C^\infty$ topology.\qed

\medskip
\begin{defi}~\label{inf-lie-defi}
We define a category $\NL$ enriched over the category of simplicial sets as follows. The objects of $\NL$ are normed $L_\infty$ algebras. For two objects $A$ and $B$, the set of $n$-simplices of the $\Hom$-simplicial set is given by
\[ \Hom(A,B)_n:=\MC^b_n(C(A,B))= \MC_n(C(A,B)) \cap C^b(A,B\ot \Om_{\Delta^n}^*).\]
Explicitly, an element $\psi$ of $\MC^b_n(C(A,B))$ is a bounded cochain
\[ \psi=\prod_{k=1}^\infty \psi_k: \prod_{k=1}^\infty A^{\ot k} \ra B\ot \Om_{\Delta^n}^*\]
which is also a $L_\infty$ homomorphism. The composition of morphisms in the category $\NL$ is well-defined by Lemma~\ref{higher-comp-lem}.
\end{defi}

For a pronilpotent $L_\infty$ algebra $L$, it was proved in~\cite{Getzler} that the simplicial set $\MC_\bullet(L)$ is a Kan complex. In our case, since $C(A,B)$ is pronilpotent, we conclude that the simplicial set $\MC_\bullet(C(A,B))$ is Kan. In the Appendix~\ref{kan-sec}, we prove the following

\medskip
\begin{prop}~\label{kan-prop}
The simplicial subset $$\Hom(A,B)_\bullet=MC^b_\bullet(C(A,B))\subset\MC_\bullet(C(A,B))$$ is also a Kan complex.
\end{prop}

Thus the category $\NL$ is a category enriched over Kan complexes. If $\CC$ is any category enriched over simplicial sets, Lurie~\cite[Definition 1.1.5.5]{Lur} defines a simplicial set $N(\CC)$, called the simplicial nerve of $\CC$. The set of $n$-simplices of $N(\CC)$ is given by
\[ \Hom\big(\mathfrak{C}[\Delta^n],\CC),\]
where $\mathfrak{C}[\Delta^n]$ is a simplicial category~\footnote{The category $\mathfrak{C}[\Delta^n]$ is a certain natural simplicial thickening of the ordinary category $[\Delta^n]$, see~\cite[Definition 1.1.5.1]{Lur} for its precise construction.}, and $\Hom$ is the set of simplicial functors from $\mathfrak{C}[\Delta^n]$ to $\CC$. Furthermore, it is shown in~\cite[Proposition 1.1.5.10]{Lur} that for a simplicial category $\CC$ whose $\Hom$ simplicial sets are Kan complexes, the resulting nerve $N(\CC)$ is always an $\infty$-category. Thus Proposition~\ref{kan-prop} implies the following

\medskip
\begin{cor}~\label{infinity-cor}
The simplicial set $N(\NL)$ is an $\infty$-category.
\end{cor}

\subsection{The simplicial category of L-infinity spaces}
We can also define a category $\LS$ enriched over simplicial sets whose objects are $L_\infty$ spaces. Namely let $(V_\alpha,\fkg_\alpha,D^\alpha)$ and $(V_\beta,\fkg_\beta,D^\beta)$ be two $L_\infty$ spaces. Define a simplicial set
\[ \Hom\big((V_\alpha,\fkg_\alpha,D^\alpha),(V_\beta,\fkg_\beta,D^\beta)\big)_\bullet\]
by setting its $n$-simplices to be the set of pairs $(F,F^\sharp)$ where
\[ F: V_\alpha \times \Delta^n \ra V_\beta\]
is a smooth morphism which is furthermore complex analytic in the $V_\alpha$-direction; and \[F^\sharp: (\pi_1)^*\fkg_\alpha \ra F^*\fkg_\beta \ot (\pi_2)^*\Om^*_{\Delta^n}\] 
is a morphism of $L_\infty$ algebra bundles which is compatible with the connections $(\pi_1)^*D^\alpha$ and $F^*D^\beta$.

We do not know if the simplicial set $\Hom\big((V_\alpha,\fkg_\alpha,D^\alpha),(V_\beta,\fkg_\beta,D^\beta)\big)_\bullet$ is a Kan complex. Nevertheless, the simplicial nerve construction applied to $\LS$ still yields a simplicial set $N(\LS)$.

Next we generalize Constructions~\ref{main-local-cons} and~\ref{1-hom-cons} to include higher morphisms.

\medskip
\begin{cons}~\label{s-functor-cons}
Let $g_\alpha$ and $g_\beta$ be two normed $L_\infty$ algebras, and let 
\[ f: g_\alpha \ra g_\beta \ot \Om^*_{\Delta^n}\]
be an element of $\Hom(g_\alpha,g_\beta)_n$. For each point $t\in \Delta^n$, the restriction
$$f|_{t}: g_\alpha\ra g_\beta$$
is a bounded $L_\infty$ homomorphism. We choose small enough open neighborhoods of origins $V_\alpha\subset (g_\alpha)_1$, and $V_\beta\subset (g_\beta)_1$ such that we get $L_\infty$ spaces $(V_\alpha,\fkg_\alpha,D^\alpha)$ and $(V_\beta,\fkg_\beta,D^\beta)$. Fruthermore, by the boundedness of $f$, we may shrink $V_\alpha$ and $V_\beta$ if necessary to ensure there is a well-defined map
\begin{align*}
F&: V_\alpha\times \Delta^n \ra V_\beta, \mbox{\; by putting\;\;}\\
F(b,t)&:= (f|_t)_*(b)= \sum_{k=1}^\infty (f|_t)_k(b^k).
\end{align*}
We then define a $L_\infty$ homomorphism
\[F^\sharp: (\pi_1)^*\fkg_\alpha \ra F^*\fkg_\beta \ot (\pi_2)^*\Om^*_{\Delta^n}\] 
by setting
\[ F^\sharp(\alpha_1,\cdots,\alpha_k):=\sum_{j\geq 0} \frac{1}{j!} f_{k+j}(b^j,\alpha_1,\cdots,\alpha_k).\]
Again we may shrink $V_\alpha$ to assure the convergence (in the $C^\infty$ topology) of the infinite sum, due to the uniform bounded of $f$. Again we denote by $\cS(f)$ the pair $(F,F^\sharp)$.
\end{cons}

\medskip
\begin{defi}~\label{over-defi}
Let $(V_\alpha,\fkg_\alpha,D^\alpha,\mathfrak{s}_\alpha)$ and $(V_\beta,\fkg_\beta,D^\beta,\mathfrak{s}_\beta)$ be two $L_\infty$ enhancements of a complex analytic space $U$. An $n$-morphism 
\[(F,F^\sharp)\in \Hom\big((V_\alpha,\fkg_\alpha,D^\alpha), (V_\beta,\fkg_\beta,D^\beta)\big)_n\]
is said to be {\sl over $U$} if the following diagram is commutative
\[\begin{CD}
V_\alpha\times \Delta^n @> F>> V_\beta\\
@A \mathfrak{s}_\alpha \times \id AA        @A \mathfrak{s}_\beta AA\\
U\times \Delta^n @> \pi >> U.
\end{CD}\]
\end{defi}

\medskip
\begin{prop}~\label{ind-ce}
Let $(V_\alpha,\fkg_\alpha,D^\alpha,\mathfrak{s}_\alpha)$ and $(V_\beta,\fkg_\beta,D^\beta,\mathfrak{s}_\beta)$ be two $L_\infty$ enhancements of a complex analytic space $U$. Let 
\[(F,F^\sharp)\in \Hom\big((V_\alpha,\fkg_\alpha,D^\alpha), (V_\beta,\fkg_\beta,D^\beta)\big)_0\]
be a $0$-morphism. Then it induces a morphism
\[ \underline{H}^*F: \mathfrak{s}_\beta^{-1} \underline{H}^*(V_\beta,\fkg_\beta,D^\beta) \ra \mathfrak{s}_\alpha^{-1} \underline{H}^*(V_\alpha,\fkg_\alpha,D^\alpha)\]
on the Chevalley-Eilenberg cohomology algebras.
\end{prop}

\Pf. By definition of a morphism, there is a morphism of algebras
\[ F^\sharp : F^*C^*\fkg_\beta \ra C^*\fkg_\alpha,\]
which is also a morphism of $D$-modules, with respect to the Grothendieck connections $F^*D_\beta$ and $D_\alpha$. Thus this map extends to the corresponding de Rham complexes
\[ F^\sharp : \Om_{V_\alpha}(F^*C^*\fkg_\beta) \ra \Om_{V_\alpha}(C^*\fkg_\alpha).\]
Since the $D$-module pull-back $F^*$ corresponds to the exact functor $F^{-1}$ under the de Rham functors, $F^\sharp$ induces a morphism
\[ F^{-1} \underline{H}^*(V_\beta,\fkg_\beta,D^\beta) \ra \underline{H}^*(V_\alpha,\fkg_\alpha,D^\alpha)\]
of algebras on $V_\alpha$. Now applying the inverse image functor $\mathfrak{s}_\alpha^{-1}$ both sides and using the identity $F\circ \mathfrak{s}_\alpha= \mathfrak{s}_\beta$ (because $F$ is a morphism of enhancements), we obtain the desired map $\underline{H}^*F$.
\ed

Next we prove that the two boundary $0$-morphisms of a $1$-morphism of enhancements induces the same map on the Chevalley-Eilenberg cohomology algebras.

\medskip
\begin{prop}~\label{inv-ce}
Let $(V_\alpha,\fkg_\alpha,D^\alpha,\mathfrak{s}_\alpha)$ and $(V_\beta,\fkg_\beta,D^\beta,\mathfrak{s}_\beta)$ be two $L_\infty$ enhancements of a complex analytic space $U$. Let 
\[(F,F^\sharp)\in \Hom\big((V_\alpha,\fkg_\alpha,D^\alpha), (V_\beta,\fkg_\beta,D^\beta)\big)_1\]
be a $1$-morphism. Denote by $(F_0,F_0^\sharp)$ and $(F_1,F_1^\sharp)$ the two boundaries of $(F,F^\sharp)$. Then we have
\[ \underline{H}^*F_0 = \underline{H}^*F_1.\]
\end{prop}

\Pf. By definition there is a map
\[ F^\sharp: F^*C^*\fkg_\beta \ra \pi_1^*C^*\fkg_\alpha\ot \pi_2^*\Omega_{\Delta^1}.\]
We may write $F^\sharp=A+B dt$ for two morphisms
\[ A, B: F^*C^*\fkg_\beta \ra \pi_1^*C^*\fkg_\alpha.\]
The fact that $F^\sharp$ is compatible with $F^*D_\beta$ and $\pi_1^*D_\alpha$ implies there are induced morphisms
\[ a, b: F^{-1} C^*(V_\beta, \fkg_\beta,D_\beta) \ra \pi_1^{-1} C^*(V_\alpha, \fkg_\alpha,D_\alpha)\]
which satisfies the identity
\begin{equation}~\label{1-mor-eq}
 \frac{da}{dt}=\pi_1^{-1} Q_\alpha \circ b+ b \circ  F^{-1} Q_\beta,\end{equation}
where $Q_\alpha$ and $Q_\beta$ are the Chevalley-Eilenberg differentials. Note that this equation is slightly different from the usual homotopy equation because the map $F: V_\alpha\times \Delta^1 \ra V_\beta$ may depend on the $t$-parameter. Nevertheless, the fact that $F$ is a morphism of enhancement (see the diagram in Definition~\ref{over-defi}) implies that \[(\mathfrak{s}_\alpha\times \id)^{-1}\circ F^{-1}= \pi^{-1}\circ \mathfrak{s}_\beta^{-1}.\]
Thus applying $(\mathfrak{s}_\alpha\times \id)^{-1}$ to Equation~\ref{1-mor-eq} yields
\[ \frac{d}{dt} a' = \pi^{-1}\mathfrak{s}_\alpha^{-1} Q_\alpha \circ b' + b' \circ \pi^{-1}\mathfrak{s}_\beta^{-1}Q_\beta,\]
where $a'=(\mathfrak{s}_\alpha\times \id)^{-1}a$ and $b'=(\mathfrak{s}_\alpha\times \id)^{-1} b$ are morphisms
\[ \pi^{-1}\mathfrak{s}_\beta^{-1} C^*(V_\beta,\fkg_\beta,D_\beta) \ra \pi^{-1}\mathfrak{s}_\alpha^{-1} C^*(V_\alpha,\fkg_\alpha,D_\alpha).\]
Now, integration over $[0,1]$ implies that
\[ a'(1)-a'(0)=[Q,\int_0^1b'dt],\]
which shows that $\underline{H}^*F_0=\underline{H}^*F_1$ since $a'(1)$ and $a'(0)$ induces $\underline{H}^*F_1$ and $\underline{H}^*F_0$ on cohomology.\ed

\medskip

\subsection{Homotopy L-infinity enhancements} In this subsection, we formulate a homotopy cocycle condition to glue local L-infinity enhancements by coherent homotopies. 

For this, it is useful to first recall the definition of a topological manifold $M$.  Let $\UU=\left\{U_i (i\in \II)\right\}$ an open covering of $M$, then an atlas on $M$ is given by homeomorphisms
\[ s_i: U_i \ra V_i \subset \R^N,\]
such that $s_j\circ s_i^{-1}: s_i(U_{ij}) \ra s_j(U_{ij})$ is continuous.

Let us translate this into the language of $\infty$-categories. Let $N(\UU)$ be the \emph{\v Cech nerve} of the covering $\UU$. The set of $0$-simplices of $N(\UU)$ is just the index set $\II$. The set of $1$-simplices of $N(\UU)$ consists of pairs of indices $(i_0,i_1)\in \II\times \II$ such that $U_{i_0}\cap U_{i_1}\neq \emptyset$. In general the set of $n$-simplices consists of $(n+1)$-tuples $(i_0,\cdots,i_n)$ such that $U_{i_0}\cap\cdots\cap U_{i_n}\neq \emptyset$. The simplicial maps are defined in the obvious way. Let $\CC$ be the category of open subsets of $\R^N$ with continuous morphisms, and denote by $N(\CC)$ its nerve. A $k-$simplex $\gamma\in N(\CC)_k$ is of the form
\[ V_0 \ra V_1 \ra \cdots \ra V_k.\]
It is said to be over an open subset $U\subset M$ if there are homeomorphisms $s_i: U\ra V_i, \; 0\leq i\leq k$ with the following diagram commutative
\[\begin{xy} 
(0,25)*+{V_0}="a"; (20,25)*+{V_1}="b";%
(60,25)*+{V_{k-1}}="c";
(80,25)*+{V_k}="d"; (40,0)*+{U}="e";
(40, 13)*+{\cdots};
(40,25)*+{\cdots}="f";
{\ar@{->}^{} "a";"b"};
{\ar@{->}_{} "b";"f"};
{\ar@{->}^{} "f";"c"};
{\ar@{->}^{s_0} "e";"a"};
{\ar@{->}^{s_1} "e";"b"};
{\ar@{->}_{s_{k-1}} "e";"c"};
{\ar@{->}_{s_k} "e";"d"};
{\ar@{->}^{} "c";"d"};
\end{xy}\]
For two simplices $\gamma$ and $\eta$ over $U$, the equality $\gamma=\eta$ means they are equal over $U$ (i.e. all homeomorphisms $s_i$'s are the same). 

For a $k-$simplex $\gamma$ over $U$, and another open subset $U'\subset U$, we define the restriction of $\gamma$ on $U'$, denoted by $\gamma\mid_{U'}$, to be the $k-$simplex
\[  s_0(U') \ra s_1(U') \ra \cdots \ra s_k(U').\]
This is naturally a simplex over $U'$. 

\begin{lem}
A topological atlas on $M$ subjected to the covering $\UU=\left\{U_i (i\in \II)\right\}$ is equivalent to an assignment
\[ \alpha\in N(\UU)_k\mapsto \Phi_k(\alpha)\in N(\CC)_k,\;\; \forall k\geq 0,\]
such that
\begin{itemize}
\item[--] the $k$-simplex $\Phi_k(\alpha)$ is over $U_\alpha$;
\item[--] the assignment is compatible with the simplicial structure in the following sense
\begin{align*}
 \partial_i (\Phi_k(\alpha)) &= \Phi_{k-1}(\partial_i\alpha)|_{U_\alpha}\mbox{ \;\; for \;\;} 0\leq i\leq k;\\
 \delta_j (\Phi_k(\alpha)) &= \Phi_{k+1}(\delta_j \alpha) \mbox{\;\; for \;\;} 0\leq j\leq k.
\end{align*}
Note that the second line makes sense because we have $U_{\alpha}=U_{\delta_j \alpha}$.
\end{itemize}
\end{lem}

We would like to formulate a global version of $L_\infty$ enhancements in a similar manner. Namely, on each open subset $U$, there is an enhancement $(V,\fkg,D,\fks)$; on double intersections, there should be transition maps; on triple intersections, there should be cocycle conditions. However, observe that the local objects that we are trying to glue admits homotopy theory, and thus it is only natural to have weak equivalences on double intersections, cocycle conditions up to homotopy on triple intersections, and coherent higher homotopies on more intersections. Moreover, instead of using the \v Cech coverings, we need to use {\sl hypercoverings}, which is also natural from a homotopy theory point of view. The differences between the two coverings are studied in~\cite{Dugger}.

\medskip
\begin{defi}
A hypercovering  of a topological space $M$ is a pair $(N_\bullet, \UU)$ where
\begin{itemize}
\item[(1.)] $N_\bullet$ is a simplicial set, called the nerve of the hypercovering;
\item[(2.)] $\UU$ is an assignment: for each $\alpha\in N_k$, there is an open subset $\UU(\alpha):=U_\alpha$ of the space $M$.
\end{itemize}
We denote by $\partial$'s and $\delta$'s the structure maps of $N_\bullet$. Then the pair $(N_\bullet,\UU)$ is subjected to the following properties
\begin{itemize}
\item[(a.)] for $0\leq j\leq k$, we have $U_\alpha\subset U_{\partial_j(\alpha)}$ for $\alpha\in N_k$;
\item[(b.)] for $0\leq j\leq k$, we have $U_\alpha=U_{\delta_j \alpha}$ for $\alpha\in N_k$;
\item[(c.)] $M=\bigcup_{\alpha \in N_0} U_\alpha$;
\item[(d.)] for $k\geq 1$ and every tuple $(\alpha_0,\cdots,\alpha_{k})\in (N_{k-1})^{k+1}$ such that $\partial_s \alpha_t=\partial_{t-1} \alpha_s$ with $0\leq s<t\leq k$, we have
\[\bigcap_{j=0}^k U_{\alpha_j}=\bigcup_{\substack{\alpha\in N_k,\\ \partial_j(\alpha)=\alpha_j, 0\leq j\leq k}} U_\alpha.\]
\end{itemize}
\end{defi}

\begin{defi}\label{h-lie-enh-defi}
Let $(M,\OO_M)$ be an analytic space. And let $(N_\bullet,\UU)$ be a hypercovering of $M$. A homotopy $L_\infty$ enhancement of $M$ subjected to the hypercovering $(N_\bullet,\UU)$ is given by an assignment
\[ \alpha\in N_k\mapsto \Phi_k(\alpha)\in N(\LS)_k, \forall k\geq 0.\]
Furthermore, we require that 
\begin{itemize}
\item[(i)] vertices of $\Phi_k(\alpha)$ are $L_\infty$ enhancements of $U_\alpha$;
\item[(ii)] all facets of the simplex $\Phi_k(\alpha)$ are defined over $U_\alpha$;
\item[(iii)] we have the following compatibility conditions 
\begin{align}~\label{cocycle}
\begin{split}
 \partial_i (\Phi_k(\alpha)) &= \Phi_{k-1}(\partial_i\alpha)|_{U_\alpha}\mbox{ \;\; for \;\;} 0\leq i\leq k;\\
 \delta_j (\Phi_k(\alpha)) &= \Phi_{k+1}(\delta_j \alpha) \mbox{\;\; for \;\;} 0\leq j\leq k.
 \end{split}
\end{align}
\end{itemize}
We use the notation ${\mathfrak M}:=(M,\Phi)$ to denote the enhanced space.
\end{defi}

\medskip
\begin{prop}~\label{derived-str-prop}
Let $(M,\OO_M)$ be an analytic space. Let $\Phi: N_\bullet \mapsto N(\LS)_\bullet$ be a homotopy $L_\infty$ enhancement of $M$ subjected to a hypercovering $(N_\bullet,\UU)$. Then the Chevalley-Eilenberg cohomology sheaves 
\[ \underline{H}^*\big( C^*( \Phi_0(\alpha))\big), \;\; \alpha \in N_0\]
glue into a global sheaf of non-positively graded algebras over $M$. We denote it by $\underline{H}^*( \fM)$. Furthermore, each of the graded components $\underline{H}^k(\fM)$ is coherent over $\OO_M$. If we assume that for each $\alpha\in N_0$, the $L_\infty$ algebra bundle $\fkg_\alpha$ is concentrated only in degree $1$ and $2$. Then the whole Chevalley-Eilenberg cohomology $\underline{H}^*(\fM)$ is coherent over $\OO_M$.
\end{prop}

\Pf. The existence of $\underline{H}^*( \fM)$ follows from Propositions~\ref{ind-ce} and~\ref{inv-ce}. The properties of $\underline{H}^*( \fM)$ follows from the local computation in Proposition~\ref{D-cohomology-prop}.\qed

\section{ Enhancements of moduli spaces}
\label{moduli-sec}

This section is devoted to the proof of the main Theorem~\ref{main-thm} and Corollary~\ref{main-cor}.

\subsection{Local constructions}
Let $E$ be a smooth complex vector bundle over a smooth projective variety $X$. Let $M$ be a moduli space of \textsl{{simple}} holomorphic structures on $E$ with fixed determinant. Assume also that $M$ is locally a universal family. We also fix a Hermitian metric on $E$ and a Kahler metric on $X$ to apply Hodge theory.

Let $q\in M$ be a point. We choose a complex structure $\db_q$ for the holomorphic vector bundle corresponding to $q$. By Corollaries~\ref{lie-cor},~\ref{i-cor},~\ref{p-cor} and Proposition~\ref{h-prop}, there exist
\begin{itemize}
\item[(A.)] a $L_\infty$ algebra structure on $\Ext_0^*(\EE_q,\EE_q)$;
\item[(B.)] a $L_\infty$ homomorphism
\[ \II : \Ext_0^*(\EE_q,\EE_q) \ra \Omega_X^{0,*}(\EEnd_0(\EE_q));\]
\item[(C.)] a $L_\infty$ homomorphism
\[ \PP: \Omega_X^{0,*}(\EEnd_0(\EE_q))_l \ra \Ext_0^*(\EE_q,\EE_q);\]
\item[(D.)] and $L_\infty$ homotopy 
\[ \HH: \Omega_X^{0,*}(\EEnd_0(\EE_q))_l \ra \Omega_X^{0,*}(\EEnd_0(\EE_q))_{l-1}\ot \Om_{[0,1]}^*\]
between $\id$ and $\II\circ \PP$.
\end{itemize}

For a point $b\in \Ext^1_0(\EE_\alpha,\EE_\alpha)$, we denote by
\[ B:= \sum_k \II_k(b,\cdots, b) \in \Omega_X^{0,1}(\EEnd_0(\EE))\]
where, by Corollary~\ref{lie-cor}, the series converges for $||b||$ small enough. We choose a small enough neighborhood  $U_q$ of the point $q\in M$ together with a Kuranishi chart
\[ {\mathfrak s}_q : U_q \hookrightarrow V_q\]
where $V_q\subset \Ext^1_0(\EE_q,\EE_q)$ is a small enough neighborhood of the origin, such that the perturbations
\[ l^b,\;\; \II^b, \;\; \PP^B, \;\; \mbox{ and } \HH^B\]
are well-defined. Note that this is possible due to the boundedness properties proved in Corollaries~\ref{lie-cor},~\ref{i-cor},~\ref{p-cor} and Proposition~\ref{h-prop}. 

By Construction~\ref{main-local-cons} we obtain an $L_\infty$ enhancement 
\[ {\mathfrak s}_q: (U_q,\OO_{U_q}) \hookrightarrow (V_q, \fkg_q, D^q).\]
Perform this construction for every point in $M$ yields the degree zero piece of the to-be-constructed homotopy $L_\infty$ enhancement of $M$. Namely we set
\[ N_0:= M,\]
and for each $q\in N_0$ we have an associated open subset $U_q$. Obviously this data satisfies the covering property
\[ M= \bigcup_{q\in N_0} U_q.\]
Define a map $\Phi: N_0 \ra \LS_0$ by the assignment
\[ \Phi_0: q\in N_0  \mapsto (V_q, \fkg_q, D^q).\]
As we have seen, by definition, the $L_\infty$ space $(V_q,\fkg_q, D^q)$ is an enhancement of $(U_q,\OO_{U_q})$ via the embedding ${\mathfrak{s_q}}$.

\subsection{Double intersections} Let $(i,j)\in N_0\times N_0$ be a pair of indices. Let $U_{ij}$ denote the intersection $U_i\bigcap U_j$. Fix a point $q\in U_{ij}$, denote by
\[ b_i:= \fks_i(q) \in V_i\subset \Ext_0^1(\EE_i,\EE_i), \mbox{\;\; and \;\;} b_j:=\fks_j(q) \in V_j\subset \Ext_0^1(\EE_j,\EE_j).\]
We need to compare the two $L_\infty$ algebras
\[ \fkg_i \mid_{b_i}= \Ext^*_0(\EE_i,\EE_i)^{b_i} \mbox{\;\; and \;\; } \fkg_j \mid_{b_j}= \Ext^*_0(\EE_j,\EE_j)^{b_j}.\]
Observe that the composition $$f_q^{ij}:=    \PP^{B_j}\circ \Ad(\psi_q^{ij}) \circ \II^{b_i}$$ defined by the diagram

\begin{equation}\label{double-diag}
\begin{CD}
\Om_X^{0,*}(\EEnd_0(\EE_i))^{B_i} @>\Ad(\psi_q^{ij})>> \Om_X^{0,*}(\EEnd_0(\EE_j))^{B_j} \\
@A \II^{b_i} AA                                                    @VV \PP^{B_j} V\\
\Ext^*_0(\EE_i,\EE_i)^{b_i}     @.         \Ext^*_0(\EE_j,\EE_j)^{b_j}
\end{CD}\end{equation}
gives a quasi-isomorphism. Here we used the notation $B_i=(\II)_* b_i$, $B_j=(\II)_* b_j$. Furthermore, the map $$\psi_q^{ij}: \big(\Om_X^{0,*}(E), \db_i+B_i\big) \ra \big(\Om_X^{0,*}(E), \db_j+B_j\big)$$ is an isomorphism whose existence is due to the fact that the two complexes both correspond to the point $q$ in the moduli space, and hence are isomorphic holomorphic structures. Note that the operator $\Ad(\psi_q^{ij})$ is independent of the choice of $\psi_q^{ij}$ by simpleness.

The $L_\infty$ morphism $f_q^{ij}$ is a bounded morphism in the $W^{l,2}$-norm if $l>\dim_\C X$, which follows from the boundedness of $\II$, $\PP$, and $\Ad(\psi_q^{ij})$. 

\medskip
\begin{lem}~\label{invertible-lem}
The morphism $f_q^{ij}$ is invertible in the homotopy category of $\NL$.
\end{lem}

\Pf. Note that this does not follow from that $f_q^{ij}$ is a quasi-isomorphism as we need to prove that it admits a bounded inverse, up to bounded homotopies. In our context, we define its homotopy inverse $f_q^{ji}$ by the composition
\[\begin{CD}
\Om_X^{0,*}(\EEnd_0(\EE_j))^{B_j} @>\Ad(\psi_q^{ji})>> \Om_X^{0,*}(\EEnd_0(\EE_i))^{B_i} \\
@A \II^{b_j} AA                                                    @VV \PP^{B_i} V\\
\Ext^*_0(\EE_j,\EE_j)^{b_j}     @.         \Ext^*_0(\EE_i,\EE_i)^{b_i}
\end{CD}\]
The morphism $f_q^{ji}$ is bounded, by the same reason as that of $f_q^{ij}$. 

Next we prove that $f_q^{ji}\circ f_q^{ij}$ is homotopic to $\id$ in $\NL$. For this, we construct a homotopy $\HH_q^{ij}$ by the composition
\begin{equation}~\label{homotopy-eq}
\begin{CD}
\Ext^*_0(\EE_i,\EE_i)^{b_i} @>\II^{b_i}>> \Om_X^{0,*}(\EEnd_0(\EE_i))^{B_i} \\
@. @VV\Ad(\psi_q^{ij})V\\
@. \Om_X^{0,*}(\EEnd_0(\EE_j))^{B_j} \\
@.   @VV\HH^{B_j}V\\
@. \Om_X^{0,*}(\EEnd_0(\EE_j))^{B_j}\ot \Om_{[0,1]}^*\\
@.   @VV\Ad(\psi_q^{ji}) V\\
\Ext^*_0(\EE_i,\EE_i)^{b_i} \ot\Om_{[0,1]}^* @<<\PP^{B_i}\ot \id< \Om_X^{0,*}(\EEnd_0(\EE_i))^{B_i} \ot \Om_{[0,1]}^*.
\end{CD}\end{equation}
The composition is bounded since also morphisms involved are bounded. Let us check the boundary conditions of $\HH_q^{ij}$. We have
\[ \HH_q^{ij}(0)= \II^{b_i}\Ad(\psi_q^{ij})\HH^{B_j}(0) \Ad(\psi_q^{ji}) \PP^{B_i}.\]
Since $\HH^{B_j}(0)=\id$, and $\Ad(\psi_q^{ij})\circ \Ad(\psi_q^{ij})=\id$ by the simpleness of $\EE_i$, we get
\[ \HH_q^{ij}(0)= \II^{b_i} \PP^{B_i}=\id,\]
as desired. The other boundary condition that
\[ \HH_q^{ij}(1)=f_q^{ji}\circ f_q^{ij}\]
follows from the identity $\HH^{B_j}(1)=\II^{b_j}\circ \PP^{B_j}$. The proof in the opposite direction that $f_q^{ij}\circ f_q^{ji}$ is homotopic to $\id$ is entirely similar.\qed

\medskip
Next we prove a property of the morphism $f_q^{ij}$ which allows us to ``move" the base point. Recall the definition of the $L_\infty$ homomorphism $$f_q^{ij}:=    \PP^{B_j}\circ \Ad(\psi_q^{ij}) \circ \II^{b_i}$$ from Diagram~\ref{double-diag}. Recall also the notations $b_i=\fks_i(q)$, $B_i=\II_*b_i$, $b_j=\fks_j(q)$, and $B_j=\II_* b_j$. 

Let $q'\in U_{ij}$ be another point which is sufficiently close to $q$ so that the perturbation
 \[ \big(f_q^{ij}\big)^{\fks_i(q')-\fks_i(q)}\]
is well defined. In the following, we use the notations $b_i':=\fks_i(q')$, $B_i':=\II_* b_i'$, $b_j':=\fks_j(q')$, and $B_j':=\II_* b_j'$.

\medskip
\begin{lem}~\label{moving-lem}
We have 
\[\big(f_q^{ij}\big)^{\fks_i(q')-\fks_i(q)}\cong f_{q'}^{ij},\]
where $\cong$ denotes homotopy equivalence in $\NL$.
\end{lem}
\Pf. The perturbation of the first morphism
\[ \II^{b_i}: \Ext_0^*(\EE_i,\EE_i)^{b_i} \ra \Om_X^{0,*}(\EEnd_0(\EE_i))^{B_i}\]
by $\Delta b_i=b_i'-b_i$ gives 
\[ \II^{b_i'}: \Ext_0^*(\EE_i,\EE_i)^{b_i'} \ra \Om_X^{0,*}(\EEnd_0(\EE_i))^{B_i'}\]
where we have used that $(\II^{b_i})_*\Delta b_i = B_i'-B_i$. Perturbing the second morphism
\[ \Ad(\psi_q^{ij}): \Om_X^{0,*}(\EEnd_0(\EE_i))^{B_i} \ra \Om_X^{0,*}(\EEnd_0(\EE_j))^{B_j}\]
by $\Delta B_i:=(\II^{b_i})_*\Delta b_i = B_i'-B_i$ yields
\[ \Ad(\psi_q^{ij}): \Om_X^{0,*}(\EEnd_0(\EE_i))^{B_i'} \ra \Om_X^{0,*}(\EEnd_0(\EE_j))^{B_j+\Ad(\psi_q^{ij})_*\Delta B_i}.\]
Note that the morphism $\Ad(\psi_q^{ij})$, being a strict morphism, does not change under perturbations. 

The key is to compare the Maurer-Cartan element $\Ad(\psi_q^{ij})_*\Delta B_i$ with $\Delta B_j= (\II^{b_j})_*\Delta b_j=B_j'-B_j$. Since they give rise to isomorphic holomorphic structures corresponding to the point $q'$ in moduli space, the two Maurer-Cartan elements are gauge-equivalent. Furthermore, by choosing $\Delta b_i $ small enough, there exists an element $\eta\in \Om_X^{0,0}(\EEnd_0(\EE_j))^{B_j}$ in the Lie algebra of the gauge group, such that
\[ e^{\eta}* [\Delta B_j] = \Ad(\psi_q^{ij})_*\Delta B_i.\]
This defines a Maurer-Cartan element 
\[ \theta:= e^{t\eta}*[\Delta B_j] + \eta dt \in \Om_X^{0,0}(\EEnd_0(\EE_j))^{B_j} \ot \Om_{[0,1]}^*.\]
Using $\theta$ to perturb the $L_\infty$ morphism 
\[ \PP^{B_j}\ot \id : \Om_X^{0,0}(\EEnd_0(\EE_j))^{B_j} \ot \Om_{[0,1]}^* \ra \Ext_0^*(\EE_j,\EE_j)^{b_j}\ot \Om_{[0,1]}^*,\]
we obtain a $L_\infty$ morphism
\[ (\PP^{B_j}\ot \id)^\theta: [\Om_X^{0,0}(\EEnd_0(\EE_j))^{B_j} \ot \Om_{[0,1]}^*]^\theta \ra \Ext_0^*(\EE_j,\EE_j)^{b_j'}\ot \Om_{[0,1]}^*.\]
We have used the fact that $(\PP^{B_j}\ot\id)_*\theta=\Delta b_j$ is independent of $t$ due to the simpleness of $\EE_j$.

Putting everything together, we consider the following composition
\begin{equation}~\label{moving-eq}\begin{CD}
\Ext^*_0(\EE_i,\EE_i)^{b_i'} @> \II^{b_i'}>> \Om_X^{0,*}(\EEnd_0(\EE_i))^{B_i'} \\
@.                            @VV \Ad(\psi_{q'}^{ij}) V\\
      @.              \Om_X^{0,*}(\EEnd_0(\EE_j))^{B_j'}  \\
@.                             @VV\Ad( e^{t\eta}) V \\
\Ext_0^*(\EE_j,\EE_j)^{b_j'}\ot \Om_{[0,1]}^*      @< (\PP^{B_j}\ot\id)^\theta <<            [\Om_X^{0,0}(\EEnd_0(\EE_j))^{B_j} \ot \Om_{[0,1]}^*]^\theta
\end{CD}\end{equation}
When the above diagram specializes to $t=0$, it gives $f_{q'}^{ij}$; while when $t=1$, it is $(f_q^{ij})^{b_i'-b_i}$. \qed

Applying the Construction~\ref{1-hom-cons} to the bounded morphism $f_q^{ij}$, we obtain a morphism of $L_\infty$ spaces
\[ \cS(f_q^{ij}): (V_i,\fkg_i,D^i)\mid_{V_{ij}} \ra (V_j,\fkg_j,D^j)\mid_{V_{ij}^\dagger},\]
where we arrange the open subsets $V_{ij}\subset V_i$ and $V_{ij}^\dagger\subset V_j$ so that
\[ \fks_i^{-1}(V_{ij}) = \fks_j^{-1} (V_{ij}^\dagger) \subset U_{ij}.\]
Denote this common open neighborhood of $q$ by $U_{(i,j,q)}$. By shrinking $U_{(i,j,q)}$ if necessary, we require that Lemma~\ref{moving-lem} above holds for any $q'\in U_{(i,j,q)}$.

\noindent Summarizing, let us set
\[ N_1:= \left\{ (i,j,q)\in N_0\times N_0\times M\mid q\in U_{ij}.\right\} \]
There are two boundary maps $N_1\ra N_0$ defined by
\[ (i,j,q) \mapsto i \mbox{\;\; and \;\; } (i,j,q) \mapsto j.\]
The degeneracy map $N_0\ra N_1$ is given by
\[ i \mapsto (i,i, i) \mbox{\;\; (recall that $N_0=M$).}\]
The covering property is obvious:
\[ U_{ij} = \bigcup_{q\in U_{ij}} U_{(i,j,q)}.\]
The construction described above gives a map
\[ \Phi_1: N_1 \ra \LS_1 \mbox{\;\; by the assignment \;\;} \Phi_1((i,j,q))= \cS(f_q^{ij}),\]
which gives a morphism of enhancements defined over $U_{(i,j,q)}$ with the desired compatibility conditions on its boundaries.

\subsection{Triple intersections} Let $\alpha,\beta,\gamma\in N(\UU)_1$ be three indices of the form
\begin{align}~\label{triple-cover}
\begin{split}
\alpha &=(i,j,q_\alpha)\\
\beta &=(j,k, q_\beta)\\
\gamma &=(i,k,q_\gamma).
\end{split}
\end{align} 
This configuration is depicted as the following picture.
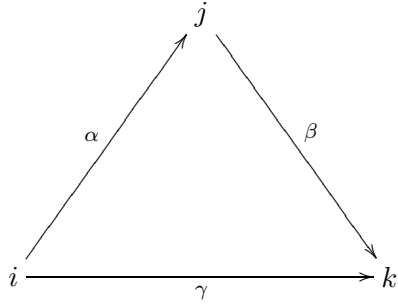
\[\begin{xy} 
(0,0)*+{i}="a"; (50,0)*+{k}="b";%
(25,35)*+{j}="c"; 
{\ar@{->}^{\alpha} "a";"c"};
{\ar@{->}_{\gamma} "a";"b"}
{\ar@{->}^{\beta} "c";"b"}
\end{xy}\]
Fix a point $q\in U_{\alpha\beta\gamma}=U_\alpha\bigcap U_\beta\bigcap U_\gamma$. We would like to prove $\Phi_1(\beta)\circ\Phi_1(\alpha)$ and $\Phi_1(\gamma)$ are homotopic after restricting on an appropriate neighborhood of $q$. 

Observe that the morphisms $\Phi_1(\alpha)$, $\Phi_1(\beta)$, and $\Phi_1(\gamma)$, when restricted to the point $q$, are three bounded $L_\infty$ homomorphisms
\begin{align*}
(f_{q_\alpha}^{ij})^{\Delta b_i}:& \Ext^*_0(\EE_i,\EE_i)^{\fks_i(q)} \ra \Ext^*_0(\EE_j,\EE_j)^{\fks_j(q)},\\
(f_{q_\beta}^{jk})^{\Delta b_j}:& \Ext^*_0(\EE_j,\EE_j)^{\fks_j(q)} \ra \Ext^*_0(\EE_k,\EE_k)^{\fks_k(q)},\\
(f_{q_\gamma}^{ik})^{\Delta b_i}:& \Ext^*_0(\EE_i,\EE_i)^{\fks_i(q)} \ra \Ext^*_0(\EE_k,\EE_k)^{\fks_k(q)}.
\end{align*}
Here the $\Delta b$'s are defined by
\begin{align*}
\Delta b_i&:= \fks_i(q)-\fks_i(q_\alpha),\\
\Delta b_j&:=\fks_j(q)-\fks_j(q_\beta)
\end{align*}

\medskip
\begin{lem}~\label{triple-lem}
The composition $(f_{q_\beta}^{jk})^{\Delta b_j}\circ (f_{q_\alpha}^{ij})^{\Delta b_i}$ is homotopic to $(f_{q_\gamma}^{ik})^{\Delta b_i}$ in the category $\NL$.
\end{lem}

\Pf. Recall that $\NL$ is a simplicial category. Denote by $N(\NL)$ its simplicial nerve. The three morphisms gives us a map of simplicial set
\[ \partial \Delta^2 \ra N(\NL).\]
By Lemma~\ref{moving-lem}, this fits into a bigger diagram
\[ \Delta^1\times \partial\Delta^2 \ra N(\NL),\]
which is depicted as
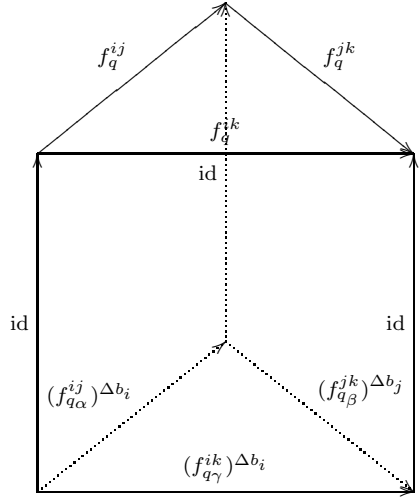
\[\begin{xy} 
(0,0)*{}="A"; (50,0)*{}="C"; (25,20)*{}="B";
(0,45)*{}="D"; (50,45)*{}="F"; (25,65)*{}="E";
{\ar@{->}^{(f_{q_\gamma}^{ik})^{\Delta b_i}} "A"; "C"};
{\ar@{->}^{\id} "A"; "D"};
{\ar@{.>}^{(f_{q_\alpha}^{ij})^{\Delta b_i}} "A"; "B"};
{\ar@{.>}^{(f_{q_\beta}^{jk})^{\Delta b_j}} "B"; "C"};
{\ar@{->}^{f_q^{ij}} "D"; "E"};
{\ar@{->}^{f_q^{jk}} "E"; "F"};
{\ar@{->}^{f_q^{ik}} "D"; "F"};
{\ar@{.>}^{\id} "B"; "E"};
{\ar@{->}^{\id} "C"; "F"};
\end{xy}\]
In the above diagram, the vertical edges are $\id$'s, and the vertical planes are filled by homotopies constructed in Lemma~\ref{moving-lem}. 

Next, we prove that the top triangle admits a filling. We set
\[ b_i:=\fks_i(q), \mbox{ \;\; and \;\; } B_i:=\II_* b_i.\]
Similarly we have $b_j$, $B_j$, $b_k$, and $B_k$. By definition, we have
\[f_q^{jk} \circ f_q^{ij}= \PP^{B_k}\Ad(\psi_q^{jk}) \II^{b_j} \PP^{B_j}\Ad(\psi_q^{ij})\II^{b_i}.\]
The middle part $\II^{b_j} \PP^{B_j}$ is homotopic to $\id$ via $\HH^{B_j}$. Thus, the composition
\[ \PP^{B_k}\Ad(\psi_q^{jk}) \HH^{B_j} \Ad(\psi_q^{ij})\II^{b_i}\]
gives the desired filling of the top triangle. Note that to check the boundary condition, we need to use the identity
\[ \Ad(\psi_q^{jk})\Ad(\psi_q^{ij})=\Ad(\psi_q^{ik}),\]
which follows from the simpleness assumption. Putting the top filling into the diagram above, we obtain a morphism of simplicial set
\[(\Delta^1\times \partial \Delta^2) \bigcup_{\left\{1\right\}\times\partial\Delta^2} (\left\{1\right\} \times \Delta^2) \ra N(\NL).\]
By Lemma~\ref{invertible-lem}, the edges of this morphism are homotopy equivalences in $\NL$, which implies that the image in fact lies in a Kan complex $N(\NL)'\subset N(\NL)$, the $\infty$-groupoid associated to the $\infty$-category $N(\NL)$. The lemma then follows from the homotopy lifting property of Kan complexes since the inclusion
\[ (\Delta^1\times \partial \Delta^2) \bigcup_{\left\{1\right\}\times\partial\Delta^2} (\left\{1\right\} \times \Delta^2)  \ra \Delta^1\times \Delta^2\]
is an acyclic cofibration.\qed

Next we prove a lemma that is similar to Lemma~\ref{moving-lem} in order to ``move" homotopies between $L_\infty$ homomorphisms. Namely, we consider another point $q'\in U_{\alpha\beta\gamma}$ sufficiently close to $q$. In the proof the previous lemma, we constructed a $2$-simplex in $N(\NL)$ depicted as
\[\scalebox{0.8} {
\begin{xy} 
(0,0)*+{\Ext_0^*(\EE_i,\EE_i)^{\fks_i(q)}}="a"; (100,0)*+{\Ext_0^*(\EE_k,\EE_k)^{\fks_k(q)}}="b";%
(50,70)*+{\Ext_0^*(\EE_j,\EE_j)^{\fks_j(q)}}="c"; 
(50,30)*+{\PP^{B_k}\Ad(\psi_q^{jk}) \HH^{B_j} \Ad(\psi_q^{ij})\II^{b_i}};
{\ar@{->}^{f_q^{ij}} "a";"c"};
{\ar@{->}_{f_q^{ik}} "a";"b"}
{\ar@{->}^{f_q^{jk}} "c";"b"}
\end{xy} }\]
Denote this by a map $f_q^{ijk}: \Delta^2\ra N(\NL)$. Similarly, this construction applied to the point $q'$ yields a map $f_{q'}^{ijk}: \Delta^2\ra N(\NL)$, corresponding to the following picture.
\[\scalebox{0.8}{\begin{xy} 
(0,0)*+{\Ext_0^*(\EE_i,\EE_i)^{\fks_i(q')}}="a"; (100,0)*+{\Ext_0^*(\EE_k,\EE_k)^{\fks_k(q')}}="b";%
(50,70)*+{\Ext_0^*(\EE_j,\EE_j)^{\fks_j(q')}}="c"; 
(50,30)*+{\PP^{B_k'}\Ad(\psi_{q'}^{jk}) \HH^{B_j'} \Ad(\psi_{q'}^{ij})\II^{b_i'}};
{\ar@{->}^{f_{q'}^{ij}} "a";"c"};
{\ar@{->}_{f_{q'}^{ik}} "a";"b"}
{\ar@{->}^{f_{q'}^{jk}} "c";"b"}
\end{xy}}\]
 
\medskip
\begin{lem}\label{moving-homotopy}
In the above setup, for $q'$ sufficiently close to $q$, we have
\[ \big( f_q^{ijk} \big)^{\fks_i(q')-\fks_i(q)} \cong f_{q'}^{ijk}.\] 
That is, there exists a morphism
\[ \mu: \Delta^1\times \Delta^2 \ra N(\NL),\]
such that $\mu(0)=\big( f_q^{ijk} \big)^{\fks_i(q')-\fks_i(q)}$ and $\mu(1)=f_{q'}^{ijk}$.
\end{lem}

\Pf. After unwinding the definition of the simplicial nerve construction, to construct the homotopy $\mu$ is equivalent to construct a map of simplicial sets
\[ \chi: \Delta^1\times \Delta^1 \ra \Hom_{\NL}\big(\Ext_0^*(\EE_i,\EE_i)^{\fks_i(q')}, \Ext_0^*(\EE_k,\EE_k)^{\fks_k(q')}\big)\]
whose boundary is the following diagram
\[\scalebox{1}{\begin{xy}
(0,0)*+{f_{q'}^{ik}}="A"; 
(50,0)*+{ (f_q^{ik})^{\fks_i(q')-\fks_i(q)}}="B";
(0,50)*+{f_{q'}^{jk}\circ f_{q'}^{ij}}="C"; 
(50,50)*+{ (f_q^{jk})^{\fks_j(q')-\fks_j(q)} \circ (f_q^{ij})^{\fks_i(q')-\fks_i(q)}} ="D";
{\ar@{->} "A"; "B"}
{\ar@{->}^{\PP^{B_k'}\Ad(\psi_{q'}^{jk}) \HH^{B_j'} \Ad(\psi_{q'}^{ij})\II^{b_i'}} "A"; "C"}
{\ar@{->} "C"; "D"}
{\ar@{->}_{\big(\PP^{B_k}\Ad(\psi_q^{jk}) \HH^{B_j} \Ad(\psi_q^{ij})\II^{b_i}\big)^{\fks_i(q')-\fks_i(q)}} "B"; "D"}
\end{xy}}\]
with the horizontal arrows defined as in Diagram~\ref{moving-eq}.

The two parameter family $\chi$ is constructed as follows. First, we observe the following three parameter family $\widetilde{\chi}$ defined by the composition
\[ \begin{CD}
\Ext_0^*(\EE_i,\EE_i)^{b_i'} @>\widetilde{\chi}>> \Ext_0^*(\EE_k, \EE_k)^{b_k'}\ot \Om_{[0,1]\times [0,1]\times [0,1]}^*\\
@V\II^{b_i'} VV                            @AA [\PP^{B_k}]^{\theta_2} A\\
\Om_X^{0,*}(\EEnd_0(\EE_i))^{B_i'}       @. \big(\Om_X^{0,*}(\EEnd_0(\EE_k))^{B_k'} \ot \Om_{[0,1]\times [0,1]\times [0,1]}^*\big)^{\theta_2}\\
@V \Ad(\psi_{q'}^{ij}) VV            @AA \Ad(e^{t_2\eta_2}) A\\
 \Om_X^{0,*}(\EEnd_0(\EE_j))^{B_j'}   @.    \Om_X^{0,*}(\EEnd_0(\EE_k))^{B_k'} \ot \Om_{[0,1]\times [0,1]}^*\\
 @V \Ad(e^{t_1\eta_1}) VV   @AA \Ad(\psi_{q'}^{jk}) A \\
[ \Om_X^{0,*}(\EEnd_0(\EE_j))^{B_j}\ot \Om_{[0,1]}^*]^{\theta_1} @>>[\HH^{B_j}]^{\theta_1}>  \Om_X^{0,*}(\EEnd_0(\EE_j))^{B_j'} \ot \Om_{[0,1]\times [0,1]}^*\\
\end{CD}\]
Here the $\eta$'s and the $\theta$'s are as defined in the proof of Lemma~\ref{moving-lem}. As shown above, total composition gives us a morphism 
\[ \widetilde{\chi}: \Delta^1\times\Delta^1\times\Delta^1 \ra \Hom_{\NL}\big(\Ext_0^*(\EE_i,\EE_i)^{\fks_i(q')}, \Ext_0^*(\EE_k,\EE_k)^{\fks_k(q')}\big).\]
We define 
\[ \chi:=\widetilde{\chi}|_{t_1=t_2}\]
where the coordinates $t_1$ and $t_2$ are as used in the above diagram. From the definition, it is ready to verify all the boundary conditions of $\chi$.\qed

\medskip
Finally, we can proceed to construct the set $N_2$, together with the map $\Phi_2: N_2\ra N(\LS)$. We define the set
\[ N_2:=\left\{ (\alpha,\beta,\gamma, q)\in N_1\times N_1\times N_1\times M \mid \mbox{$\alpha$, $\beta$, $\gamma$ are as in~\ref{triple-cover}, and \;} q\in U_{\alpha\beta\gamma}\right\}.\]
The three boundary maps $N_2\ra N_1$ are the three projections. The two degeneracy maps $N_1\ra N_2$ are given by
\begin{align*}
\alpha=(i,j,q) &\mapsto \big(\alpha,\alpha,(j,j,j), q\big), \mbox{\;\; and \;\;} \\
\alpha=(i,j,q) &\mapsto \big((i,i,i), \alpha,\alpha, q \big).
\end{align*}
For each $(\alpha,\beta,\gamma, q) \in N_2$, let 
$$f_{(\alpha,\beta,\gamma,q)}: \Delta^2 \ra \Hom_{\NL}\big(\Ext_0^*(\EE_i,\EE_i)^{\fks_i(q)}, \Ext_0^*(\EE_k,\EE_k)^{\fks_k(q)}\big)$$ be any filling of the diagram
\[\begin{xy} 
(0,0)*{\Ext_0^*(\EE_i,\EE_i)^{\fks_i(q)}}="A"; (70,0)*{\Ext_0^*(\EE_k,\EE_k)^{\fks_k(q)}}="C"; (35,60)*{\Ext_0^*(\EE_j,\EE_j)^{\fks_j(q)}}="B";
(35,30)*{f_{(\alpha,\beta,\gamma,q)}}="D";
{\ar@{->}^{\big(f_{q_\gamma}^{ik}\big)^{\fks_i(q)-\fks_i(q_\gamma)}} "A"; "C"};
{\ar@{->}^{\big(f_{q_\alpha}^{ij}\big)^{\fks_i(q)-\fks_i(q_\alpha)}} "A"; "B"};
{\ar@{->}^{\big(f_{q_\beta}^{jk}\big)^{\fks_j(q)-\fks_j(q_\beta)}} "B"; "C"};
\end{xy}\]
constructed in the proof of Lemma~\ref{triple-lem}. We may choose a small enough open subset $U_{(\alpha,\beta,\gamma,q)}\subset M$, together with small enough open subsets
\[ V_{(\alpha,\beta,\gamma,q)}\subset V_i, \;\; V_{(\alpha,\beta,\gamma,q)}^\dagger \subset V_j, \;\; V_{(\alpha,\beta,\gamma,q)}^{\dagger\dagger}\subset V_k,\]
such that 
\[\fks_i^{-1}(V_{(\alpha,\beta,\gamma,q)})=\fks_j^{-1}(V_{(\alpha,\beta,\gamma,q)}^\dagger)=\fks_k^{-1}(V_{(\alpha,\beta,\gamma,q)}^{\dagger\dagger})=U_{(\alpha,\beta,\gamma,q)},\]
and set 
\[ \Phi_2\big((\alpha,\beta,\gamma,q)\big):=\cS(f_{(\alpha,\beta,\gamma,q)})\]
with the morphism $\cS(f)$ defined over $U_{(\alpha,\beta,\gamma,q)}$.

\subsection{Higher homotopies}
The constructions on double and triple intersections illustrate how we should proceed on multiple intersections. Indeed, assume that we have constructed the nerve of a hypercovering up to degree $k-1$ inductively by
\begin{align*}
 N_{j}:=& \big\{ \alpha=(\alpha_0,\cdots,\alpha_j,q_\alpha)\in (N_{j-1})^{j+1}\times M\mid \\
 &\partial_s \alpha_t=\partial_{t-1} \alpha_s, \;\forall 0\leq s<t\leq j, \mbox{\; and \;} q_\alpha\in U_{\alpha_0}\cap\cdots\cap U_{\alpha_j}.\big\}
 \end{align*}
for $0\leq j \leq k-1$. We refer to $q_\alpha$ as the base point of $U_\alpha$. Assume also that required maps
\[ \Phi_{j}: N_{j} \ra N(\LS)_{j}\]
which satisfy the homotopy version of cocycle condition~\ref{cocycle} are also defined for $0\leq j\leq k-1$. Furthermore, for each index $\alpha\in N_j$ whose vertices are $i_0, i_1,\cdots, i_{j-1}, i_j\in N_0$, we have that 
\[ \Phi_{j}(\alpha)=\cS(f_\alpha),\]
for some appropriate simplex $f_\alpha\in N(\NL)_j$ whose vertices are given by normed $L_\infty$ algebras
\[ \Ext_0^*(\EE_{i_0},\EE_{i_0})^{\fks_{i_0}(q_\alpha)}, \cdots, \Ext_0^*(\EE_{i_j},\EE_{i_j})^{\fks_{i_j}(q_\alpha)}.\]
For each point $q\in U_\alpha$, we introduce a {\sl standard} $j$-simplex $f_q^{i_0,\cdots,i_j}\in N(\NL)$\footnote{This notation is due to the fact that the standard $j$-simplex only depends on the point $q$, and the vertices of $\alpha$.} whose
\medskip
\begin{itemize}
\item $0$-simplices are $\Ext_0^*(\EE_{i_0},\EE_{i_0})^{\fks_{i_0}(q)}, \cdots, \Ext_0^*(\EE_{i_j},\EE_{i_j})^{\fks_{i_j}(q)}$;
\medskip
\item $1$-simplices are of the form $\PP\Ad(\psi)\II$ 
\medskip
\item $2$-simplices are of the form $\PP\Ad(\psi)\HH\Ad(\psi)\II$;
\medskip
\item $k$-simplices are of the form $\PP\Ad(\psi)\HH\Ad(\psi)\cdots\Ad(\psi)\HH\Ad(\psi)\II$.
\end{itemize}

\medskip
\noindent In the above definition, we have omitted relevant superscripts and subscripts from $\II$'s, $\PP$'s, $\psi$'s and $\HH$'s, since it is clear from the context. For example, we have
\begin{align*}
f_q^{i_0,i_1} &= \PP^{B_{i_1}} \Ad(\psi_q^{i_0i_1}) \II^{b_{i_0}} \\
f_q^{i_0,i_1,i_2} &= \PP^{B_{i_2}} \Ad(\psi_q^{i_1i_2}) \HH^{B_{i_1}} \Ad(\psi_q^{i_0i_1}) \II^{b_{i_0}}
\end{align*}
where $b_{i_0}=\fks_{i_0}(q)$, $B_{i_1}=\II_* b_{i_1}$, and similarly for other subscripts.

\medskip
Analogous to Lemma~\ref{moving-lem} and~\ref{moving-homotopy}, we impose the following additional condition for the inductive construction.

\medskip
\begin{cond}~\label{moving-cond}
Let $0\leq j\leq k-1$, and $\alpha\in N_j$. Then for any point $q'\in U_\alpha$, we have that
\[ (f_\alpha)^{\fks_{i_0}(q')-\fks_{i_0}(q_\alpha)} \cong f_{q'}^{i_0,\cdots,i_j}\]
where the symbol $\cong$ denotes that the two simplices are homotopic in $N(\NL)_\bullet$.
\end{cond}

To prove the induction can be proceeded, we first show that liftings exist on the $k+1$-intersections.
We set
\begin{align*}
 N_{k}:=& \big\{ \alpha=(\alpha_0,\cdots,\alpha_k,q_\alpha)\in (N_{k-1})^{k+1}\times M\mid \\
 &\partial_s \alpha_t=\partial_{t-1} \alpha_s, \;\forall 0\leq s<t\leq k, \mbox{\; and \;} q_\alpha\in U_{\alpha_0}\cap\cdots\cap U_{\alpha_k}.\big\}
 \end{align*}

\medskip
\begin{prop}~\label{existence}
Assume the above induction hypothesis. For any index $\alpha\in N_k$ with vertices $i_0,\cdots,i_k$, set
\begin{align*}
\Delta_j&:= \fks_{i_0}(q_\alpha)-\fks_{i_0}(q_{\alpha_{j}}), \;\; \forall 1\leq j\leq k\\
\Delta_0&:=\fks_{i_1}(q_\alpha)-\fks_{i_1}(q_{\alpha_0}).
\end{align*}
Then the map
\[ \lambda_\alpha: \partial \Delta^k \ra (\NL)_\bullet\]
defined by the $k+1$ $(k-1)$-dimensional simplices $(f_{\alpha_0})^{\Delta_0}, (f_{\alpha_1})^{\Delta_1}, \cdots, (f_{\alpha_k})^{\Delta_k}$ admits a filling $f_\alpha\in (\NL)_k$. Furthermore, the filling $f_\alpha$ constructed is homotopic to the standard $k$-simplex $f_{q_\alpha}^{i_0,\cdots,i_k}$.
\end{prop}

\Pf. The proof is similar to that of Lemma~\ref{triple-lem}. Indeed, we first use Condition~\ref{moving-cond} to obtain a morphism
\[ \widetilde{\lambda_\alpha}: \Delta^1\times \partial \Delta^k \ra (\NL)_\bullet,\]
which restricts to $\lambda_\alpha$ on $0\times \partial \Delta^k$. On $1\times \partial \Delta^k$, the map $\widetilde{\lambda_\alpha}$ is given by the standard simplices
\[ f_{q_\alpha}^{i_0,\cdots,\widehat{i_j},\cdots,i_k}, \;\; 0\leq j\leq k.\]
Observe that the standard simplex $f_{q_\alpha}^{i_0,\cdots,i_k}$ is a filling for the map restricted to $1\times \partial \Delta^k$. Hence we obtain a further extension of the map $\widetilde{\lambda_\alpha}$ to a map
\[ (\Delta^1\times \partial \Delta^k) \bigcup_{\left\{1\right\}\times\partial\Delta^k} (\left\{1\right\} \times \Delta^k) \ra N(\NL)_\bullet.\]
Again, by Lemma~\ref{invertible-lem} the edges of the above diagram are weak equivalences in the $\infty$-category $N(\NL)_\bullet$. The proposition then follows from the homotopy lifting property of Kan complexes since the inclusion
\[ (\Delta^1\times \partial \Delta^k) \bigcup_{\left\{1\right\}\times\partial\Delta^k} (\left\{1\right\} \times \Delta^k)  \ra \Delta^1\times \Delta^k\]
is an acyclic cofibration.\qed

\medskip
To finish the inductive construction, we next prove Condition~\ref{moving-cond} for $j=k$. The proof is similar to that of Lemma~\ref{moving-homotopy}.

\medskip
\begin{prop}~\label{moving-prop}
For any point $q'$ sufficiently close to $q_\alpha$, we have that
\[ \big(f_{q_\alpha}^{i_0,\cdots,i_k}\big)^{\fks_{i_0}(q')-\fks_{i_0}(q_\alpha)} \cong f_{q'}^{i_0,\cdots,i_k}.\]
\end{prop}

\Pf. We need to construct a morphism
\[ \mu: \Delta^1\times \Delta^k \ra N(\NL)\]
which bounds $\big(f_{q_\alpha}^{i_0,\cdots,i_k}\big)^{\fks_{i_0}(q')-\fks_{i_0}(q_\alpha)}$ and $f_{q'}^{i_0,\cdots,i_k}$. This is equivalent to construct a morphism
\[ \chi: \Delta^1\times (\Delta^1)^{k-1} \ra \Hom_{\NL}\big(\Ext_0^*(\EE_{i_0},\EE_{i_0})^{\fks_{i_0}(q')}, \Ext_0^*(\EE_{i_k},\EE_{i_k})^{\fks_{i_k}(q')}\big)\]
with prescribed boundary condition induced from $\big(f_{q_\alpha}^{i_0,\cdots,i_k}\big)^{\fks_{i_0}(q')-\fks_{i_0}(q_\alpha)}$ and $f_{q'}^{i_0,\cdots,i_k}$. We first construct a family $\widetilde{\chi}$ parametrized by $([0,1])^{2k-1}$ by the following long composition:
\[\begin{tikzcd}[column sep=5.0em, row sep=4.0em]
\Ext_0^*(\EE_{i_0},\EE_{i_0})^{b_{i_0}'} \arrow{d}{\II^{b_{i_0}'}}\\
\Om_X^{0,*}(\EEnd_0(\EE_{i_i}))^{B_{i_0}'}\arrow{r}{\Ad(e^{t_1\eta_1})\Ad(\psi_{q'}^{i_0i_1})}& \big(\Om_X^{0,*}(\EEnd_0(\EE_{i_1}))^{B_{i_1}}\ot \Om_{[0,1]}^*\big)^{\theta_1}\arrow{ld}[description]{[\HH^{B_{i_1}}]^{\theta_1}}\\
\Om_X^{0,*}(\EEnd_0(\EE_{i_1}))^{B_{i_1}'}\ot \Om_{([0,1])^2}^*\arrow{r}{\Ad(e^{t_2\eta_2})\Ad(\psi_{q'}^{i_1i_2})}& \big(\Om_X^{0,*}(\EEnd_0(\EE_{i_2}))^{B_{i_2}}\ot \Om_{([0,1])^3}^*\big)^{\theta_2}\arrow{ld}[description]{[\HH^{B_{i_2}}]^{\theta_2}}\\
\cdots \cdots &\cdots \cdots\arrow{ld}[description]{[\HH^{B_{i_{k-1}}}]^{\theta_{k-1}}}\\
\Om_X^{0,*}(\EEnd_0(\EE_{i_{k-1}}))^{B_{i_{k-1}}'}\ot \Om_{([0,1])^{2k-2}}^*\arrow{r}{\Ad(e^{t_k\eta_k})\Ad(\psi_{q'}^{i_{k-1}i_k})}& \big(\Om_X^{0,*}(\EEnd_0(\EE_{i_k}))^{B_{i_k}}\ot \Om_{([0,1])^{2k-1}}^*\big)^{\theta_k}\arrow{d}[description]{[\PP^{B_{i_k}}]^{\theta_k}}\\
& \Ext_0^*(\EE_{i_k},\EE_{i_k})^{b_{i_k}'}\ot \Om_{([0,1])^{2k-1}}^*
\end{tikzcd}\]
This total composition gives us a morphism 
\[ \widetilde{\chi}: (\Delta^1)^{2k-1} \ra \Hom_{\NL}\big(\Ext_0^*(\EE_{i_0},\EE_{i_0})^{\fks_{i_0}(q')}, \Ext_0^*(\EE_{i_k},\EE_{i_k})^{\fks_{i_k}(q')}\big).\]
We define 
\[ \chi:=\widetilde{\chi}|_{t_1=t_2=\cdots=t_k}\]
where the coordinates $t_1, t_2, \cdots, t_k$ are as shown in the above diagram. It is straightforward to verify all the boundary conditions of $\chi$.\qed

\medskip
\begin{cor}
Condition~\ref{moving-cond} holds for $j=k$.
\end{cor}

\Pf. For each $\alpha\in N_k$, we can choose a lifting $f_\alpha$ of $\lambda_\alpha$ as in Proposition~\ref{existence}. Note that by construction, we have that
\[ f_\alpha \cong f_{q_\alpha}^{i_0,\cdots,i_k}.\]
Furthermore, by Proposition~\ref{moving-prop},  we can choose a small enough neighborhood $q_\alpha\in U_\alpha\subset U_{\alpha_0\cdots\alpha_k}$ so that for all points $q'\in U_\alpha$, we have that
\[ \big(f_{q_\alpha}^{i_0,\cdots,i_k}\big)^{\fks_{i_0}(q')-\fks_{i_0}(q_\alpha)} \cong f_{q'}^{i_0,\cdots,i_k}.\]
It follows that
\[ \big( f_\alpha \big)^{\fks_{i_0}(q')-\fks_{i_0}(q_\alpha)} \cong f_{q'}^{i_0,\cdots,i_k},\]
which proves Condition~\ref{moving-cond} for the case $j=k$.\qed

\medskip
Finally, we define the $k$-th morphism $\Phi_k: N_k\ra N(\LS)_k$ by putting
\[ \Phi_k(\alpha):= \cS(f_\alpha)\]
where $\cS(f_\alpha)$ is defined over appropriate open subsets
\begin{align*}
V_\alpha&\subset V_{i_0},\\
V_\alpha^\dagger &\subset V_{i_1},\\
\cdots &\\
V_{\alpha}^{(\dagger)^k} &\subset V_{i_k}
\end{align*}
such that
\[ U_\alpha=\fks_{i_0}^{-1}(V_\alpha)=\cdots=\fks_{i_k}^{-1}(V_\alpha^{(\dagger)^k}).\]

\subsection{Proof of Corollary~\ref{main-cor}} The existence of global sheaf of algebras $\OO_M^\bullet$ and its properties follow immediately from Proposition~\ref{derived-str-prop}. In the following we sketch a proof that the sheaf $\OO_M^\bullet$ is canonical, independent of all choices made in the construction of the homotopy $L_\infty$ enhancement $\Phi$. For this, it is convenient to introduce the notion of a refinement of a $L_\infty$ enhancement.

\begin{defi}
Let $\Phi$ be a homotopy $L_\infty$ enhancement of an analytic space $M$, subjected to a hypercovering $(N_\bullet,\UU)$. A refinement of the hypercovering is another hypercovering $(N_\bullet', \UU')$ such that there is a morphism of simplicial sets $\mu: N_\bullet' \ra N_\bullet$, and
\[ U_\alpha \subset U_{\mu(\alpha)}, \;\forall \alpha\in N_k', \;\forall k.\]
A refinement of $\Phi$ is another homotopy $L_\infty$ enhancement $\Phi'$, subjected to a refinement hypercovering $(N_\bullet', \UU')$ such that
\[ \Phi'(\alpha) \in \Phi(\mu(\alpha))\mid_{U_\alpha},\; \forall \alpha\in N_k', \;\forall k.\]
\end{defi}

If $\Phi'$ is a refinement of $\Phi$, then it is immediate that
\[ \underline{H}^*(M,\Phi') \cong \underline{H}^*(M,\Phi),\]
because local pieces of $\underline{H}^*(M,\Phi')$ are restrictions of the sheaf $\underline{H}^*(M,\Phi)$, and transition maps of $\underline{H}^*(M,\Phi')$ are restrictions of the transition maps of $\underline{H}^*(M,\Phi)$. From the construction of $\Phi$, different choices of $U$'s and $V$'s amount to only refinements of $\Phi$. Hence the cohomology sheaf $\OO_M^\bullet$ is independent of these choices.

Next, we show that the construction of $\OO_M^\bullet$ is also independent of the metric data chosen in order to perform homological perturbations. Let us analyze this at a point $q\in M$. The metrics data $g$ and $g'$ give rise to two different $L_\infty$ structures $\Phi_g$ and $\Phi_{g'}$ resulted from two sets of perturbation data $(i,p,h)$ and $(i',p',h')$. Consider the following composition
\[\begin{CD}
\big(\Om_X^{0,*}(\EEnd_0(E)), \db_q\big) @>\Ad>> \big(\Om_X^{0,*}(\EEnd_0(E)), \db_q'\big)\\
@A I AA       @VV P' V \\
\Ext^*_0(\EE_q,\EE_q)  @.   \Ext^*_0(\EE_q,\EE_q)'.
\end{CD}\]
Again, with respect to the $W^{l,2}$-norm ($l>>0$), all three morphisms are bounded, and hence their composition is a bounded $L_\infty$ morphisms. Applying the functor $\cS$ gives a morphism between the associated $L_\infty$ spaces in a neighborhood of the point $q$.

We need to show that these locally define morphisms glue on intersections. Let $U_i$ and $U_j$ be two open subsets, and let $p\in U_{ij}$. As before, set $b_i=\mathfrak{s}_i(p)$, $b_j=\mathfrak{s}_j(p)$, and similarly for $b_i'$, $b_j'$. Using Lemma~\ref{moving-lem}, we have
\begin{align*}
[P_j'\Ad_{q_i',q_j'} I_i']^{b_i'}\circ [P'_i \Ad_{q_i,q_i'} I_i]^{b_i}&\cong (P'_j)^{B_j'}\Ad_{B_i',B_j'}(I'_i)^{b_i'}\circ (P'_i)^{B_i'}\Ad_{B_i,B_i'} I^{b_i}\\
&\cong (P'_j)^{B_j'}\Ad_{B_i',B_j'}\Ad_{B_i,B_i'} I^{b_i}\\
&\cong (P'_j)^{B_j'}\Ad_{B_i,B_j'} I^{b_i},
\end{align*}
and the other composition gives
\begin{align*}
[P_j'\Ad_{q_j,q_j'} I_j]^{b_j}\circ [P_j \Ad_{q_i,q_j} I_i]^{b_i}&\cong (P'_j)^{B_j'}\Ad_{B_j,B_j'}(I_j)^{b_j}\circ (P_j)^{B_j}\Ad_{B_i,B_j} I^{b_i}\\
&\cong (P'_j)^{B_j'}\Ad_{B_j,B_j'}\Ad_{B_i,B_j} I^{b_i}\\
&\cong (P'_j)^{B_j'}\Ad_{B_i,B_j'}I^{b_i}.
\end{align*}
Thus the locally defined morphisms are compatible with the gluing data, from which we obtain a global morphism
\[ \Lambda:  \underline{H}^*(M,\Phi_g) \cong \underline{H}^*(M,\Phi_{g'}).\]
Since $\Lambda$ is locally induced by $\cS([P'_i \Ad_{q_i,q_i'} I_i])$, it is a quasi-isomorphism because all morphisms in $[P'_i \Ad_{q_i,q_i'} I_i]$ are. This completes the proof of Corollary~\ref{main-cor}.

\appendix
\section{Homological transfer with norms}~\label{ht-sec}

Let $E$ be a smooth complex vector bundle over a smooth projective variety $X$.  Throughout the section we fix an even positive integer $l$ such that $l>>0$ so that the multiplication map on $\Omega_X^{0,*}(\EEnd(E))$ extends to a bounded morphism on the corresponding order $l$ Sobolev completion. For example, any $l>\dim_\C X$ is fine. We denote the Sobolev completion of order $l$ by a subscript $l$. 
\medskip

\subsection{A-infinity case}~\label{a-pert-subsec}
Let $\EE$ be a holomorphic structure on $E$. By choosing a K\" ahler metric on $X$ and a Hermitian metric on the underlying smooth bundle $E$, one can define, and use the Hodge Laplacian operator $\Delta$ to realize the cohomology groups $\Ext^*(\EE,\EE)$ as the subspace of $\Omega_X^{0,*}(\EEnd(\EE))$ consisting of harmonic elements. Moreover, there is a homotopy retraction
\medskip
\begin{align}~\label{per-data-equ}
\begin{split}
i: \;\;\;& \Ext^*(\EE,\EE) \hookrightarrow \Omega_X^{0,*}(\EEnd(\EE)),\\
p: \;\;\;& \Omega_X^{0,*}(\EEnd(\EE)) \twoheadrightarrow \Ext^*(\EE,\EE),\\
h :\;\;\; & \Omega_X^{0,*}(\EEnd(\EE)) \ra \Omega_X^{0,*}(\EEnd(\EE)).
\end{split}
\end{align}
Here $h:=-\ad(\db_\EE^*)\circ G$, and $G$ is the Green's operator defined as
\[ G:=\int_0^\infty  e^{-t\Delta}dt.\] 
Note that the Green's operator is of order $-2$, which implies that $h$ is of order $-1$. Thus in particular, $h$ is a bounded operator on the Sobolev completion. These operators satisfy 
\begin{align*}
p\circ i &=\id\\
i\circ p &=\id+[\ad(\db_\EE),h],
\end{align*}
realizing $i$ and $p$ as homotopy inverses. Using $(i,p,h)$ we can transfer the differential graded algebra structure on $\Omega_X^{0,*}(\EEnd(\EE))$ to an $A_\infty$ algebra structure on $\Ext^*(\EE,\EE)$ via the planar tree formula (see for example~\cite{KS} and~\cite{Markl}). Recall that the Kontsevich-Soibelman's planar tree formula is defined by
\[ m_k := \sum_{T\in \OO(k)} \pm m_{k,T}\]
where $\OO(k)$ denotes the set of isomorphism classes of binary planar rooted trees with $k$ leaves. For each multi-linear map $m_{k,T}$ is the operadic composition by putting the product map of the algebra $\Omega_X^{0,*}(\EEnd(\EE))$ on internal vertices, the homotopy $h$ on internal edges, and the harmonic projection $p$ on the root of $T$. The trees involved in $m_2$ and $m_3$ are illustrated in the following picture.

\[\begin{tikzpicture}
  \tikzset{VertexStyle/.style = {
                                 text           = black,
                                 inner sep      = 0pt,
                                 outer sep      = 0pt,
                                 }}
  \tikzset{EdgeStyle/.style   = {thick,
                                 double          = black,
                                 double distance = 1pt
                               }}
  \tikzset{LabelStyle/.style =   {draw,
                                  shape = circle,
                                  inner sep      = .5pt,
                                 outer sep      = 0pt,
                                  fill           = white,
                                  text           = black,
                                  radius   = 2 pt}}
     \node[VertexStyle](i1) at (-10, 2) {};
     \node[VertexStyle](i2) at (-8, 2) {};
     \node[VertexStyle](m) at (-9, 1) {};
     \node[VertexStyle](r) at (-9, 0) {};     
     \node (m2) at (-9, -1) {$m_2$};
     \node[VertexStyle](j1) at (-6, 2) {};
          \node[VertexStyle](j2) at (-5, 2) {};
     \node[VertexStyle](j3) at (-4, 2) {};
     \node[VertexStyle](j4) at (-5.5,1) {};
          \node[VertexStyle](j5) at (-5, 0) {};
               \node[VertexStyle](j6) at (-5, -1) {};
\node[VertexStyle](k1) at (-2, 2) {};
          \node[VertexStyle](k2) at (-1, 2) {};
     \node[VertexStyle](k3) at (0, 2) {};
     \node[VertexStyle](k4) at (-0.5,1) {};
          \node[VertexStyle](k5) at (-1, 0) {};
               \node[VertexStyle](k6) at (-1, -1) {};
\node[VertexStyle](m3) at (-3,-2) {$m_3$};

     \draw[EdgeStyle](i1) to node[LabelStyle]{i} (m);
     \draw[EdgeStyle](i2) to node[LabelStyle]{i} (m);
     \draw[EdgeStyle](m) to node[LabelStyle]{p} (r);
     \draw[EdgeStyle](j1) to node[LabelStyle]{i} (j4);
          \draw[EdgeStyle](j2) to node[LabelStyle]{i} (j4);
               \draw[EdgeStyle](j3) to node[LabelStyle]{i} (j5);
     \draw[EdgeStyle](j4) to node[LabelStyle]{h} (j5);
     \draw[EdgeStyle](j5) to node[LabelStyle]{p} (j6);
\draw[EdgeStyle](k1) to node[LabelStyle]{i} (k5);
          \draw[EdgeStyle](k2) to node[LabelStyle]{i} (k4);
               \draw[EdgeStyle](k3) to node[LabelStyle]{i} (k4);
     \draw[EdgeStyle](k4) to node[LabelStyle]{h} (k5);
     \draw[EdgeStyle](k5) to node[LabelStyle]{p} (k6);

\end{tikzpicture}\]

The following lemma was obtained in~\cite{Fukaya}. We provide a more detailed proof. The norm used here is the Sobolev $W^{l, 2}$-norm which shall be denoted by $||\bullet||_{W^{l,2}}$.

\begin{lem}\label{norm-m-lem}
There exists a constant $C>0$, independent of $k$, such that
\[ ||m_k||_{W^{l,2}}\leq C^k.\]
\end{lem}

\Pf. Since both the product operator $\circ$, and the homotopy operator $h$ are bounded, there exist $D>1$ such that $||\circ||_{W^{l,2}}\leq D$ and $||h||_{W^{l,2}}\leq D$. On a binary planar tree with $k$ leaves we have $k-1$ internal vertices and $k-2$ internal edges, thus we get
\[ ||m_k||_{W^{l,2}} \leq \sum_{T\in \OO(k)} D^{k-1}\cdot D^{k-2}\]
as the projection operator $p$ on the root has norm $1$. Now the number of planar binary trees with $k$ leaves is given by the Catalan numbers $\frac{[2(k-1)]!}{(k-1)!k!}$, which is well-known to be bounded by, say $4^{k-1}$. This implies that
\[ ||m_k||_{W^{l,2}} \leq 4^{k-1}D^{k-1}D^{k-2} \leq (4D^2)^k.\]
Thus we may choose any $C\geq 4D^2$.
\ed

The lemma implies that, for $b\in \Ext^*(\EE,\EE)$ such that $||b||_{W^{l,2}}$ small enough, the infinite series
\[ m^b_k(\alpha_1,\cdots,\alpha_k):=\sum_{j_0\geq 0,\cdots,j_k\geq 0} I_{j_0+\cdots+j_k+k}(b^{j_0},\alpha_1,\cdots,\alpha_k,b^{j_k})\]
is absolutely convergent in the $W^{l,2}$-norm.

The homotopy retraction data $(i,p,h)$ also induces an $A_\infty$ homomorphism
\[ I: \Ext^*(\EE,\EE) \ra \Omega_X^{0,*}(\EEnd(\EE)).\]
The tree formula of this morphism is given by
\[ I_k:= \sum_{T\in \OO(k)} \pm I_{k,T}\]
where $I_{k,T}$ is defined in the same way as in the case of $m_{k,T}$ except on the root of $T$, instead of putting $p$, one puts $h$. The following pictures are trees in $I_2$ and $I_3$.

\[\begin{tikzpicture}
  \tikzset{VertexStyle/.style = {
                                 text           = black,
                                 inner sep      = 0pt,
                                 outer sep      = 0pt,
                                 }}
  \tikzset{EdgeStyle/.style   = {thick,
                                 double          = black,
                                 double distance = 1pt
                               }}
  \tikzset{LabelStyle/.style =   {draw,
                                  shape = circle,
                                  inner sep      = .5pt,
                                 outer sep      = 0pt,
                                  fill           = white,
                                  text           = black,
                                  radius   = 2 pt}}
     \node[VertexStyle](i1) at (-10, 2) {};
     \node[VertexStyle](i2) at (-8, 2) {};
     \node[VertexStyle](m) at (-9, 1) {};
     \node[VertexStyle](r) at (-9, 0) {};     
     \node (I2) at (-9, -1) {$I_2$};
     \node[VertexStyle](j1) at (-6, 2) {};
          \node[VertexStyle](j2) at (-5, 2) {};
     \node[VertexStyle](j3) at (-4, 2) {};
     \node[VertexStyle](j4) at (-5.5,1) {};
          \node[VertexStyle](j5) at (-5, 0) {};
               \node[VertexStyle](j6) at (-5, -1) {};
\node[VertexStyle](k1) at (-2, 2) {};
          \node[VertexStyle](k2) at (-1, 2) {};
     \node[VertexStyle](k3) at (0, 2) {};
     \node[VertexStyle](k4) at (-0.5,1) {};
          \node[VertexStyle](k5) at (-1, 0) {};
               \node[VertexStyle](k6) at (-1, -1) {};
\node[VertexStyle](I3) at (-3,-2) {$I_3$};

     \draw[EdgeStyle](i1) to node[LabelStyle]{i} (m);
     \draw[EdgeStyle](i2) to node[LabelStyle]{i} (m);
     \draw[EdgeStyle](m) to node[LabelStyle]{h} (r);
     \draw[EdgeStyle](j1) to node[LabelStyle]{i} (j4);
          \draw[EdgeStyle](j2) to node[LabelStyle]{i} (j4);
               \draw[EdgeStyle](j3) to node[LabelStyle]{i} (j5);
     \draw[EdgeStyle](j4) to node[LabelStyle]{h} (j5);
     \draw[EdgeStyle](j5) to node[LabelStyle]{h} (j6);
\draw[EdgeStyle](k1) to node[LabelStyle]{i} (k5);
          \draw[EdgeStyle](k2) to node[LabelStyle]{i} (k4);
               \draw[EdgeStyle](k3) to node[LabelStyle]{i} (k4);
     \draw[EdgeStyle](k4) to node[LabelStyle]{h} (k5);
     \draw[EdgeStyle](k5) to node[LabelStyle]{h} (k6);

\end{tikzpicture}\]

The proof of Lemma~\ref{norm-m-lem} also yields the following boundedness property of the morphisms $I_k$.

\medskip
\begin{lem}~\label{norm-i-lem}
There exists a constant $C>0$,  independent of $k$, such that
\[ ||I_k||_{W^{l,2}}\leq C^k.\]
\end{lem}

The lemma implies that, for $b\in \Ext^*(\EE,\EE)$ such that $||b||_{W^{l,2}}$ small enough, the infinite series
\[ I^b_k(\alpha_1,\cdots,\alpha_k):=\sum_{j_0\geq 0,\cdots,j_k\geq 0} I_{j_0+\cdots+j_k+k}(b^{j_0},\alpha_1,\cdots,\alpha_k,b^{j_k})\]
is absolutely convergent in the $W^{l,2}$-norm.

In the following we set 
\[ \kappa(b):=\sum_{k\geq 1} m_k(b^k), \mbox{\;\; and \;\; } B:=\sum_{k\geq 1} I_k(b^k)\]

\begin{lem}~\label{smooth-lem}
The section $I_k^b(\alpha_1,\cdots,\alpha_k)$ is smooth, i.e. it lies inside $\Om_X^{0,*}(\EEnd(\EE))$.
\end{lem}

\Pf. We first prove that $B$ is smooth. Using that $I$ is an $A_\infty$ homomorphism, we have
\[ \db B = I_1^b(\kappa(b))+ B\cdot B.\]
Since by construction the morphism $I_k$ has image in the kernel of $\db^*$, we get
\[ \db^*\db B= \db^* (B\cdot B).\]
But $B$ itself also lies in the image of $\db^*$, which implies that the above equation is equivalent to
\[ \Delta B= \db^*(B\cdot B).\]
Now the smoothness of B follows from the ellipticity of $\Delta$ together with the {\sl a priori} convergence in $W^{l,2}$-norm.

We prove the lemma by induction on $k\geq 1$. 
If $k=1$, that $I$ forms an $A_\infty$ homomorphism implies that
\[ \db I^b_1(\alpha)= I^b_2(\kappa(b),\alpha)+I^b_2(\alpha,\kappa(b))+I^b_1(m^b_1(\alpha))+ B\cdot I^b_1(\alpha)+I^b_1(\alpha) \cdot B .\]
Note that . Thus, applying $\db^*$ to the above equation, we get
\[ \db^*\db \big(I_1^b(\alpha)\big) +\db^*\big([B,I_1^b(\alpha)]\big)=0.\]
Since we also have $\db^* I_1^b(\alpha)=0$, this is equivalent to
\[ \Delta \big(I_1^b(\alpha)\big) + \db^*\big([B,I_1^b(\alpha)]\big)=0.\]
The regularity of $I_1^b(\alpha)$ follows from the ellipticity and the {\sl a priori} convergence.

Assume that for all $j<k$, the morphism $I_j^b$ has smooth image. Again, $I$ being an $A_\infty$ homomorphism implies that
\[ \db I^b_k(\alpha_1,\cdots,\alpha_k)\in [B,I_k^b(\alpha_1,\cdots,\alpha_k)]+\sum_{i+j=k} I_i^b(\alpha_1,\cdots,\alpha_i)\cdot I_j^b(\alpha_{i+1},\cdots,\alpha_k)+ \ker (\db^*).\]
Applying $\db^*$ yields
\[ \Delta I^b_k(\alpha_1,\cdots,\alpha_k)=\db^*\big([B,I_k^b(\alpha_1,\cdots,\alpha_k)]+\sum_{i+j=k} I_i^b(\alpha_1,\cdots,\alpha_i)\cdot I_j^b(\alpha_{i+1},\cdots,\alpha_k)\big).\]
Again, the regularity of $I_k^b(\alpha_1,\cdots,\alpha_k)$ follows from the ellipticity and the {\sl a priori} convergence.\ed

Following~\cite{Markl}, there exists another $A_\infty$ homomorphism 
$$P: \Omega_X^{0,*}(\EEnd(\EE)) \ra \Ext^*(\EE,\EE)$$
in the opposite direction. We also need to have certain bounded properties for $P$. For this it is essential to have an explicit formula of $P$, as in the case of $m$ and $I$. Fortunately, such formulas are provided in~\cite{Markl}. Namely, on a planar binary tree $T$ with $k$ leaves, one puts certain decoration $\tau$ on it by assigning white dots or black dots to the edges of $T$ (including the leaves). The decorations need to satisfy certain conditions, which we refer to ${\sl loc. cit.}$ Section $4$ for a precise treatment. We denote by $\QQ(k)$ the set of all such decorated binary planar rooted trees. Then one defines an operator 
\[ P_k:= \sum_{T\in\QQ(k)} \pm P_{k,T},\]
where $P_{k,T}$ is the operadic composition along $T$ on which we put the product map on vertices of $T$; the operator $i\circ p$ on an edge decorated by a black dot; the operator $h$ if the edge is decorated by a white dot; and finally on the root of $T$ put the operator $p$. The decorated trees involved in $P_2$ and $P_3$ are illustrated below.

\[\begin{tikzpicture}[font=\footnotesize]
  \tikzset{VertexStyle/.style = {
                                 text           = black,
                                 inner sep      = 0pt,
                                 outer sep      = 0pt,
                                 }}
  \tikzset{EdgeStyle/.style   = {thick,
                                 double          = black,
                                 double distance = 1pt
                               }}
  \tikzset{LabelStyle/.style =   {draw,
                                  shape = circle,
                                  inner sep      = .5pt,
                                 outer sep      = 0pt,
                                  fill           = white,
                                  text           = black,
                                  radius   = 2 pt}}
   \tikzset{LabelStyle2/.style =   {draw,
                                  shape = circle,
                                  inner sep      = .5pt,
                                 outer sep      = 0pt,
                                  fill           = black,
                                  text           = white,
                                  radius   = 2 pt}}
     \node[VertexStyle](i1) at (-10, 2) {};
     \node[VertexStyle](i2) at (-8, 2) {};
     \node[VertexStyle](m) at (-9, 1) {};
     \node[VertexStyle](r) at (-9, 0) {};   
      \node[VertexStyle](j1) at (-6, 2) {};
     \node[VertexStyle](j2) at (-4, 2) {};
     \node[VertexStyle](mm) at (-5, 1) {};
     \node[VertexStyle](rr) at (-5, 0) {};
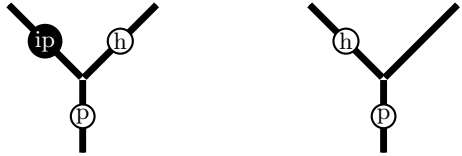
     \node[VertexStyle](P2) at (-7, -1) {trees in $P_2$};

     \draw[EdgeStyle](i1) to node[LabelStyle2]{ip} (m);
     \draw[EdgeStyle](i2) to node[LabelStyle]{h} (m);
     \draw[EdgeStyle](m) to node[LabelStyle]{p} (r);
     \draw[EdgeStyle](j1) to node[LabelStyle]{h} (mm);
     \draw[EdgeStyle](j2) to node{} (mm);
     \draw[EdgeStyle](mm) to node[LabelStyle]{p} (rr);

\end{tikzpicture}\]

\[\begin{tikzpicture}[font=\footnotesize]
  \tikzset{VertexStyle/.style = {
                                 text           = black,
                                 inner sep      = 0pt,
                                 outer sep      = 0pt,
                                 }}
  \tikzset{EdgeStyle/.style   = {thick,
                                 double          = black,
                                 double distance = 1pt
                               }}
  \tikzset{LabelStyle/.style =   {draw,
                                  shape = circle,
                                  inner sep      = .5pt,
                                 outer sep      = 0pt,
                                  fill           = white,
                                  text           = black,
                                  radius   = 2 pt}}
   \tikzset{LabelStyle2/.style =   {draw,
                                  shape = circle,
                                  inner sep      = .5pt,
                                 outer sep      = 0pt,
                                  fill           = black,
                                  text           = white,
                                  radius   = 2 pt}}
     \node[VertexStyle](1) at (-10, 2) {};
     \node[VertexStyle](2) at (-9, 2) {};
     \node[VertexStyle](3) at (-8, 2) {};
     \node[VertexStyle](4) at (-9.5, 1) {};
     \node[VertexStyle](5) at (-9, 0) {};          
    \node[VertexStyle](6) at (-9, -1) {};
    \node[VertexStyle](i1) at (-7, 2) {};
     \node[VertexStyle](i2) at (-6, 2) {};
     \node[VertexStyle](i3) at (-5, 2) {};
     \node[VertexStyle](i4) at (-6.5, 1) {};
     \node[VertexStyle](i5) at (-6, 0) {};          
    \node[VertexStyle](i6) at (-6, -1) {};
    \node[VertexStyle](j1) at (-4, 2) {};
     \node[VertexStyle](j2) at (-3, 2) {};
     \node[VertexStyle](j3) at (-2, 2) {};
     \node[VertexStyle](j4) at (-3.5, 1) {};
     \node[VertexStyle](j5) at (-3, 0) {};          
    \node[VertexStyle](j6) at (-3, -1) {};
    \node[VertexStyle](k1) at (-1, 2) {};
     \node[VertexStyle](k2) at (0, 2) {};
     \node[VertexStyle](k3) at (1, 2) {};
     \node[VertexStyle](k4) at (-0.5, 1) {};
     \node[VertexStyle](k5) at (0, 0) {};          
    \node[VertexStyle](k6) at (0, -1) {};
    \node[VertexStyle](l1) at (2, 2) {};
     \node[VertexStyle](l2) at (3, 2) {};
     \node[VertexStyle](l3) at (4, 2) {};
     \node[VertexStyle](l4) at (2.5, 1) {};
     \node[VertexStyle](l5) at (3, 0) {};          
    \node[VertexStyle](l6) at (3, -1) {};
    
     \draw[EdgeStyle](1) to node[LabelStyle]{h} (4);
          \draw[EdgeStyle](2) to node{} (4);
               \draw[EdgeStyle](3) to node[LabelStyle]{h} (5);
                    \draw[EdgeStyle](4) to node[LabelStyle2]{ip} (5);
                         \draw[EdgeStyle](5) to node[LabelStyle]{p} (6);
\draw[EdgeStyle](i1) to node[LabelStyle2]{ip} (i4);
          \draw[EdgeStyle](i2) to node[LabelStyle]{h} (i4);
               \draw[EdgeStyle](i3) to node[LabelStyle]{h} (i5);
                    \draw[EdgeStyle](i4) to node[LabelStyle2]{ip} (i5);
                         \draw[EdgeStyle](i5) to node[LabelStyle]{p} (i6);
 \draw[EdgeStyle](j1) to node[LabelStyle]{h} (j4);
          \draw[EdgeStyle](j2) to node{} (j4);
               \draw[EdgeStyle](j3) to node{} (j5);
                    \draw[EdgeStyle](j4) to node[LabelStyle]{h} (j5);
                         \draw[EdgeStyle](j5) to node[LabelStyle]{p} (j6);
\draw[EdgeStyle](k1) to node[LabelStyle2]{ip} (k4);
          \draw[EdgeStyle](k2) to node[LabelStyle]{h} (k4);
               \draw[EdgeStyle](k3) to node{} (k5);
                    \draw[EdgeStyle](k4) to node[LabelStyle]{h} (k5);
                         \draw[EdgeStyle](k5) to node[LabelStyle]{p} (k6);

     \draw[EdgeStyle](l1) to node[LabelStyle2]{ip} (l4);
          \draw[EdgeStyle](l2) to node[LabelStyle2]{ip} (l4);
               \draw[EdgeStyle](l3) to node[LabelStyle]{h} (l5);
                    \draw[EdgeStyle](l4) to node[LabelStyle]{h} (l5);
                         \draw[EdgeStyle](l5) to node[LabelStyle]{p} (l6);

\end{tikzpicture}\]

\[\begin{tikzpicture}[font=\footnotesize]
  \tikzset{VertexStyle/.style = {
                                 text           = black,
                                 inner sep      = 0pt,
                                 outer sep      = 0pt,
                                 }}
  \tikzset{EdgeStyle/.style   = {thick,
                                 double          = black,
                                 double distance = 1pt
                               }}
  \tikzset{LabelStyle/.style =   {draw,
                                  shape = circle,
                                  inner sep      = .5pt,
                                 outer sep      = 0pt,
                                  fill           = white,
                                  text           = black,
                                  radius   = 2 pt}}
   \tikzset{LabelStyle2/.style =   {draw,
                                  shape = circle,
                                  inner sep      = .5pt,
                                 outer sep      = 0pt,
                                  fill           = black,
                                  text           = white,
                                  radius   = 2 pt}}
     \node[VertexStyle](1) at (-10, 2) {};
     \node[VertexStyle](2) at (-9, 2) {};
     \node[VertexStyle](3) at (-8, 2) {};
     \node[VertexStyle](4) at (-8.5, 1) {};
     \node[VertexStyle](5) at (-9, 0) {};          
    \node[VertexStyle](6) at (-9, -1) {};
    \node[VertexStyle](i1) at (-7, 2) {};
     \node[VertexStyle](i2) at (-6, 2) {};
     \node[VertexStyle](i3) at (-5, 2) {};
     \node[VertexStyle](i4) at (-5.5, 1) {};
     \node[VertexStyle](i5) at (-6, 0) {};          
    \node[VertexStyle](i6) at (-6, -1) {};
    \node[VertexStyle](j1) at (-4, 2) {};
     \node[VertexStyle](j2) at (-3, 2) {};
     \node[VertexStyle](j3) at (-2, 2) {};
     \node[VertexStyle](j4) at (-2.5, 1) {};
     \node[VertexStyle](j5) at (-3, 0) {};          
    \node[VertexStyle](j6) at (-3, -1) {};
    \node[VertexStyle](k1) at (-1, 2) {};
     \node[VertexStyle](k2) at (0, 2) {};
     \node[VertexStyle](k3) at (1, 2) {};
     \node[VertexStyle](k4) at (0.5, 1) {};
     \node[VertexStyle](k5) at (0, 0) {};          
    \node[VertexStyle](k6) at (0, -1) {};
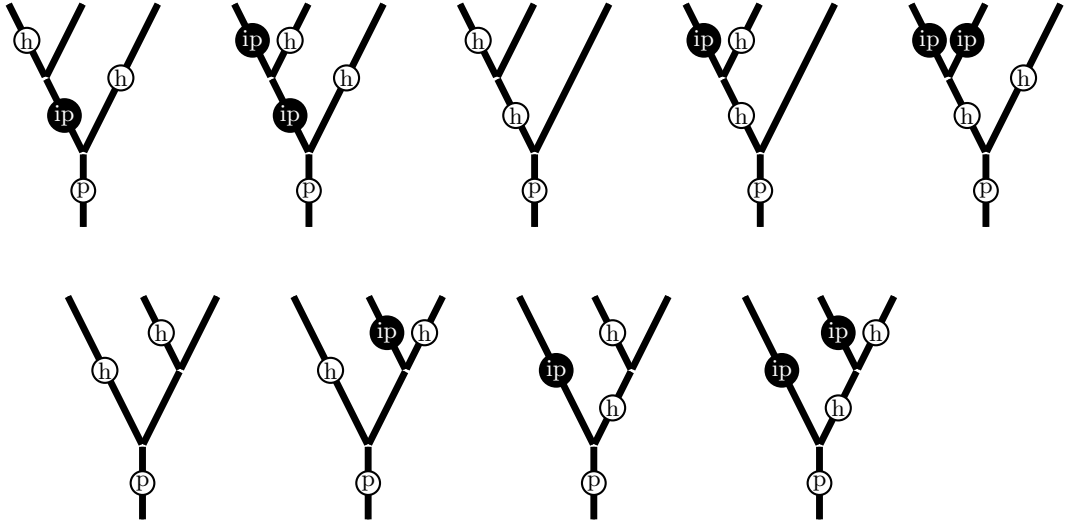
    \node(P3) at (-4,-2) {trees in $P_3$};
    
    \draw[EdgeStyle](k1) to node[LabelStyle2]{ip} (k5);
          \draw[EdgeStyle](k2) to node[LabelStyle2]{ip} (k4);
               \draw[EdgeStyle](k3) to node[LabelStyle]{h} (k4);
                    \draw[EdgeStyle](k4) to node[LabelStyle]{h} (k5);
                         \draw[EdgeStyle](k5) to node[LabelStyle]{p} (k6);
\draw[EdgeStyle](j1) to node[LabelStyle2]{ip} (j5);
          \draw[EdgeStyle](j2) to node[LabelStyle]{h} (j4);
               \draw[EdgeStyle](j3) to node{} (j4);
                    \draw[EdgeStyle](j4) to node[LabelStyle]{h} (j5);
                         \draw[EdgeStyle](j5) to node[LabelStyle]{p} (j6);
\draw[EdgeStyle](i1) to node[LabelStyle]{h} (i5);
          \draw[EdgeStyle](i2) to node[LabelStyle2]{ip} (i4);
               \draw[EdgeStyle](i3) to node[LabelStyle]{h} (i4);
                    \draw[EdgeStyle](i4) to node{} (i5);
                         \draw[EdgeStyle](i5) to node[LabelStyle]{p} (i6);    
\draw[EdgeStyle](1) to node[LabelStyle]{h} (5);
          \draw[EdgeStyle](2) to node[LabelStyle]{h} (4);
               \draw[EdgeStyle](3) to node{} (4);
                    \draw[EdgeStyle](4) to node{} (5);
                         \draw[EdgeStyle](5) to node[LabelStyle]{p} (6);

\end{tikzpicture}\]

The following result was proved in~\cite{Markl}, which gives a homotopy inverse of the morphism $I$.

\medskip
\begin{prop}~\label{formula-p}
The maps $P_k: \Omega_X^{0,*}(\EEnd(\EE))^{\otimes k} \ra \Ext^*(\EE,\EE)$ form an $A_\infty$ homomorphism.
\end{prop}

From the explicit formula of the $P_k$'s, we can deduce the following

\medskip
\begin{lem}~\label{norm-p-lem}
There exists a constant $C>0$, independent of $k$, such that
\[ ||P_k||_{W^{l,2}}\leq C^k.\]
\end{lem}

\Pf. We only observe the fact that, given a binary planar tree $T\in \OO(k)$, the number of edges including leaves is $2k-2$, which implies that the possible decorations on a fixed $T$ is bounded by $2^{2k-2}\leq 4^k$, which further implies that
\[ |\QQ(k)| \leq 4^k\cdot 4^k =16^k.\]
The rest of the proof is similar to that of Lemma~\ref{norm-m-lem}.\ed

The lemma implies that, for $B\in \Om_X^{0,*}(\EEnd(\EE))_l$ such that $||B||_{W^{l,2}}$ small enough, the infinite series
\[ P^B_k(\alpha_1,\cdots,\alpha_k):=\sum_{j_0\geq 0,\cdots,j_k\geq 0} P_{j_0+\cdots+j_k+k}(B^{j_0},\alpha_1,\cdots,\alpha_k,B^{j_k})\]
is absolutely convergent in the $W^{l,2}$-norm.

The $A_\infty$ homomorphisms $I$ and $P$ are inverse homotopy equivalences. To see this, first we note that the $A_\infty$ composition
\[ P\circ I =\id.\]
This follows from the interpretation of $I$ and $P$ as homological perturbation of differential graded coalgebras, see for example~\cite{Man}. 

Next we prove the other direction that $I\circ P$ is homotopic $\id$. This was proved in~\cite{Markl} with the {\sl coalgebra definition} of a homotopy between $A_\infty$ morphisms (~\cite[Section 2]{Markl}). For purposes of this paper, we need another definition of a homotopy between $A_\infty$ homomorphisms, due to Sullivan. 

\medskip
\begin{defi}~\label{a-inf-homotopy}
Two $A_\infty$ homomorphisms $f, g: A\ra B$ are homotopic if there exists an $A_\infty$ homomorphism
\[ H: A\ra B\ot\Omega_{\Delta^1}^*\]
where $\Omega_{\Delta^1}^*$ denotes the set of {\sl smooth} differential forms on the interval $\Delta^1=[0,1]$, such that
\[ H(0)=f \mbox{\;\;and\;\;} H(1)=g.\] 
Furthermore, if $B$ is a nuclear Frechet space, then the tensor product on the right hand side is taken in the category of nuclear Frechet spaces.
\end{defi}

Next we construct a homotopy in the sense of Definition~\ref{a-inf-homotopy} between $I\circ P$ and $\id$ on the differential graded algebra $\Omega_X^{0,*}(\EEnd(\EE))$.

We shall first construct a homotopy $H^\dagger$ parametrized by $t\in [0,\infty]$, using the heat kernel. Observe that the perturbation data $(i,p,h)$ in~\ref{per-data-equ} fits into the following family of perturbation data.
\begin{align}
\begin{split}
\id: \;\;\;& \Omega_X^{0,*}(\EEnd(\EE)) \rightarrow \Omega_X^{0,*}(\EEnd(\EE)),\\
K_t=e^{-t\Delta}: \;\;\;& \Omega_X^{0,*}(\EEnd(\EE)) \rightarrow \Om_X^{0,*}(\EEnd(\EE)),\\
h_t=\int_0^t -\ad(\db^*)e^{-s\Delta} ds :\;\;\; & \Omega_X^{0,*}(\EEnd(\EE)) \ra \Omega_X^{0,*}(\EEnd(\EE)).
\end{split}
\end{align}

Let $T\in \OO(s)$, and $(T_1,\tau_1)\in \QQ(l_1), \cdots, (T_s,\tau_s)\in \QQ(l_s)$ be $s$ decorated trees. We use the abbreviation $T(T_1,\cdots,T_s)$ to denote the planar rooted tree with $k=l_1+\cdots+l_s$ leaves obtained by joining the root of $T_j$ to the $j$-th leaf of $T$; and on it we put a decoration defined as follows:
\begin{itemize}
\item[(i)] on the edges of $T_1,\cdots, T_s$, we keep the decorations $\tau_1,\cdots, \tau_s$;
\item[(ii)] on the joining edges, we put black dots;
\item[(iii)] on the root and internal edges of $T$, we put white dots.
\end{itemize}
Denote by $\WW(k)$ the set of decorated trees obtained with the above procedure for various different possible $1\leq s\leq k$.

The constructions in~\cite{Markl} implies a family of $A_\infty$ homomorphism
\[ R: \Omega_X^{0,*}(\EEnd(\EE)) \ra \Omega_X^{0,*}(\EEnd(\EE))\ot C^\infty_{[0,\infty]}.\]
The $k$-th component of $R$ is defined by
\begin{equation}~\label{r-equ}
 R_k:= \sum_{T\in \WW(k)} R_{k,T}
\end{equation}
where the operator $R_{k,T}$ is the operadic composition by putting product map on vertices, the map $K_t$ on black dots, and the map $h_t$ on white dots. Explicitly, we have $ R_1= K_t$, and the maps $R_2$, $R_3$ are given by
\[\begin{tikzpicture}[font=\tiny]
  \tikzset{VertexStyle/.style = {
                                 text           = black,
                                 inner sep      = 0pt,
                                 outer sep      = 0pt,
                                 }}
  \tikzset{EdgeStyle/.style   = {thick,
                                 double          = black,
                                 double distance = 1pt
                               }}
  \tikzset{LabelStyle/.style =   {draw,
                                  shape = circle,
                                  inner sep      = .5pt,
                                 outer sep      = 0pt,
                                  fill           = white,
                                  text           = black,
                                  }}
   \tikzset{LabelStyle2/.style =   {draw,
                                  shape = circle,
                                  inner sep      = .5pt,
                                 outer sep      = 0pt,
                                  fill           = black,
                                  text           = white,
                                  }}
     \node[VertexStyle](i1) at (-10, 2) {};
     \node[VertexStyle](i2) at (-8, 2) {};
     \node[VertexStyle](m) at (-9, 1) {};
     \node[VertexStyle](r) at (-9, 0) {};   
      \node[VertexStyle](j1) at (-6, 2) {};
     \node[VertexStyle](j2) at (-4, 2) {};
     \node[VertexStyle](mm) at (-5, 1) {};
     \node[VertexStyle](rr) at (-5, 0) {};
      \node[VertexStyle](k1) at (-2, 2) {};
     \node[VertexStyle](k2) at (0, 2) {};
     \node[VertexStyle](mmm) at (-1, 1) {};
     \node[VertexStyle](rrr) at (-1, 0) {};
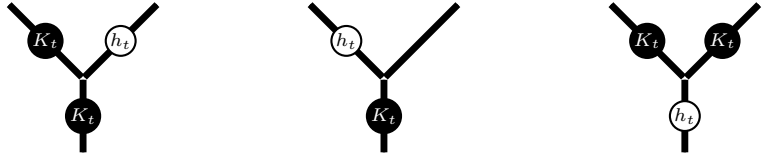
     \node[VertexStyle](P2) at (-5, -1) {\mbox{\normalsize Decorated trees in  $R_2$}};

     \draw[EdgeStyle](i1) to node[LabelStyle2]{$K_t$ } (m);
     \draw[EdgeStyle](i2) to node[LabelStyle]{$h_t$} (m);
     \draw[EdgeStyle](m) to node[LabelStyle2]{$K_t$} (r);
     \draw[EdgeStyle](j1) to node[LabelStyle]{$h_t$} (mm);
     \draw[EdgeStyle](j2) to node{} (mm);
     \draw[EdgeStyle](mm) to node[LabelStyle2]{$K_t$} (rr);
          \draw[EdgeStyle](k1) to node[LabelStyle2]{$K_t$} (mmm);
                    \draw[EdgeStyle](k2) to node[LabelStyle2]{$K_t$} (mmm);
     \draw[EdgeStyle](mmm) to node[LabelStyle]{$h_t$} (rrr);

\end{tikzpicture}\]

\[\begin{tikzpicture}[font=\tiny]
  \tikzset{VertexStyle/.style = {
                                 text           = black,
                                 inner sep      = 0pt,
                                 outer sep      = 0pt,
                                 }}
  \tikzset{EdgeStyle/.style   = {thick,
                                 double          = black,
                                 double distance = 1pt
                               }}
  \tikzset{LabelStyle/.style =   {draw,
                                  shape = circle,
                                  inner sep      = .5pt,
                                 outer sep      = 0pt,
                                  fill           = white,
                                  text           = black,
                                  radius   = 2 pt}}
   \tikzset{LabelStyle2/.style =   {draw,
                                  shape = circle,
                                  inner sep      = .5pt,
                                 outer sep      = 0pt,
                                  fill           = black,
                                  text           = white,
                                  radius   = 2 pt}}
     \node[VertexStyle](1) at (-10, 2) {};
     \node[VertexStyle](2) at (-9, 2) {};
     \node[VertexStyle](3) at (-8, 2) {};
     \node[VertexStyle](4) at (-9.5, 1) {};
     \node[VertexStyle](5) at (-9, 0) {};          
    \node[VertexStyle](6) at (-9, -1) {};
    \node[VertexStyle](i1) at (-7, 2) {};
     \node[VertexStyle](i2) at (-6, 2) {};
     \node[VertexStyle](i3) at (-5, 2) {};
     \node[VertexStyle](i4) at (-6.5, 1) {};
     \node[VertexStyle](i5) at (-6, 0) {};          
    \node[VertexStyle](i6) at (-6, -1) {};
    \node[VertexStyle](j1) at (-4, 2) {};
     \node[VertexStyle](j2) at (-3, 2) {};
     \node[VertexStyle](j3) at (-2, 2) {};
     \node[VertexStyle](j4) at (-3.5, 1) {};
     \node[VertexStyle](j5) at (-3, 0) {};          
    \node[VertexStyle](j6) at (-3, -1) {};
    \node[VertexStyle](k1) at (-1, 2) {};
     \node[VertexStyle](k2) at (0, 2) {};
     \node[VertexStyle](k3) at (1, 2) {};
     \node[VertexStyle](k4) at (-0.5, 1) {};
     \node[VertexStyle](k5) at (0, 0) {};          
    \node[VertexStyle](k6) at (0, -1) {};
    \node[VertexStyle](l1) at (2, 2) {};
     \node[VertexStyle](l2) at (3, 2) {};
     \node[VertexStyle](l3) at (4, 2) {};
     \node[VertexStyle](l4) at (2.5, 1) {};
     \node[VertexStyle](l5) at (3, 0) {};          
    \node[VertexStyle](l6) at (3, -1) {};
    
     \draw[EdgeStyle](1) to node[LabelStyle]{$h_t$} (4);
          \draw[EdgeStyle](2) to node{} (4);
               \draw[EdgeStyle](3) to node[LabelStyle]{$h_t$} (5);
                    \draw[EdgeStyle](4) to node[LabelStyle2]{$K_t$} (5);
                         \draw[EdgeStyle](5) to node[LabelStyle2]{$K_t$} (6);
\draw[EdgeStyle](i1) to node[LabelStyle2]{$K_t$} (i4);
          \draw[EdgeStyle](i2) to node[LabelStyle]{$h_t$} (i4);
               \draw[EdgeStyle](i3) to node[LabelStyle]{$h_t$} (i5);
                    \draw[EdgeStyle](i4) to node[LabelStyle2]{$K_t$} (i5);
                         \draw[EdgeStyle](i5) to node[LabelStyle2]{$K_t$} (i6);
 \draw[EdgeStyle](j1) to node[LabelStyle]{$h_t$} (j4);
          \draw[EdgeStyle](j2) to node{} (j4);
               \draw[EdgeStyle](j3) to node{} (j5);
                    \draw[EdgeStyle](j4) to node[LabelStyle]{$h_t$} (j5);
                         \draw[EdgeStyle](j5) to node[LabelStyle2]{$K_t$} (j6);
\draw[EdgeStyle](k1) to node[LabelStyle2]{$K_t$} (k4);
          \draw[EdgeStyle](k2) to node[LabelStyle]{$h_t$} (k4);
               \draw[EdgeStyle](k3) to node{} (k5);
                    \draw[EdgeStyle](k4) to node[LabelStyle]{$h_t$} (k5);
                         \draw[EdgeStyle](k5) to node[LabelStyle2]{$K_t$} (k6);

     \draw[EdgeStyle](l1) to node[LabelStyle2]{$K_t$} (l4);
          \draw[EdgeStyle](l2) to node[LabelStyle2]{$K_t$} (l4);
               \draw[EdgeStyle](l3) to node[LabelStyle]{$h_t$} (l5);
                    \draw[EdgeStyle](l4) to node[LabelStyle]{$h_t$} (l5);
                         \draw[EdgeStyle](l5) to node[LabelStyle2]{$K_t$} (l6);

\end{tikzpicture}\]

\[\begin{tikzpicture}[font=\tiny]
  \tikzset{VertexStyle/.style = {
                                 text           = black,
                                 inner sep      = 0pt,
                                 outer sep      = 0pt,
                                 }}
  \tikzset{EdgeStyle/.style   = {thick,
                                 double          = black,
                                 double distance = 1pt
                               }}
  \tikzset{LabelStyle/.style =   {draw,
                                  shape = circle,
                                  inner sep      = .5pt,
                                 outer sep      = 0pt,
                                  fill           = white,
                                  text           = black,
                                  radius   = 2 pt}}
   \tikzset{LabelStyle2/.style =   {draw,
                                  shape = circle,
                                  inner sep      = .5pt,
                                 outer sep      = 0pt,
                                  fill           = black,
                                  text           = white,
                                  radius   = 2 pt}}
     \node[VertexStyle](1) at (-10, 2) {};
     \node[VertexStyle](2) at (-9, 2) {};
     \node[VertexStyle](3) at (-8, 2) {};
     \node[VertexStyle](4) at (-8.5, 1) {};
     \node[VertexStyle](5) at (-9, 0) {};          
    \node[VertexStyle](6) at (-9, -1) {};
    \node[VertexStyle](i1) at (-7, 2) {};
     \node[VertexStyle](i2) at (-6, 2) {};
     \node[VertexStyle](i3) at (-5, 2) {};
     \node[VertexStyle](i4) at (-5.5, 1) {};
     \node[VertexStyle](i5) at (-6, 0) {};          
    \node[VertexStyle](i6) at (-6, -1) {};
    \node[VertexStyle](j1) at (-4, 2) {};
     \node[VertexStyle](j2) at (-3, 2) {};
     \node[VertexStyle](j3) at (-2, 2) {};
     \node[VertexStyle](j4) at (-2.5, 1) {};
     \node[VertexStyle](j5) at (-3, 0) {};          
    \node[VertexStyle](j6) at (-3, -1) {};
    \node[VertexStyle](k1) at (-1, 2) {};
     \node[VertexStyle](k2) at (0, 2) {};
     \node[VertexStyle](k3) at (1, 2) {};
     \node[VertexStyle](k4) at (0.5, 1) {};
     \node[VertexStyle](k5) at (0, 0) {};          
    \node[VertexStyle](k6) at (0, -1) {};
     \node[VertexStyle](R31) at (-5, -2) {\mbox{\normalsize Decorated trees in $R_3$ of the form $T(T_1)$}};

    \draw[EdgeStyle](k1) to node[LabelStyle2]{$K_t$} (k5);
          \draw[EdgeStyle](k2) to node[LabelStyle2]{$K_t$} (k4);
               \draw[EdgeStyle](k3) to node[LabelStyle]{$h_t$} (k4);
                    \draw[EdgeStyle](k4) to node[LabelStyle]{$h_t$} (k5);
                         \draw[EdgeStyle](k5) to node[LabelStyle2]{$K_t$} (k6);
\draw[EdgeStyle](j1) to node[LabelStyle2]{$K_t$} (j5);
          \draw[EdgeStyle](j2) to node[LabelStyle]{$h_t$} (j4);
               \draw[EdgeStyle](j3) to node{} (j4);
                    \draw[EdgeStyle](j4) to node[LabelStyle]{$h_t$} (j5);
                         \draw[EdgeStyle](j5) to node[LabelStyle2]{$K_t$} (j6);
\draw[EdgeStyle](i1) to node[LabelStyle]{$h_t$} (i5);
          \draw[EdgeStyle](i2) to node[LabelStyle2]{$K_t$} (i4);
               \draw[EdgeStyle](i3) to node[LabelStyle]{$h_t$} (i4);
                    \draw[EdgeStyle](i4) to node{} (i5);
                         \draw[EdgeStyle](i5) to node[LabelStyle2]{$K_t$} (i6);    
\draw[EdgeStyle](1) to node[LabelStyle]{$h_t$} (5);
          \draw[EdgeStyle](2) to node[LabelStyle]{$h_t$} (4);
               \draw[EdgeStyle](3) to node{} (4);
                    \draw[EdgeStyle](4) to node{} (5);
                         \draw[EdgeStyle](5) to node[LabelStyle2]{$K_t$} (6);
\end{tikzpicture}\]

\[\begin{tikzpicture}[font=\tiny]
  \tikzset{VertexStyle/.style = {
                                 text           = black,
                                 inner sep      = 0pt,
                                 outer sep      = 0pt,
                                 }}
  \tikzset{EdgeStyle/.style   = {thick,
                                 double          = black,
                                 double distance = 1pt
                               }}
  \tikzset{LabelStyle/.style =   {draw,
                                  shape = circle,
                                  inner sep      = .5pt,
                                 outer sep      = 0pt,
                                  fill           = white,
                                  text           = black,
                                  radius   = 2 pt}}
   \tikzset{LabelStyle2/.style =   {draw,
                                  shape = circle,
                                  inner sep      = .5pt,
                                 outer sep      = 0pt,
                                  fill           = black,
                                  text           = white,
                                  radius   = 2 pt}}
     \node[VertexStyle](1) at (-10, 2) {};
     \node[VertexStyle](2) at (-9, 2) {};
     \node[VertexStyle](3) at (-8, 2) {};
     \node[VertexStyle](4) at (-9.5, 1) {};
     \node[VertexStyle](5) at (-9, 0) {};          
    \node[VertexStyle](6) at (-9, -1) {};
    \node[VertexStyle](i1) at (-7, 2) {};
     \node[VertexStyle](i2) at (-6, 2) {};
     \node[VertexStyle](i3) at (-5, 2) {};
     \node[VertexStyle](i4) at (-6.5, 1) {};
     \node[VertexStyle](i5) at (-6, 0) {};          
    \node[VertexStyle](i6) at (-6, -1) {};
     \node[VertexStyle](j1) at (-4, 2) {};
     \node[VertexStyle](j2) at (-3, 2) {};
     \node[VertexStyle](j3) at (-2, 2) {};
     \node[VertexStyle](j4) at (-2.5, 1) {};
     \node[VertexStyle](j5) at (-3, 0) {};          
    \node[VertexStyle](j6) at (-3, -1) {};
   \node[VertexStyle](k1) at (-1, 2) {};
     \node[VertexStyle](k2) at (0, 2) {};
     \node[VertexStyle](k3) at (1, 2) {};
     \node[VertexStyle](k4) at (0.5, 1) {};
     \node[VertexStyle](k5) at (0, 0) {};          
    \node[VertexStyle](k6) at (0, -1) {}; 
    \node[VertexStyle](R31) at (-4.5, -2) {\mbox{\normalsize Decorated trees in $R_3$ of the form $T(T_1,T_2)$}};
    
    \draw[EdgeStyle](1) to node[LabelStyle2]{$K_t$} (4);
          \draw[EdgeStyle](2) to node[LabelStyle]{$h_t$} (4);
               \draw[EdgeStyle](3) to node[LabelStyle2]{$K_t$} (5);
                    \draw[EdgeStyle](4) to node[LabelStyle2]{$K_t$} (5);
                         \draw[EdgeStyle](5) to node[LabelStyle]{$h_t$} (6);
                           \draw[EdgeStyle](i1) to node[LabelStyle]{$h_t$} (i4);
          \draw[EdgeStyle](i2) to node{} (i4);
               \draw[EdgeStyle](i3) to node[LabelStyle2]{$K_t$} (i5);
                    \draw[EdgeStyle](i4) to node[LabelStyle2]{$K_t$} (i5);
                         \draw[EdgeStyle](i5) to node[LabelStyle]{$h_t$} (i6);
  \draw[EdgeStyle](j1) to node[LabelStyle2]{$K_t$} (j5);
          \draw[EdgeStyle](j2) to node[LabelStyle2]{$K_t$} (j4);
               \draw[EdgeStyle](j3) to node[LabelStyle]{$h_t$} (j4);
                    \draw[EdgeStyle](j4) to node[LabelStyle2]{$K_t$} (j5);
                         \draw[EdgeStyle](j5) to node[LabelStyle]{$h_t$} (j6);
            \draw[EdgeStyle](k1) to node[LabelStyle2]{$K_t$} (k5);
          \draw[EdgeStyle](k2) to node[LabelStyle]{$h_t$} (k4);
               \draw[EdgeStyle](k3) to node{} (k4);
                    \draw[EdgeStyle](k4) to node[LabelStyle2]{$K_t$} (k5);
                         \draw[EdgeStyle](k5) to node[LabelStyle]{$h_t$} (k6);
\end{tikzpicture}\]

\[\begin{tikzpicture}[font=\tiny]
  \tikzset{VertexStyle/.style = {
                                 text           = black,
                                 inner sep      = 0pt,
                                 outer sep      = 0pt,
                                 }}
  \tikzset{EdgeStyle/.style   = {thick,
                                 double          = black,
                                 double distance = 1pt
                               }}
  \tikzset{LabelStyle/.style =   {draw,
                                  shape = circle,
                                  inner sep      = .5pt,
                                 outer sep      = 0pt,
                                  fill           = white,
                                  text           = black,
                                  radius   = 2 pt}}
   \tikzset{LabelStyle2/.style =   {draw,
                                  shape = circle,
                                  inner sep      = .5pt,
                                 outer sep      = 0pt,
                                  fill           = black,
                                  text           = white,
                                  radius   = 2 pt}}
     \node[VertexStyle](1) at (-10, 2) {};
     \node[VertexStyle](2) at (-9, 2) {};
     \node[VertexStyle](3) at (-8, 2) {};
     \node[VertexStyle](4) at (-9.5, 1) {};
     \node[VertexStyle](5) at (-9, 0) {};          
    \node[VertexStyle](6) at (-9, -1) {};
    
     \node[VertexStyle](i1) at (-6, 2) {};
     \node[VertexStyle](i2) at (-5, 2) {};
     \node[VertexStyle](i3) at (-4, 2) {};
     \node[VertexStyle](i4) at (-4.5, 1) {};
     \node[VertexStyle](i5) at (-5, 0) {};          
    \node[VertexStyle](i6) at (-5, -1) {};
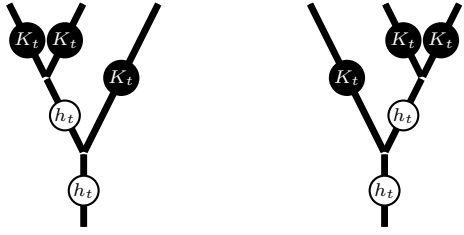
\node[VertexStyle](R31) at (-7, -2) {\mbox{\normalsize Decorated trees in $R_3$ of the form $T(T_1,T_2,T_3)$}};

\draw[EdgeStyle](1) to node[LabelStyle2]{$K_t$} (4);
          \draw[EdgeStyle](2) to node[LabelStyle2]{$K_t$} (4);
               \draw[EdgeStyle](3) to node[LabelStyle2]{$K_t$} (5);
                    \draw[EdgeStyle](4) to node[LabelStyle]{$h_t$} (5);
                         \draw[EdgeStyle](5) to node[LabelStyle]{$h_t$} (6);
 \draw[EdgeStyle](i1) to node[LabelStyle2]{$K_t$} (i5);
          \draw[EdgeStyle](i2) to node[LabelStyle2]{$K_t$} (i4);
               \draw[EdgeStyle](i3) to node[LabelStyle2]{$K_t$} (i4);
                    \draw[EdgeStyle](i4) to node[LabelStyle]{$h_t$} (i5);
                         \draw[EdgeStyle](i5) to node[LabelStyle]{$h_t$} (i6);                            
\end{tikzpicture}\]    
The conditions 
\[ K_0=\id \mbox{ \;\; and \;\; } h_0=0\]
imply that $R|_{t=0}=\id$. On the other end, the boundary conditions
\[ K_\infty=i\circ p \mbox{\;\; and \;\;} h_\infty= h\]
imply that $R|_{t=\infty}= I \circ P$.

Next we extend the family of homomorphisms $R$ to a homotopy between $\id$ and $I\circ P$. For this, we shall slightly modify the set $\WW(k)$ as follows. For each decorated tree $T$ in $\WW(k)$, we choose a black dot, and switch it to a third colored dot: say a blue dot. We denote by the set of all such decorated trees by $\VV(k)$. For each tree $T\in \VV(k)$, we associate an operator $S_{k,T}$ as the operadic composition along $T$ by putting the product map on vertices, $K_t$ on black dots, $h_t$ on white dots, and the operator $-\ad(\db^*)e^{-t\Delta}$ on the unique blue dot. Then we define a sequence $S$ of multilinear maps $S_k$ for $k\geq 1$ by formula
\begin{equation}~\label{s-equ}
 S_k:= \sum_{T\in \VV(k)} S_{k,T}.
\end{equation}

\begin{lem}
The morphism $H^\dagger:=R+S dt: \Omega_X^{0,*}(\EEnd(\EE)) \ra \Omega_X^{0,*}(\EEnd(\EE))\ot \Om_{[0,\infty]}^*$ is an $A_\infty$ homomorphism such that $H^\dagger\mid_{t=0}=\id$ and $H^\dagger\mid_{t=\infty}=I\circ P$.
\end{lem}

\Pf. Let us denote by
\[ \epsilon(R_{k,T})=\sum_{T'} S_{k,T'},\]
where the right hand side summation is over decorated trees $T'\in \VV(k)$ such that after switching the unique blue dot of $T'$ to a black dot, the resulting decorated tree is just $T$. Then using the identities
\begin{align}~\label{commutator}
\frac{d}{dt} K_t&=[\ad(\db),-\ad(\db^*)e^{-t\Delta}] \mbox{ and }\\
\frac{d}{dt} h_t&= -\ad(\db^*)e^{-t\Delta},
\end{align} 
we can deduce that
\begin{equation}~\label{eq-a-hom}
[\ad(\db), \epsilon(R_{k,T})]+\frac{d}{dt}R_{k,T}= \epsilon \big( [\ad(\db), R_{k,T}]\big).\end{equation} 
This follows from the observation that 
\begin{itemize}
\item[(1)] Bracketing with $\ad(\db)$ on a blue dot (which is assigned the operator $-\ad(\db^*)e^{-t\Delta}$) yields the derivative of $R_{k,T}$ on black dots, by the first formula in~\ref{commutator}.
\item[(2)] Bracketing with $\ad(\db)$ on a white dot in $R_{k,T}$ yields a black dot plus the identity map. This extra black dot turned into a blue dot, after the $\epsilon$ operation, which cancels exactly the derivative of $R_{k,T}$ on white dots by the second formula in~\ref{commutator}.
\end{itemize}
Since the maps $R_k$'s form an $A_\infty$ homomorphism, we have
\[ [\ad(\db), R_k] = \sum_{i,j} m_2(R_i\ot R_j) + \sum R_{k-1}(\id^l\ot m_2\ot \id^{k-l-2}).\]
Applying $\epsilon$ to this equation, and using identity~\ref{eq-a-hom}, we get
\[ [\ad(\db), S_k]+\frac{d}{dt}R_k=\sum_{i,j} m_2(S_i\ot R_j)+m_2(R_i\ot S_j) + \sum S_{k-1}(\id^l\ot m_2\ot \id^{k-l-2}).\]
This proves that $R+Sdt$ forms an $A_\infty$ homorphism. The boundary conditions have already been checked.\ed

\begin{lem}~\label{exp-decay}
Let $f(u):= \tan\frac{\pi u}{2}$. Then the pull-back $A_\infty$ homomorphism $$f^*H^\dagger: \Omega_X^{0,*}(\EEnd(\EE)) \ra \Omega_X^{0,*}(\EEnd(\EE))\ot \Om_{[0,1]}^*$$ is well-defined. It gives a homotopy (in the sense of Definition~\ref{a-inf-homotopy}) between $\id$ and $I\circ P$, which we denote by $H$.
\end{lem}

\Pf. The problem is that the derivatives of $f(u)$ blows up at $u=1$. However, they are all of polynomial order, while the heat kernel's derivative operators $(-\Delta)^me^{-t\Delta}$ for $m\geq 1$ are of exponential decay at infinity. Thus the pull-back operator $f^*H^\dagger$ remains smooth over $[0,1]$.\ed

\subsection{Boundedness of the homotopy $H$} In this subsection, we prove some boundedness property of the homotopy $H$ and its derivatives in the $t$-direction. To simplify the notations, in the following we set $A:=\Omega_X^{0,*}(\EEnd(\EE))$.

\begin{lem}
The operator $K_t=e^{-t\Delta}$ on $A$ extends to a bounded linear operator on $A_l$. The operator norm $||K_t||_{W^{l,2}}$ is uniformly bounded by a constant $C$ on the interval $[0,\infty]$. (The constant $C$ depends on $l$.)
\end{lem}

\Pf. By the spectral decomposition theorem, on the $L^2$-completion of $A$, the operator $e^{-t\Delta}$ acts on an eigen-section $s_\lambda$ by $e^{-t\lambda}$. In this case, it is clear that $||K_t||_{L^2}$ is uniformly bounded by $1$, since eigenvalues of $\Delta$ is positive.

For general $l$, Garding's inequality implies that there exists a constant $C>0$, depending on $l$, such that for any section $a\in A$, we have
\[ || a ||_{W^{l,2}} \leq C (||a ||_{L^2} + || \Delta^{l/2} a||_{L^2}).\]
Thus we get
\begin{align*}
|| e^{-t\Delta} a ||_{W^{l,2}} &\leq C( || e^{-t\Delta} a ||_{L^2}+ || \Delta^{l/2} e^{-t\Delta} a ||_{L^2})\\
&\leq C( || e^{-t\Delta} a||_{L^2} + || e^{-t\Delta} (\Delta^{l/2} a) ||_{L^2})\\
&\leq C( || a ||_{L^2}+ || (\Delta^{l/2} a) ||_{L^2})\\
&\leq C || a ||_{W^{l,2}}.
\end{align*}
It is clear from the proof that the constant $C$ only depends on $l$.\ed

\begin{lem}~\label{bound-2-lem}
Let $t>0$ be strictly positive. Then for any positive integer $j$, the operator $\Delta^j e^{-t\Delta}$ extends to a bounded linear operator on $A_l$. Furthermore, there exists a constant $C$, depending on $j$ and $l$, such that
\[ ||\Delta^j e^{-t\Delta}||_{W^{l,2}}\leq C/t^j.\]
\end{lem}

\Pf. On the $L^2$-completion of $A$, the operator $\Delta^j e^{-t\Delta}$ acts on an eigen-section $s_\lambda$ by $\lambda^j e^{-t\lambda}$. Observe that the function
\[ \lambda \mapsto \lambda^j e^{-t\lambda}\]
is uniformly bounded by its maximal value $\frac{j^j\cdot e^{-j}}{t^j}$, realized at $\lambda=j/t$. We conclude that the operator norm
\[ ||\Delta^j e^{-t\Delta}||_{L^2} \leq D/t^j\]
for some constant $D>0$ which depends on $j$.

In general, again by Garding's inequality we have
\begin{align*}
|| \Delta^j e^{-t\Delta} a ||_{W^{l,2}} &\leq E ( || \Delta^j e^{-t\Delta} a ||_{L^2}+||\Delta^{l/2+j} e^{-t\Delta} a ||_{L^2})\\
&\leq E\cdot (D/t^j) ( ||a||_{L^2} + || \Delta^{l/2} a ||_{L^2})\\
&\leq \frac{E\cdot D}{t^j} ||a||_{W^{l,2}}.
\end{align*}
The constant $C:= E\cdot D$ depends on both $j$ and $l$.\ed

\begin{lem}~\label{bound-3-lem}
Let $t>0$ be strictly positive. Then the operator $\ad(\db^*) e^{-t\Delta}$ extends to a bounded linear operator on $A_l$. Furthermore, there exists a constant $C>0$, depending on $l$, such that 
\[ || \ad(\db^*) e^{-t\Delta} ||_{W^{l,2}} \leq C/\sqrt{t}.\]
\end{lem}

\Pf. Again we first consider the $L^2$-norm. We have
\begin{align*}
||\ad(\db^*)e^{-t\Delta} a ||^2_{L^2}&=\langle \ad(\db^*) e^{-t\Delta} a, \ad(\db^*)e^{-t\Delta} a\rangle\\
&=\langle \ad(\db)\ad(db^*)e^{-t\Delta} a, e^{-t\Delta} a \rangle\\
&\leq \langle \Delta e^{-t\Delta}  a, e^{-t\Delta}  a\rangle\\
&\leq D/t ||a||_{L^2}^2    \mbox{\;\;\;(by the previous lemma).}
\end{align*}
Taking square root both sides, we get that
\[ || \ad(\db^*)e^{-t\Delta} ||_{L^2} \leq C/\sqrt{t}\]
with $C=\sqrt{D}$.

For general the bound in the case of $W^{l,2}$-norm, the proof is similar to the previous Lemma, using the commutativity $[\ad(\db^*),\Delta]=0$.\ed

\begin{cor}
The operator $h_t=\int_0^t -\ad(\db^*) e^{-s\Delta} ds$ extends to a bounded linear operator on $A_l$. Moreover, there exists a constant $C>0$, depending on $l$, which gives a uniform bound $||h_t||_{W^{l,2}} \leq C$.
\end{cor}

\Pf. By the previous Lemma, we have
\[ || h_t ||_{W^{l,2}} \leq C \cdot \int_0^t \frac{1}{\sqrt{s}} ds \leq 2C \sqrt{t}.\]
This implies the uniform boundedness on any bounded interval $[0,N]$. Thus it remains to bound the limit at infinity. This is well-known: since the Green's operator $\int_0^\infty e^{-s\Delta}ds$ is of order $-2$, the limit operator $h_\infty$ is of order $-1$, hence bounded.\ed

\begin{cor}
There exists a constant $C>0$, independent of $k$, such that it gives a uniform bound of the form
\[ || R_k ||_{W^{l,2}} \leq C^k\]
on the interval $[0,\infty]$.
\end{cor}

\Pf. Recall that the condition $l>\dim_\C X$ implies that the product morphism is bounded in the $W^{l,2}$-norm. Furthermore, for each decorated tree $T\in \WW(k)$, by degree reasons, there are exactly $k-1$ white dots and at most $k-1$ black dots. In the definition of the operator $R_{k,T}$ we associated $K_t$ to black dots, and $h_t$ to white dots. By the previous Lemmas, they are uniformly bounded operators on $[0,\infty]$. Thus there exists a constant $D>0$ so that
\[ ||R_{k,T} ||_{l} \leq D^k.\]
Furthermore, the cardinality of $\WW(k)$ is bounded by $4^k$ because there are at most $2k-2$ edges (including the leaves of $T$) to be decorated by either white or black dots. Putting these together, we conclude that there exists a constant $C>0$ (which can be taken to be $4D$) which gives a uniform bound
\[ || R_k ||_{W^{l,2}} \leq C^k\]
in the operator norm.\ed

The above corollary implies that, for $B\in \Omega_X^{0,*}(\EEnd(\EE))_l$ such that $||B||_{W^{2,l}}$ is small enough, the infinite series
\[ R^B_k(\alpha_1,\cdots,\alpha_k):=\sum_{j_0\geq 0,\cdots,j_k\geq 0} R_{j_0+\cdots+j_k+k}(B^{j_0},\alpha_1,\cdots,\alpha_k,B^{j_k})\]
is absolutely convergent in the $C^0$-topology of $$C^0\big( [0,\infty], \Omega_X^{0,*}(\EEnd(\EE))_l\big),$$ the space of continuous maps from $[0,\infty]$ to the Hilbert space $\Omega_X^{0,*}(\EEnd(\EE))_l$. 

As we needed $C^\infty$ convergence in $t$-direction, we have to bound the derivatives of $R_k$ in the $t$-direction. For this purpose, it is essential to have $B\in \Omega_X^{0,*}(\EEnd(\EE))$, rather than its $l$-th Soblev completion.

\medskip
\begin{prop}~\label{pert-h-prop}
There exists a constant $r>0$, such that for $B\in \Omega_X^{0,*}(\EEnd(\EE))$ and $||B||_{W^{2,l}}<r$, the infinite series
\[ R^B_k(\alpha_1,\cdots,\alpha_k):=\sum_{j_0\geq 0,\cdots,j_k\geq 0} R_{j_0+\cdots+j_k+k}(B^{j_0},\alpha_1,\cdots,\alpha_k,B^{j_k})\]
is absolutely convergent in the $C^\infty$-topology of $C^\infty\big([0,\infty], \Omega_X^{0,*}(\EEnd(\EE))_l\big)$.
\end{prop}

\Pf. First we observe that each summand $R_{j_0+\cdots+j_k+k}(B^{j_0},\alpha_1,\cdots,\alpha_k,B^{j_k})$ is inside $C^\infty\big([0,\infty], \Omega_X^{0,*}(\EEnd(\EE))_l\big)$ because of the smoothness of $B$. This follows from the fact that for a smooth section $s$, then heat propagation $K_t s=e^{-t\Delta} s$ is smooth in $t$-direction on $[0,\infty]$.

Fix an integral $m>0$, we use a superscript $[m]$ to denote the $m$-th order derivative in the $t$-direction. Thus we have
\[ R_N^{[m]}= \sum_{T\in \WW(N)} R_{N,T}^{[m]}.\]
Recall that $T$ is decorated by $N-1$ white dots where we put $h_t$, and at most $N-1$ black dots where we put $K_t$. Thus when differentiating $m$ times, the number of different ways of differentiation is given by the number of nonnegative integer solutions of
\[ r_1+\cdots + r_{2N-2} =m,\]
which is ${2N-2+m-1 \choose m}$, which is a polynomial $f_m(N)$ in $N$ of degree $m$.

The $r$-th derivative of $K_t$ is $(-\Delta)^r e^{-t\Delta}$, and for $h_t$ this is $(-\db^*)(-\Delta)^{r-1}e^{-t\Delta}$. Thus, by Lemma~\ref{bound-2-lem} and Lemma~\ref{bound-3-lem}, near $t=0$, we have
\[ ||K_t^{[r]}||\leq C_r / t^r, \mbox{\;\; and \;\;} ||h_t^{[r]}||\leq C_r/ t^{r-1/2}.\]
However, there are at least $(N-1-m)$ white dots that are not differentiated, which implies  the estimate
\begin{align*}
 ||R_{N,T}^{[m]} ||_{W^{l,2}} &\leq f_m(N) \cdot \prod_{i=1}^{2N-2} \frac{C_{r_i}}{t^{r_i}}\cdots (\frac{C}{t^{1/2}})^{N-1-m} \cdot E^N\\
& \leq f_m(N) \cdot C_m \cdot (CE)^N \cdot t^{\frac{N-1-3m}{2}}.
\end{align*}
Here we used $E^N$ to denote the bound that results from the product map on vertices of $T$ and the remaining edges that are not differentiated which are uniformly bounded. Since the term $f_m(N)$, being a polynomial of $N$, is always bounded after dividing by, say $2^N$), we obtain that for $N>1+3m$ there exists a constant $D>0$, independent of $N$, such that
\[ ||R_{N,T}^{[m]}||_{W^{l,2}} \leq C_m \cdot D^N.\]
Finally, because the number $|\WW(N)|$ is bounded by $4^N$, we get
\[ ||R_N^{[m]}||_{W^{l,2}} \leq C_m\cdot (4D)^N.\]
This implies the $C^\infty$ convergence of the series in the lemma if we set $r=\frac{1}{4D}$.\ed

For the other component $S$ of $H^\dagger=R+Sdt$, the proof that the series
\[S^B_k(\alpha_1,\cdots,\alpha_k):=\sum_{j_0\geq 0,\cdots,j_k\geq 0} S_{j_0+\cdots+j_k+k}(B^{j_0},\alpha_1,\cdots,\alpha_k,B^{j_k})\]
is absolutely convergent in the $C^\infty$-topology is almost identical, except that due to the extra blue colored vertex where we put $-\db^*K_t$, we need to use different norms. Namely, on the source $\Om_X^{0,*}(\EEnd(\EE))$ we again use the $W^{l,2}$-norm, while on the target we use the $W^{l-1,2}$-norm. Thus we obtain the following

\medskip
\begin{cor}~\label{bound-h-cor}
There exists a constant $r>0$, such that for $B\in \Omega_X^{0,*}(\EEnd(\EE))$ and $||B||_{W^{2,l}}<r$, the infinite series
\[ (H^{\dagger})^B_k(\alpha_1,\cdots,\alpha_k):=\sum_{j_0\geq 0,\cdots,j_k\geq 0} H^\dagger_{j_0+\cdots+j_k+k}(B^{j_0},\alpha_1,\cdots,\alpha_k,B^{j_k})\]
is absolutely convergent in the $C^\infty$-topology of $\Omega^*_{[0,\infty]}\big(\Omega_X^{0,*}(\EEnd(\EE))_{l-1}\big)$, the space of smooth differential forms on $[0,\infty]$ with values in the Banach space $\Omega_X^{0,*}(\EEnd(\EE))_{l-1}$. The pull-back $H=f^*H^\dagger$ via $f(u)=\tan\frac{\pi u}{2}$ enjoys the same property by Lemma~\ref{exp-decay}.
\end{cor}

\medskip
\subsection{L-infinity case}~\label{tra-lie-sub} There is a functor $$\Lie: A_\infty \rightarrow L_\infty$$ from the category of $A_\infty$ algebras to $L_\infty$ algebras generalizing the usual functor sending an associative algebra to the underlying Lie algebra with the commutator bracket. For an $A_\infty$ algebra $A$, we denote by $A^{\Lie}$ the corresponding $L_\infty$ algebra. Explicitly, if the $k$-th structure map of $A$ is
\[ m_k: (A[1])^{\otimes k} \ra A[1]\]
of degree one, then the corresponding $k$-th bracket map $$ l_k: S^k (A[1]) \ra A[1] $$ is given by $$l_k(a_1,\cdots,a_k):= \sum_\sigma \epsilon_\sigma m_k(a_{\sigma(1)},\cdots,a_{\sigma(k)})$$ where the sum is over all permutations of $(1,\cdots,k)$, and $\epsilon_\sigma= \pm 1$ is the Koszul sign of permuting the tensors $a_1,\cdots, a_k$ to $a_{\sigma(1)},\cdots, a_{\sigma(k)}$. Similarly, on the level of morphisms, the symmetrization of an $A_\infty$ homomorphism induces an $L_\infty$ homomorphism on the corresponding $L_\infty$ algebras.

\medskip
\begin{rem}
The symmetrization of a homotopy defined using coalgebra structure~\cite{Markl} does not yield a homotopy on the corresponding cocommutative coalgebra, as was pointed out, for example in~\cite{Man}. This is the main reason we need to use Sullivan's definition of a homotopy since it only involves homomorphisms which can be symmetrized. On the other hand, one might ask the question why do we need to use $A_\infty$ structures to begin with, if we are ultimately dealing with $L_\infty$ structures. The reason for this is that we do not know how the tree formula works in the backward direction (the direction of $P$)
when performing homological transfer of $L_\infty$ algebras. Thus it is not clear how to deduce the corresponding norm property in this direction.
\end{rem}

By symmetrization, we obtain a $L_\infty$ algebra $\Ext^*(\EE,\EE)^{\Lie}$. We are mainly interested in its subspace $\Ext^*_0(\EE,\EE)$ consisting of traceless classes (see Subsection~\ref{loc-lie-subsec}). 

\medskip
\begin{lem}~\label{lie-sub-lem}
Let $A$ be a differential graded algebra with a direct sum decomposition
\[ A^{\Lie} \cong C\oplus A_0\]
as differential graded Lie algebras, with the Lie bracket of $C$ trivial. Let $(i,p,h)$ be a homotopy retraction of $A$ onto its cohomology. Then we have a direct sum decomposition
\[ H^*(A)^{\Lie} \cong H^*(C) \oplus H^*(A_0)\]
as $L_\infty$ algebras (with $H^*(A_0)$ endowed with the sub-$L_\infty$ algebra structure).
\end{lem}

\Pf. Observe that the symmetrization of the associative tree formula gives precisely the tree formula in the Lie case. Thus, in the Lie case, if one of the leaves of the tree contains an element $\alpha\in H^*(C)\stackrel{i}{\hookrightarrow} C$, at the nearest vertex to this leaf we must apply the bracket operator with $\alpha$, which gives zero.\qed

Applying this lemma to the case $A=\Omega_X^{0,*}(\EEnd(\EE))$ with the decomposition
\[ A\cong \Omega_X^{0,*}\oplus \Omega_X^{0,*}(\EEnd_0(\EE))\]
yields the following 

\medskip
\begin{cor}~\label{lie-cor}
The Lie structure on $\Ext^*(\EE,\EE)^{\Lie}$ restricts to the subspace $\Ext^*_0(\EE,\EE)$. By Lemma~\ref{norm-m-lem}, the $L_\infty$ algebra $\Ext^*_0(\EE,\EE)$ is a normed $L_\infty$ algebra in the sense of Definition~\ref{normed-alg-defi}.
\end{cor}

Since the decomposition $A\cong \Omega_X^{0,*}\oplus \Omega_X^{0,*}(\EEnd_0(\EE))$ is as differential graded Lie algebras, the projection map
\[ A \twoheadrightarrow \Omega_X^{0,*}(\EEnd_0(\EE))\]
is a morphism of differential graded Lie algebras. Since projection is bounded, we obtain the following

\medskip
\begin{cor}~\label{i-cor}
The composition
\[ \II: \Ext^*_0(\EE,\EE) \stackrel{I^{\Lie}}{\longrightarrow} \Omega_X^{0,*}(\EEnd(\EE))^{\Lie} \twoheadrightarrow \Omega_X^{0,*}(\EEnd_0(\EE))\]
is a $L_\infty$ quasi-isomorphism. By Lemma~\ref{norm-i-lem}, the homomorphism $\II$ is bounded in the sense of Definition~\ref{normed-hom-defi}.
\end{cor}

The backward $L_\infty$ homomorphism can be constructed in the same way. More precisely, by Lemma~\ref{lie-sub-lem}, the projection morphism
\[ \Ext^*(\EE,\EE)^{\Lie} \twoheadrightarrow \Ext_0^*(\EE,\EE)\]
is also an $L_\infty$ morphism. Hence we obtain

\medskip
\begin{cor}~\label{p-cor}
The composition
\[ \PP: \Omega_X^{0,*}(\EEnd_0(\EE)) \stackrel{P^{\Lie}}{\longrightarrow} \Ext^*(\EE,\EE)^{\Lie} \twoheadrightarrow \Ext^*_0(\EE,\EE)\]
is a $L_\infty$ quasi-isomorphism. Moreover, we have $\PP\circ \II=\id$. By Lemma~\ref{norm-p-lem}, the homomorphism $\PP$ is also bounded.
\end{cor}

\Pf. The identity $\PP\circ \II=\id$ easily follows from $P\circ I=\id$.\qed

\medskip
Using the splitting of the differential graded Lie structure, the composition
\[ \HH: \Omega_X^{0,*}(\EEnd_0(\EE))\hookrightarrow \Omega_X^{0,*}(\EEnd(\EE)) \stackrel{H}{\longrightarrow} \Omega_X^{0,*}(\EEnd(\EE))\ot\Omega_{[0,1]}^* \twoheadrightarrow \Omega_X^{0,*}(\EEnd_0(\EE))\ot \Om_{[0,1]}^*\]
defines a homotopy between $\id$ and $\II\circ \PP$. Corollary~\ref{bound-h-cor} then implies the desired boundedness property, which we summarize in the following

\medskip
\begin{prop}~\label{h-prop}
There exists a $L_\infty$ homomorphism
\[ \HH: \Omega_X^{0,*}(\EEnd_0(\EE)) \ra \Omega_X^{0,*}(\EEnd_0(\EE))\ot \Omega_{[0,1]}^*\]
such that $\HH|_{u=0}=\id$ and $\HH|_{u=1}=\II\circ \PP$. Moreover, There exists a constant $r>0$, such that for $B\in \Omega_X^{0,*}(\EEnd_0(\EE))$ and $||B||_{W^{2,l}}<r$, the infinite series
\[ \HH^B_k(\alpha_1,\cdots,\alpha_k):=\sum_{j\geq 0} \frac{1}{j!}\cdot \HH_{j+k}(B^j,\alpha_1,\cdots,\alpha_k)\]
is absolutely convergent in the $C^\infty$-topology of $\Omega^*_{[0,1]}\big(\Omega_X^{0,*}(\EEnd(\EE))_{l-1}\big)$, the space of smooth differential forms on $[0,1]$ with values in the Banach space $\Omega_X^{0,*}(\EEnd(\EE))_{l-1}$.
\end{prop}

\section{Proof of Proposition~\ref{kan-prop}}~\label{kan-sec}

Let $A$ and $B$ be two normed $L_\infty$ algebras. The purpose of this appendix is to prove the Kan property of the simplicial subset
\[ \MC^b_\bullet(C(A,B)) \subset \MC_\bullet(C(A,B))\]
defined in Subsection~\ref{simp-lie-subsec}. The proof is essentially the same as that of the Kan property of $\MC_\bullet(C(A,B))$ in~\cite[Section 4]{Getzler}, with additional work on the boundedness property. We follow the notation used in {\sl loc. cit.}.

Fix integers $n\geq 0$, and $0\leq i\leq n$. In~\cite{Dup}, Dupont introduced an explicit homotopy operator
\[ h_n^i: \Om_{\Delta^n}^* \ra \Om_{\Delta^n}^{*-1}\]
for the de Rham complex $\Om_{\Delta^n}^*$. That is,  there is an identity 
\[ \id = \epsilon_n^i +  d h_n^i +h_n^i d\]
where $d$ stands for the de Rham differential, and $\epsilon_n^i: \Om_{\Delta^n}^* \ra \R\hookrightarrow \Om_{\Delta^n}^*$ is the evaluation map at the $i$-th vertex of $\Delta^n$ followed by the inclusion of constant functions. Dupont's homotopy induces a homotopy on the tensor product $C(A,B)\ot \Om_{\Delta^n}^*$ space:
\begin{equation}~\label{dupont-eq}
 \id = \epsilon_n^i +(d+\delta)h_n^i +h_n^i(d+\delta)
\end{equation}
where $\delta=L_1$ is the differential on $C(A,B)$. In the last equation, $\epsilon_n^i$ really stands for $$\id\ot\epsilon_n^i: C(A,B)\ot \Om_{\Delta^n}^* \ra C(A,B)\ot \Om_{\Delta^n}^*,$$and similarly for $h_n^i$. But we choose to use this abbreviation as no confusion can arise this way.

Recall the definition of the operator $h_n^i$ as the composition
\[ h_n^i:= \pi^n_* (\phi^n_i)^*,\]
using the diagram
\[ \begin{CD} \Delta^n @< \phi^n_i << [0,1]\times\Delta^n @>\pi^n>> \Delta^n.\end{CD}\]
In the diagram, the map $\phi_i^n$ is defined by 
\[ (u,t) \mapsto ut+(1-u)e^n_i.\]
To avoid long notations, we set
\[ \MC_n:=\MC\big( C(A,B\ot\Om_{\Delta^n}^*)\big), \mbox{ and } \mc_n:= \left\{ (d+\delta)\beta\mid \beta\in  C(A,B\ot\Om_{\Delta^n}^*)_0\right\},\]
and the bounded version to be
\[ \MC_n^b:=\MC_n\cap C^b(A,B\ot \Om_{\Delta^n}^*), \mbox{ and } \mc^b_n:= \mc_n \cap C^b(A,B\ot\Om_{\Delta^n}^*).\]
In~\cite{Getzler} it was shown that there is a bijection $\MC_n\cong \MC_0 \times \mc_n$ defined by the map
\begin{equation}~\label{bijection-eq}
 \alpha \mapsto (\epsilon_n^i\alpha, (d+\delta)h_n^i\alpha).\end{equation}
We first prove this bijection restricts the bounded subsets.

\medskip
\begin{lem}~\label{h-lem}
Assume that $\psi\in C^b(A,B\ot \Om_{\Delta^n}^*)$, then we have $h_n^i\psi\in C^b(A,B\ot \Om_{\Delta^n}^*)$.
\end{lem}

\Pf. Let us write
\[ \psi_k=\sum_J (\psi_k)_J dt_J\]
according to the decomposition of differential forms. By definition $h_n^i= (\pi^n)_* (\phi^n_i)^*$, thus we have 
\[h_n^i[(\psi_k)_Jdt_J]= \sum_{j\in J} \pm \int_0^1 \psi_k(ut+(1-u)e_i^n)(t_j-\delta_{ij}) u^{|J|-1} du dt_{J/j}.\]
Differentiate in the $t$-direction and using the fact that
\[ |t_j-\delta_{ij}|<1, \mbox{\;\; and \;\; } |u|<1\]
proves that
\[ ||h_n^i \psi_k ||_m \leq 2 ||\psi_k||_m.\]
Hence if $\psi$ is bounded, so is $h_n^i\psi$.\ed

\medskip
\begin{lem}~\label{d-delta-lem}
Assume that $\psi\in C^b(A,B\ot \Om_{\Delta^n}^*)$, then we have $d\psi\in C^b(A,B\ot \Om_{\Delta^n}^*)$, and that $\delta\psi=L_1\psi\in C^b(A,B\ot \Om_{\Delta^n}^*)$.
\end{lem}

\Pf. Indeed, we have
\[ ||d\psi_k||_m \leq || \psi_k ||_{m+1} \leq D_{m+1} C^k.\]
Next we verify that $L_1(\psi)$ is bounded. Let $C>0$ be a uniform constant bounding the $L_\infty$ structures of $A$, $B$, and $\psi$. Then we have
\begin{align*}
||[L_1(\psi)]_k ||_{C^m} &= ||\sum_{i=1}^k \sum_{\sigma\in\Sh (i,k-i)} \pm\psi_{k-i+1}\circ (l_i^A\ot \id^{k-i}) \circ \sigma + l_1^B \circ \psi_k ||_{C^m}\\
&\leq \sum_{i=1}^k \frac{k!}{i!(k-i)!} D_m\cdot (k-i+1)!\cdot C^{k-i+1}\cdot i!\cdot C^i + D_m \cdot k! \cdot C^{k+1}\\
&\leq D_m\cdot k!\cdot C^{k+1} [\sum_{i=1}^k (k-i+1) + 1]\\
&\leq D_m\cdot k!\cdot C^k \cdot \big( C\cdot \frac{k^2+k+2}{2}\big)\\
& \leq D_m\cdot k! \cdot E^k.
\end{align*}
The last equality follows from that polynomial growth is always bounded by exponential growth. Clearly the constant $E$ is independent of $m$ and $k$. \ed

\begin{lem}~\label{bijection-lem}
The bijection defined in~\ref{bijection-eq} restricts to a bijection
\[ \MC_n^b\cong \MC^b_0\times \mc^b_n\]
between the bounded subsets. 
\end{lem}

\Pf. Let $\alpha$ be in $\MC_n^b$. By the previous two Lemmas, we have that
\[ (d+\delta)h_n^i\alpha\in \mc_n^b.\]
Since the map $\epsilon_n^i$ is the evaluation map at the $i$-th vertex, it certainly previous the boundedness property. Hence we have shown that
\[ (\epsilon_n^i \alpha, (d+\delta)h_n^i\alpha) \in \MC_0\times \mc_n^b.\]
In the reverse direction, it was shown in~\cite{Getzler} that given $(\mu,\nu)\in \MC^b_0\times \mc^b_n$, the associated Maurer-Cartan element $\alpha\in \MC_n$ is defined inductively by 
\begin{align*}
\alpha_0&:= \mu+\nu,\\
\alpha_{k+1}&:= \alpha_0 - \sum_{ j\geq 2} \frac{1}{j!} h_n^i L_j(\alpha_k,\cdots,\alpha_k),\\
\alpha&:= \lim_{k\ra \infty} \alpha_k.
\end{align*}
For an element 
\[ \alpha\in C(A,B)\ot \Om_{\Delta^n}^* = \Hom( \overline{S}^c(A[1]),B)\ot \Om_{\Delta^n}^*,\]
denote by $[\alpha]_N$ the component in $\Hom( S^N(A[1]),B)\ot \Om_{\Delta^n}^*$ where $N\geq 1$. By induction on $N$, one can show that $N$-th components of $a_k$'s stabilize as $k$ goes to infinity:
\[ [\alpha_k]_N=[\alpha_{k+1}]_N=[\alpha_{k+2}]_N=\cdots.\]
Hence the definition $\alpha = \lim_{k\ra \infty} \alpha_k$ makes sense.

We need to prove the boundedness of $\alpha$. Recall the term $L_j(\alpha_k,\cdots,\alpha_k)$ in the definition of $\alpha$ is given by the composition
\[ \overline{S}^c(A[1]) \xrightarrow{\substack{\text{iterated} \\\text{ coproduct}}} \overline{S}^c(A[1])^{\ot j} \stackrel{\alpha_k^{\ot j}}{\longrightarrow} (B\ot \Om_{\Delta^n}^*)^{\ot j} \xrightarrow{l_j^{B\ot \Om_{\Delta^n}^*}} B\ot \Om_{\Delta^n}^*.\]
Here $l_j^{B\ot \Om_{\Delta^n}^*}$ is the structure map of the $L_\infty$ algebra $B$, extended $\Om_{\Delta^n}^*$-linearly to the tensor product $B\ot \Om_{\Delta^n}^*$. Using the formula for the iterated coproduct on symmetric coalgebras, we get the following inductive formula
\[ [\alpha]_N= [\alpha_0]_N-h_n^i\sum_{j=2}^N \frac{1}{j!} \sum_{\substack{s_1,\cdots,s_j\geq 1\\ s_1+\cdots+s_j=N}} \sum_{\sigma\in \Sh(s_1,\cdots,s_j)} l^{B\ot\Om_{\Delta^n}^*}_j\circ ([\alpha]_{s_1}\ot\cdots\ot[\alpha]_{s_j})\circ \sigma.\]
This recurrence relation can be solved by summing over decorated trees. 
Indeed, Let $T$ be a rooted tree with $N$ leaves, and with internal vertices of at least valency $3$. We put decorations on the set of vertices of $T$ by coloring it black or white, such that the root vertex is not colored, a leaf vertex (the unique vertex of a leaf) may be colored or not, and  internal vertices  must be colored. Furthermore, a decoration must satisfy the following condition: for any leaf vertex , the decoration of the unique path from it to the root vertex is either one of the following two types:
\[\begin{tikzpicture}
\node[circle,fill=white,text=black] (z) at (-10.7, 0) {leaf};
\node[circle,fill=white,text=black] (z) at (-3.3, 0) {root};
\node[circle,fill=black] (a) at (-10, 0) {};
\node[circle,draw] (b) at (-9, 0) {};
\node[circle,draw] (c) at (-8, 0) {};
\node[circle,draw] (d) at (-7, 0) {};
\node (e) at (-6, 0) {$\cdots$};
\node[circle,draw] (f) at (-5, 0) {};
\node (g) at (-4, 0) {};

\draw (a) -- (b);
\draw (b) -- (c);
\draw (c) -- (d);
\draw (d) -- (e);
\draw (e) -- (f);
\draw (f) -- (g);
\end{tikzpicture}\]

\[\begin{tikzpicture}
\node[circle,fill=white,text=black] (z) at (-10.7, 0) {leaf};
\node[circle,fill=white,text=black] (z) at (-3.3, 0) {root};
\node (a) at (-10, 0) {};
\node[circle,fill=black] (b) at (-9, 0) {};
\node[circle,draw] (c) at (-8, 0) {};
\node[circle,draw] (d) at (-7, 0) {};
\node (e) at (-6, 0) {$\cdots$};
\node[circle,draw] (f) at (-5, 0) {};
\node (g) at (-4, 0) {};

\draw (a) -- (b);
\draw (b) -- (c);
\draw (c) -- (d);
\draw (d) -- (e);
\draw (e) -- (f);
\draw (f) -- (g);
\end{tikzpicture}\]
We denote by $\OO^\dec(N)$ the set of isomorphism classes of such decorated trees. For each decorated tree $(T,\sigma)\in \OO^\dec(N)$ we define an operator 
\[ \rho_{(T,\sigma)}: S^k(A[1]) \ra B\ot \Om_{\Delta^n}^*\]
in the following way. First we choose a planar realizations $\widetilde{T}$ of $T$. For a leaf vertex of $\widetilde{T}$ which is decorated by a black dot, then we assign to it the operator $[\alpha_0]_1$; otherwise it is just the $\id$ operator. For a vertex $v$ of valency $j+1$ which is decorated by a white dot, we assign the operator $l^{B\ot\Om_{\Delta^n}^*}_j$. For a vertex $v$ of valency $s+1$ which is decorated by a black dot, we put the operator $[\alpha_0]_{s}$. Finally, for an edge which is the unique outgoing edge of a white dot, we assign the operator $h_n^i$. Then we define $\rho_{(\widetilde{T},\sigma)}$ to be the operadic composition by the associated morphisms, and set
\[ \rho_{(T,\sigma)}:=  \rho_{(\widetilde{T},\sigma)}\circ\Sym\]
to be the composition of $\rho_{(\widetilde{T},\sigma)}$ with  the symmetrization map $\Sym: S^N(A[1]) \ra T^N(A[1])$. The map $\rho_{(T,\sigma)}$ is independent of the choice of $\widetilde{T}$ because the symmetric group acts transitively on the set of planar realizations of $T$.

With this notations, one can prove that
\begin{equation}~\label{nonind-eq}
 [\alpha]_N= \sum_{(T,\sigma)\in \OO^\dec(N)} (-1)^\epsilon \frac{1}{|\Aut T|} \rho_{(T,\sigma)}
 \end{equation}
where $\epsilon$ is the number of white dots in $T$. The proof is similar to that of~\cite[Proposition 4.2]{FM} which boils down to the combinatorics of trees~\cite{Fio}.

Using the non-inductive formula~\ref{nonind-eq}, the boundedness of $\alpha$ follows from the boundedness of the operators $\alpha_0$, $l^{B\ot\Om_{\Delta^n}^*}$, $h_n^i$, and the fact that the number of decorated trees is bounded by an exponential function in $N$.
\ed

Next we prove the main result of the section, Proposition~\ref{kan-prop}. We have a simplicial $L_\infty$ algebra defined by
\[ C(A,B\ot \Om_{\Delta^\bullet}^*).\]
As a simplicial set, it is always a Kan complex because any simplicial group is. By definition there are inclusions
\[ \Hom(A,B)_\bullet=\MC^b_\bullet(C(A,B))\subset \MC_\bullet(C(A,B)) \subset C(A,B\ot \Om_{\Delta^\bullet}^*).\]
Let $$\beta:\Lambda_n^j\ra \Hom(A,B)_\bullet$$ be a horn. The composition
\[ \Lambda_n^j \ra \Hom(A,B)_\bullet\subset C(A,B\ot \Om_{\Delta^\bullet}^*)\]
gives a horn in $C(A,B\ot \Om_{\Delta^\bullet}^*)$. Since the latter is Kan, there exists a filling
\[ \widetilde{\beta}: \Delta^n \ra C(A,B\ot \Om_{\Delta^\bullet}^*).\]
We claim that $\widetilde{\beta}$ can be chosen so that 
\[\widetilde{\beta}\in C^b(A,B\ot\Om_{\Delta^n}^*).\]
To prove the claim, we follow the inductive construction of $\widetilde{\beta}$ in~\cite[Lemma 8.2.8]{Weibel}. Recall that to give a horn $\Lambda_n^j \ra \Hom(A,B)_\bullet\subset C(A,B\ot \Om_{\Delta^\bullet}^*)$ is equivalent to give elements $$x_0,\cdots,x_{j-1}, x_{j+1},\cdots,x_{n} \in C(A,B\ot \Om_{\Delta^{n-1}}^*)$$ such that $\partial_s^* x_t=\partial_{t-1}^* x_s$ for all $s<t$ ($s$ and $t$ not equal to $j$). Here $\partial_s^*$, for each $0\leq s\leq n$, is the pull-back morphism associated to the simplicial map
\[ \partial_s: \Delta^{n-1}\ra \Delta^n.\]
A filling of this horn can be inductively constructed as follows. To begin, set $g_{-1}=0$ and suppose that $g_{r-1}$ is given. If $r=j$ we set $g_r:=g_{r-1}$. If $r\neq j$, then set
\begin{equation}~\label{induction-eq}
 g_r:= g_{r-1} -\sigma^*_j (\partial_r^*g_{r-1}- x_r).\end{equation}
Here $\sigma_j^*$ is defined by the pull back morphism of the simplicial map
\[  \sigma_j: \Delta^{n}\ra \Delta^{n-1}.\]
Then one can prove that $\widetilde{\beta}:=g_n\in C(A,B\ot \Om_{\Delta^n}^*)$ is a filling of $\beta$. 

If all the $x_r$'s are in $C^b(A,B\ot\Om_{\Delta^{n-1}}^*)$, then the lifting $\widetilde{\beta}$, as defined inductively above, lies inside $C^b(A,B\ot\Om_{\Delta^{n}}^*)$ since pull-backs preserves the boundedness property, which proves the claim.

Proposition~\ref{kan-prop} now follows from Lemma~\ref{bijection-lem}. Indeed we choose a lifting of the horn 
\[ \beta: \Lambda_n^j \ra \Hom(A,B)_\bullet\]
given by the unique element $\alpha\in \MC_n^b$ such that 
\[ (\epsilon_n^i \alpha, (d+\delta)h_n^i \alpha)=(\epsilon_n^i\widetilde{\beta},(d+\delta)h_n^i \widetilde{\beta}).\]
The check that $\alpha$ is a lifting of $\beta$ is the same as in~\cite{Getzler}.

\end{document}